\documentclass{amsart}
\usepackage{amsmath}
\usepackage{amsthm}
\usepackage{amsfonts}
\usepackage{amssymb}
\usepackage{graphicx}
\usepackage{mathrsfs}

\usepackage{a4wide,color,eucal,enumerate,mathrsfs}
\usepackage[normalem]{ulem}
\usepackage{amsmath,amssymb,epsfig,bbm}
\usepackage{pdfsync}
\numberwithin{equation}{section}
\usepackage[pdfborder={0 0 0}]{hyperref}

\newcommand{\B}{\mathbb{B}}

\newcommand{\D}{\mathbb{D}}

\newcommand{\M}{\mathbb{M}}
\newcommand{\N}{\mathbb{N}}
\newcommand{\Q}{\mathbb{Q}}
\newcommand{\R}{\mathbb{R}}

\newcommand{\mm}{{\mbox{\boldmath$m$}}}
\newcommand{\nn}{{\mbox{\boldmath$n$}}}

\newcommand{\sfd}{{\sf d}}

\newcommand{\sfr}{{\sf r}}

\newcommand{\restr}[1]{\lower3pt\hbox{$|_{#1}$}}

\newcommand{\la}{{\langle}}                  
\newcommand{\ra}{{\rangle}}

\newcommand{\eps}{\varepsilon}  
\newcommand{\nchi}{{\raise.3ex\hbox{$\chi$}}}

\newtheorem{theorem}{Theorem}[section]

\newtheorem{corollary}[theorem]{Corollary}
\newtheorem{lemma}[theorem]{Lemma}
\newtheorem{proposition}[theorem]{Proposition}
\newtheorem{definition}[theorem]{Definition}

\newtheorem{remark}[theorem]{Remark}

\newcommand{\dist}{\mathrm{dist}}
\newcommand{\LIP}{\mathrm{LIP}}
\newcommand{\Lip}{\mathrm{Lip}}
\newcommand{\lip}{\mathrm{lip}}
\newcommand{\dir}{\mathrm{dir}}

\newcommand{\dive}{\mathrm{div}}

\newcommand{\diam}{\mathrm{diam}}

\newcommand{\inv}{^{-1}}

\newcommand{\fr}{\hfill$\blacksquare$}   

\renewcommand{\mm}{\mathfrak m}

\newcommand{\Cat}[1]{{\rm CAT}(#1)}

\newcommand{\lims}{\varlimsup}
\newcommand{\limi}{\varliminf}
\renewcommand{\limsup}{\varlimsup}
\renewcommand{\liminf}{\varliminf}
\renewcommand{\d}{{\rm d}}
\newcommand{\Ba}{{\sf Bar}}
\newcommand{\Geo}{{\sf Geo}}
\newcommand{\X}{{\rm X}}
\newcommand{\Y}{{\rm Y}}
\newcommand{\T}{{\rm T}}
\newcommand{\sn}{{\rm sn}}
\newcommand{\cn}{{\rm cn}}
\newcommand{\G}{{\sf G}}
\renewcommand{\div}{{\rm div}}
\renewcommand{\nn}{{\mathfrak n}}
\newcommand{\Dertwotwo}{{\rm Der}^{2,2}(\Y;\mu)}

\title{Infinitesimal Hilbertianity of locally $\Cat\kappa$-spaces}

\author{Simone Di Marino}
\address{Istituto Nazionale di Alta Matematica, Unit\`a INdAM SNS Pisa, Piazza dei Cavalieri 7, 56126 Pisa}
\email{simone.dimarino@altamatematica.it}

\author{Nicola Gigli}
\address{SISSA, Via Bonomea 265, 34136 Trieste}
\email{ngigli@sissa.it}

\author{Enrico Pasqualetto}
\address{SISSA, Via Bonomea 265, 34136 Trieste, and \newline
Department of Mathematics and Statistics,
P.O.\ Box 35 (MaD), FI-40014 University of Jyv\"{a}skyl\"{a}}
\email{enpasqua@jyu.fi}

\author{Elefterios Soultanis}
\address{SISSA, Via Bonomea 265, 34136 Trieste, and\newline
University of Fribourg, Chemin du Musee 23, CH-1700 Fribourg}
\email{elefterios.soultanis@gmail.com}

\begin{document}
\maketitle

\begin{abstract}
	
We show that, given a metric space $(\Y,\sfd)$ of curvature bounded from above in the sense of Alexandrov, and a positive Radon measure $\mu$ on $\Y$ giving finite mass to bounded sets, the resulting metric measure space $(\Y,\sfd,\mu)$ is infinitesimally Hilbertian, i.e. the Sobolev space $W^{1,2}(\Y,\sfd,\mu)$ is a Hilbert space. 

The result is obtained by constructing an isometric embedding of the `abstract and analytical' space of derivations into the   `concrete and geometrical'  bundle whose fibre at $x\in\Y$ is the tangent cone at $x$ of $\Y$. The conclusion then follows from the fact that for every $x\in\Y$ such a cone is a $\Cat0$ space and, as such, has a Hilbert-like structure.
\end{abstract}

\tableofcontents

\section{Introduction}

A metric space $(\Y,\sfd)$ is said to be a $\Cat\kappa$ space if, roughly said, it is geodesic and geodesic triangles are `thinner' than triangles in the model space $\M_\kappa$ of constant sectional curvature $=\kappa$. Typical examples of $\Cat\kappa$ spaces are simply connected Riemannian manifolds with sectional curvature $\leq\kappa$ and their Gromov-Hausdorff limits. Despite the absence of any a priori smooth structure, $\Cat\kappa$ spaces are quite regular and carry a solid calculus  resembling that on manifolds with curvature $\leq\kappa$.  We refer to \cite{Gromov07}, \cite{Bac14}, \cite{BH99}, \cite{BBI01}, \cite{AKP}  for overviews on the topic and a more detailed bibliography.

In the particular case $\kappa=0$ the $\Cat0$ condition reads as follows: for any points $x_0,x_1\in\Y$ and any geodesic $\gamma:[0,1]\to\Y$ connecting them it holds that
\begin{equation}
\label{eq:cat0}
\sfd^2(\gamma_t,y)\leq(1-t)\sfd^2(x_0,y)+t\sfd^2(x_1,y)-t(1-t)\sfd^2(x_0,x_1)\qquad\forall y\in\Y,\,t\in[0,1].
\end{equation}
This can be regarded as a parallelogram inequality and from this point of view it is perhaps not surprising that several aspects of $\Cat0$ spaces strongly resemble properties of Hilbert spaces; this perspective is emphasised e.g.\ in \cite{Bac14}. For instance, from \eqref{eq:cat0} it directly follows that if a normed vector space is a $\Cat0$ space, then the norm comes from a scalar product. Equivalently,
\begin{equation}
\label{eq:normcat}
\begin{array}{l}
\text{if a normed vector space isometrically embeds in a $\Cat0$ space,}\\
\text{then the norm comes from a scalar product.}
\end{array}
\end{equation}
Given that $\Cat0$ spaces naturally arise as tangent cones to generic $\Cat\kappa$ spaces, these analogies with Hilbert structures appear also at small scales on $\Cat\kappa$ spaces.

\bigskip

A metric measure space $(\Y,\sfd,\mu)$ is called \emph{infinitesimally Hilbertian} provided the Sobolev space $W^{1,2}(\Y,\sfd,\mu)$ is Hilbert (see \cite{Cheeger00} and then also \cite{Shanmugalingam00}, \cite{AmbrosioGigliSavare11} for the definition of Sobolev spaces in this context). The concept of infinitesimal Hilbertianity, introduced in \cite{Gigli12}, aims at detecting Hilbert structures at small scales in the non-smooth setting. The motivating example in the smooth category is the following: if $\Y$ is a smooth Finsler manifold and $\mu$ is a smooth measure on it (i.e. with smooth density when seen in charts), then the $W^{1,2}$-norm can be written as
\begin{equation}
\label{eq:w12norm}
\|f\|_{W^{1,2}}^2=\int|f|^2(x)+\|\d f(x)\|_x^2\,\d\mu(x).
\end{equation}
Since $f\mapsto \int |f|^2\,\d\mu$ always satisfies the parallelogram identity, we see that $f\mapsto\|f\|_{W^{1,2}}^2$ has the same property if and only if $f\mapsto \int \|\d f(x)\|_x^2\,\d\mu(x)$ satisfies the parallelogram identity. With a little bit of work it is possible to check that this is the case if and only if $\|\cdot\|_x^2$ satisfies the parallelogram identity for every $x$, i.e.\ if and only if $\Y$ is in fact a Riemannian manifold.

In the smooth category one could run the above consideration also with smooth functions, rather than with Sobolev ones, but this is obviously not possible on a metric measure space. In this direction let us emphasise that in the non-smooth environment it is crucial to work with Sobolev functions rather than, say, with Lipschitz ones. To see why, recall that the local Lipschitz constant $\lip f:\Y\to[0,\infty]$ of a function $f:\Y\to\R$ is defined as
\[
\lip f(x):=\lims_{y\to x}\frac{|f(y)-f(x)|}{\sfd(x,y)}\quad\text{if $x$ is not isolated, 0 otherwise}
\]
and consider the following example. Let $(\Y,\sfd)$ be the Euclidean space $(\R^d,\sfd_{\rm Eucl})$ and $\mu$ be a positive Radon measure. Then:
\begin{itemize}
\item[a)] The map
\[
\LIP_c(\R^d)\ni f\quad\mapsto\quad\int (|f|^2+\lip^2 f)\,\d\mu
\]
is not a quadratic form, in general.
\item[b)] The map
\[
W^{1,2}(\R^d,\sfd_{\rm Eucl},\mu)\ni f\quad\mapsto\quad \|f\|^2_{W^{1,2}}
\]
is a quadratic form, i.e.\ $(\R^d,\sfd_{\rm Eucl},\mu)$ is infinitesimally Hilbertian.
\end{itemize}
To see why (a) holds simply consider $\mu$ to be a Dirac delta at a point $o$ and $f,g\in \LIP_c(\R^d)$ generic functions not differentiable at $o$: for these the parallelogram identity for $f\mapsto\lip^2 f(o)$ typically fails. Intuitively, this is due to the fact that, if $f$ and $g$ are not differentiable at $o$, they are not (close to being) linear in the vicinity of $o$ and thus their local Lipschitz constants fail to capture the Hilbert structure of the cotangent space $\T^*_o\R^d$ at $o$.

The statement in (b) is non-trivial and is one of the results proved in \cite{GP16-2}. The crucial aspect of the proof is the possibility of approximating Sobolev functions with $C^1$ functions: these are by nature differentiable everywhere, and thus also $\mu$-a.e., and hence are suitable to identify the Hilbertian structure of the cotangent spaces.

Hence the idea behind the notion of infinitesimal Hilbertianity is to exploit the fact that `by nature' Sobolev functions are a.e.\ differentiable in some sense, regardless of the regularity of the metric and of the measure in consideration (for instance, if $\mu$ is a Dirac delta as above, it turns out that Sobolev functions have 0 differential, so that the claim (b) is trivially true in this case). This makes them  suitable for detecting Hilbert structures at an infinitesimal scale.  Let us emphasise that even though this is an analytic notion, it is strictly related to -- and its introduction has been motivated by -- the study of geometric properties of metric measure spaces, in particular those satisfying a curvature-dimension bound in the sense of  Lott-Sturm-Villani. An example of this link is the validity of the non-smooth splitting theorem \cite{Gigli13,Gigli13over}, which states that under the appropriate geometric rigidity given by a LSV condition the weak and `differential' notion of infinitesimal Hilbertianity implies the validity of a kind of Pythagora's theorem for the `integrated' object $\sfd$.

\bigskip

These considerations about Sobolev functions, together with the fact that tangents of $\Cat\kappa$-spaces are $\Cat 0$-spaces and thus exhibit behaviour akin to Hilbert spaces, might lead one to suspect that a $\Cat\kappa$-space equipped with any measure is infinitesimally Hilbertian.

This is indeed the case and is the main result of this manuscript:
 \begin{theorem}[Universal infinitesimal Hilbertianity of local $\Cat\kappa$ spaces]\label{main}
	Let $\kappa\in\R$, $(\Y,\sfd)$ be a local $\Cat \kappa$-space and $\mu$ a non-negative and non-zero Radon measure on $\Y$ giving finite mass to bounded sets.
	
	Then $(\Y,\sfd,\mu)$ is infinitesimally Hilbertian, i.e.\ the Sobolev space $W^{1,2}(\Y,\sfd,\mu)$ is a Hilbert space. 
\end{theorem}
Let us collect some comments:
\begin{itemize}
\item[i)] Sobolev functions on metric measure spaces are typically studied either on generic mm-spaces, mostly for foundational purposes, or on spaces which are either \emph{doubling}, support a \emph{Poincar\'e inequality}, or have \emph{Ricci curvature bounded from below}. In these contexts, Sobolev spaces constitute a key ingredient for the development of a non-smooth calculus (see \cite{Cheeger00}, \cite{Bjorn-Bjorn11}, \cite{Heinonen07}, \cite{HKST15}, \cite{AmbrosioGigliSavare11} and the references therein). All these conditions are in strong contrast to the \emph{upper sectional curvature bound} encoded by the $\Cat\kappa$ notion as they all more-or-less point to a lower (Ricci) curvature bound. 

In this direction it is worth mentioning that $\Cat\kappa$ spaces do not carry any natural reference measure (unlike, for instance, finite-dimensional Alexandrov spaces with curvature bounded from below) and perhaps for this reason they have been investigated mostly as metric spaces, rather than as metric measure spaces.

To the best of our knowledge, this manuscript contains the first result about the structure of Sobolev functions on  $\Cat\kappa$ spaces.

\item[ii)] A particular case of Theorem \ref{main} has been obtained in the recent  paper \cite{KK18} by Kapovitch and Ketterer. There the authors  consider a metric measure space $(\X,\sfd,\mm)$ which is a ${\sf CD}(K,N)$ space in the sense of Lott-Sturm-Villani (\cite{Lott-Villani09}, \cite{Sturm06I,Sturm06II}) when seen as a metric measure space and a $\Cat\kappa$ space when seen as metric space. Among other things, they prove  that $(\X,\sfd,\mm)$ is infinitesimally Hilbertian, thus giving another instance of the fact that a $\Cat\kappa$ condition forces $W^{1,2}$ to be Hilbert.  Their proof is based on the strong rigidity which comes from having both a `lower Ricci' and an `upper sectional' curvature bound (in fact the study of such rigidity, and of the regularity it enforces, is their main goal) and cannot be adapted to our case.

\item[iii)] We have mentioned that, in \cite{GP16-2}, to prove the result stated in (b) above the use of $C^1$ functions is crucial. Something similar happens here, where  we make extensive use of the fact that on $\Cat\kappa$ spaces there are many semiconvex Lipschitz functions (e.g.\ distance functions) and they have a well-defined notion of differential at every point; see Subsection \ref{se:diff}.

\item[iv)] This manuscript is part of a broader program aiming at stating and proving the Bochner-Eells-Sampson inequality
\begin{equation}
\label{eq:BES}
\Delta\frac{|\d u|^2_{\sf HS}}2\geq\la\d u,\Delta\d u\ra_{\sf HS}+K|\d u|_{\sf HS}^2
\end{equation}
for maps from a ${\sf RCD}(K,N)$ space $(\X,\sfd_\X,\mm_\X)$ to a $\Cat0$ space $(\Y,\sfd)$. Notice that inequality (\ref{eq:BES})  would immediately imply Lipschitz regularity of harmonic maps, by well known elliptic regularity theory in the non-smooth setting.

The role of this manuscript, to be used in conjunction with \cite{GPS18}, is to ensure that $L^2(\T^*\Y;u_*(|\d u|^2\mm_\X))$ is a Hilbert module, so that the same holds for the tensor product $L^2(\T^*\X;\mm_\X)\otimes \big(u^*L^2(\T^*\Y;u_*(|\d u|^2\mm_\X))\big)$ and thus the `pointwise Hilbert-Schmidt norm' appearing in \eqref{eq:BES} makes sense. We refer to \cite{GPS18} and \cite{GT18} for more details on this.

\item[v)] $\Cat\kappa$ spaces are not necessarily separable (for instance, the $\Cat0$ space obtained by glueing uncountably many copies of $[0,1]$ at 0 is not separable), as opposed to finite-dimensional spaces with curvature bounded from \emph{below}. For this reason separability is not an assumption in Theorem \ref{main}. Still, given that Sobolev spaces on metric measure spaces are typically studied in a separable environment, we first prove our main result for separable spaces and postpone the technical details needed to handle the general case until the final section.
\end{itemize}

\bigskip

Let us briefly describe the proof of Theorem \ref{main}. The basic intuition  is given by \eqref{eq:normcat} and the fact that the tangent cone of a local $\Cat\kappa$ space is a $\Cat0$ space. More precisely, we consider:
\begin{itemize}
\item[(1)] The space $\Dertwotwo$ of \emph{derivations} (with divergence),
as introduced by the first author in \cite{DiM14a,DiMarino14}
(see Section \ref{se:mms}). These are in duality with Sobolev functions.
\item[(2)] The collection $L^2(\T_G\Y;\mu)$ of `$L^2(\mu)$ Borel sections of the bundle $\T_G\Y$ on $\Y$ whose fibre at $x$ is the tangent cone $\T_x\Y$' (see Section \ref{se:gtb}). 
\end{itemize}
In Theorem \ref{thm:main} and Corollary \ref{cor:main}, we construct an isometric embedding $$\mathscr F: \Dertwotwo\hookrightarrow L^2(\T_G\Y;\mu)$$ which respects distances fibrewise. From this fact, the arguments behind \eqref{eq:normcat} and the aforementioned duality between derivations and Sobolev functions easily imply the main Theorem \ref{main}.

To construct the embedding $\mathscr F$, recall that a derivation $b\in \Dertwotwo$ gives rise to a \emph{normal 1-current} $T_b$ in the sense of Ambrosio-Kirchheim \cite{AmbrosioKirchheim00} (Lemma \ref{le:dercurr}). Using Paolini--Stepanov's version \cite{PaolStep12,PaolStep13} of Smirnov's \emph{superposition principle} (see Theorem \ref{thm:sup2}) we express the 1-current $T_b$ as a superposition $\int[\![\gamma]\!]\,\d\pi_{T_b}(\gamma)$, where $\pi_{T_b}$ is a finite measure on the space of absolutely continuous curves and $[\![\gamma]\!]$ is the current induced by $\gamma$.

Inspired by \cite{Lyt04}, we see that if $\gamma$ is an absolutely continuous curve then the \emph{right and left derivatives} $\dot\gamma_t^+$ and $\dot\gamma^-_t$ exist as elements of $\T_{\gamma_t}\Y$, and satisfy $\dot\gamma_t^+\oplus\dot\gamma_t^-=0$, for almost every $t\in [0,1]$ (Proposition \ref{prop:velo}, Remark \ref{rmk:velo} and Lemma \ref{lem:velo}).

Given the measures $(\pi_{T_b}\times{\mathcal L}^1\restr{[0,1]})_x$, obtained by disintegrating $\pi_{T_b}\times{\mathcal L}^1\restr{[0,1]}$ with respect to the evaluation map $(\gamma,t)\mapsto \gamma_t$, we consider their push-forward by the `right-derivative' map (cf. Proposition \ref{prop:rightder}), thus obtaining measures $\nn_x$ supported in $\T_x\Y$.

The Borel section $\mathscr F(b)$ is defined to be, at almost every $x\in \Y$, the barycenter of $\nn_x$. The barycenter lies in the tangent cone $\T_x\Y$. By a rigidity property of barycenters (Lemma \ref{lem:rigidity}), and convexity properties of tangent cones, the measure $\mathfrak{n}_x$ is concentrated on a half-line for almost every $x\in \Y$.

\bigskip Theorem \ref{main2} below is an improved version of the embedding result (Theorem \ref{thm:main} and Corollary \ref{cor:main}), and follows from it by Theorem \ref{main} and Proposition \ref{prop:isom-module}. It states that the \emph{tangent module $L^2(\T\Y;\mu)$}, introduced by the second named author in \cite{Gigli14} (see also \cite{Gigli12}), admits an isometric embedding into $L^2(\T_G\Y;\mu)$ that is compatible with the fibrewise $\Cat 0$-structure on the target side. We refer to \cite{Gigli14,Gigli12} for the theory of tangent modules, and to Section \ref{sec:tangent} for the notation (below Theorem \ref{conecat}).

\begin{theorem}\label{main2}
	Let $\Y$ be a complete and separable locally $\Cat \kappa$-space ($\kappa\in\R$) and $\mu$ a Borel measure on $\Y$ that is finite on bounded sets. Then there is a map $\mathscr F:L^2({\rm T}\Y;\mu)\hookrightarrow L^2({\rm T}_G\Y;\mu)$ such that for $X,Y\in L^2({\rm T}\Y;\mu)$
	\begin{itemize}
		\item[(1)] $\mathscr F(X+Y)=\mathscr F(X)\oplus\mathscr F(Y)$, 
		\item[(2)] $\sfd_\cdot(\mathscr F(X),\mathscr F(Y))=|X-Y|$, and
		\item[(3)] $2(|\mathscr F(X)|_\cdot^2+|\mathscr F(Y)|_\cdot^2)=\sfd_\cdot^2(\mathscr F(X),\mathscr F(Y))+|\mathscr F(X)\oplus\mathscr F(Y)|_\cdot^2$
	\end{itemize}
	pointwise $\mu$-almost everywhere.
\end{theorem}

Both main results, along with Proposition \ref{prop:isom-module}, are proven in the end of Section 6.

\bigskip{\bf Acknowledgement.} This research has been supported by the MIUR SIR-grant `Nonsmooth Differential Geometry' (RBSI147UG4).
We would also like to thank S.\ Wenger for some useful comments on a preliminary version of this paper.

\section{ $\Cat \kappa$-spaces and basic  calculus on them}
\subsection{Definition of $\Cat \kappa$-spaces and basic properties}\label{se:cat}

In this paper \emph{geodesics} will always be assumed to be minimizing and with constant speed. If, for two given points $x,y$ in a metric space $(\Y,\sfd)$, there is only one (up to reparametrization) geodesic connecting them, the one defined on $[0,1]$  will be denoted by $\G_x^y$. Given a point $x\in\Y$, we denote by $\dist_x:\Y\to\R$ the function $y\mapsto \sfd(x,y)$.

For $\kappa\in\R$ the \emph{model space} $\M_\kappa$ is the connected, simply connected, complete 2-dimensional manifold with constant curvature $\kappa$, and $\sfd_\kappa$ is the distance induced by the metric tensor. Thus $(\M_\kappa,\sfd_\kappa)$ is (a) the hyperbolic space $\mathbb H^2_\kappa$ of constant sectional curvature $\kappa$, if $\kappa<0$, (b) $\R^2$ with the usual Euclidean metric, if $\kappa=0$, and (c) the sphere $S^2_\kappa$ of constant sectional curvature $\kappa$, if $\kappa>0$. 

 We set $D_\kappa:=\diam(\M_\kappa)$, i.e.
\begin{align*}
D_\kappa=\left\{
\begin{array}{ll}
\infty&\quad\text{ is }\kappa\le 0,\\
\frac{\pi}{\sqrt\kappa}&\quad\text{ if }\kappa>0.
\end{array}
\right.
\end{align*}
We refer to \cite[Chapter I.2]{BH99} for a detailed study of the model spaces $\M_\kappa$.

\bigskip

$\Cat\kappa$ spaces are geodesic  spaces where geodesic triangles are `thinner' than in $\M_\kappa$: they offer a metric counterpart to the notion of `having sectional curvature bounded from above by $\kappa$'.

To define them we start by recalling that if $a,b,c\in\Y$ is a triple of points satisfying $\sfd(a,b)+\sfd(b,c)+\sfd(c,a)<2D_\kappa$, then there are points, called \emph{comparison points}, $\bar a,\bar b,\bar c\in\M_\kappa$ such that
\[
\sfd_\kappa(\bar a,\bar b)=\sfd(a,b),\qquad\qquad\sfd_\kappa(\bar b,\bar c)=\sfd(b,c),\qquad\qquad\sfd_\kappa(\bar c,\bar a)=\sfd(c,a).
\]
A point $d\in\Y$ is said to be intermediate between $b,c\in\Y$ provided $\sfd(b,d)+\sfd(d,c)=\sfd(b,c)$ (if $\Y$ is geodesic, as we shall always assume, this means that $d$ lies on a geodesic joining $b$ and $c$). A \emph{comparison point of $d$} is a point $\bar d\in\M_\kappa$, such that
\[
\sfd_\kappa(\bar d,\bar b)=\sfd(d,b),\qquad\qquad\sfd_\kappa(\bar d,\bar c)=\sfd(d,c).
\]
\begin{definition}[$\Cat\kappa$ spaces]\label{cat}
A  metric space $(\Y,\sfd)$ is called a $\Cat\kappa $-space if it is geodesic and satisfies the following triangle comparison principle: for any $a,b,c\in \Y$, satisfying $\sfd(a,b)+\sfd(b,c)+\sfd(c,a)<2D_\kappa$, and any intermediate point $d$ between $b,c$, there are comparison points $\bar a,\bar b,\bar c,\bar d\in\M_\kappa$ as above such that
\begin{equation}
\label{eq:defcat}
\sfd(a,d)\le \sfd_\kappa(\bar  a,\bar d).
\end{equation}
A  metric space $(\Y,\sfd)$ is said to be locally $\Cat\kappa$ (or of curvature $\le \kappa$) if every point in $\Y$ has a neighbourhood which is a $\Cat\kappa$-space with the inherited metric. 
\end{definition}
It is worth noting that balls of radius  $<D_\kappa/2$ in the model space $\M_\kappa$ are convex, cf. Definition \ref{def:convex_sets}. Hence the  comparison property \eqref{eq:defcat} grants that the same is true on $\Cat\kappa$ spaces (see \cite[Proposition II.1.4.(3)]{BH99} for the rigorous proof of this fact). It is then easy to see that, for the same reasons, $(\Y,\sfd)$ is locally $\Cat\kappa$  provided every point has a neighbourhood $U$ where the comparison inequality (\ref{eq:defcat}) holds for every triple of points $a,b,c\in U$, where the geodesics connecting the points (and thus the intermediate points) are allowed to exit the neighbourhood $U$.

Let us fix the following notation: if $(\Y,\sfd)$ is a local $\Cat\kappa$ space, for every $x\in\Y$ we set
\[
\sfr_x:=\sup\big\{r\leq D_\kappa/2\ :\ \bar B_r(x)\ \text{is a $\Cat\kappa$ space}\big\}.
\]
Notice that in particular $B_{\sfr_x}(x)$ is a $\Cat\kappa$ space. The definition trivially grants that $\sfr_y\geq\sfr_x-\sfd(x,y)$ and thus in particular $x\mapsto \sfr_x$ is continuous.

We mention in passing that restricting attention to \emph{complete} $\Cat\kappa$-spaces presents no loss of generality, since the completion of a $\Cat\kappa$-space is a $\Cat\kappa$-space; see \cite[Corollary 3.11]{BH99}.

In a $\Cat\kappa$ space, points at distance $<D_\kappa$ are connected by a unique (up to parametrization) geodesic and these geodesics vary continuously with the endpoints. The following lemma is a quantitative version of this statement, and directly implies the uniqueness and continuous dependence of geodesics between points of distance $<D_\kappa$.
\begin{lemma}\label{le:quantmidpoint}
Let  $\kappa\in\R$ and let $\Y$ be a $\Cat\kappa$-space. For every $\lambda<D_\kappa$,  there are constants $C=C(\kappa,\lambda)>0$ and $\eps_0=\eps_0(\kappa,\lambda)>0$ such that the following holds: if $x,y\in \Y$ satisfy $\sfd(x,y)\le \lambda$, and $m$ is the midpoint of $x,y$, we have, for any $\eps\in (0,\eps_0)$ and $m'\in \Y$, that 
\[
\sfd(m,m')\le C\eps,\quad\textrm{ whenever }\quad\sfd^2(x,m'),\sfd^2(y,m')\leq\tfrac14\sfd^2(x,y)+\eps^2.
\]
\end{lemma}
\begin{proof}
By the definition of $\Cat\kappa$ space, using the triangle comparison property with the points $x,y,m'$, we see that it is sufficient to prove the claim when $\Y$ is the model space $\M_\kappa$. Since $\Cat\kappa$ spaces are $\Cat{\kappa'}$ spaces for $\kappa'\geq\kappa$ (see \cite[Part II, Chapter 1]{BH99}), we can assume that $\kappa>0$. Thus we may assume $\Y=S_\kappa^2$. In this case the conclusion follows by direct computations, one possible line of thought being the following.

Let $\eps_0$ be such that
\begin{equation}\label{eps0}
\frac{\eps^2_0}{2\cos(\sqrt\kappa\lambda/2)}<1\quad\textrm{ and }\quad \sqrt{\kappa(\eps_0^2+(\lambda/2)^2)}<\pi/2.
\end{equation}
Let $\eps\in (0,\eps_0)$, and let $x,y, m$ and $m'$ be as in the claim. 

Set $r_\eps:=\sqrt{\sfd(x,y)^2/4+\eps^2}$ and consider the set
\[
S:=\overline B_{S_\kappa^2}(x,r_\eps)\cap \overline B_{S_\kappa^2}(y,r_\eps)
\]
(note that $m'\in S$). The distance 
\[
\max_{s\in S}\sfd(s,m)
\]
is maximized at a point $s\in \partial S$ where the geodesic segment $[m,s]$ makes a right angle with the geodesic segment $[x,y]$. The spherical cosine law, applied to the triangle $\Delta(x,m,s)$ (resp. $\Delta(y,m,s)$) yields
\begin{equation}\label{cosinelaw}
\cos(\sqrt\sfd(m,s))\cos\frac{\sqrt\kappa\sfd(x,y)}{2}=\cos(\sqrt\kappa r_\eps).
\end{equation}
Denote $a:=\sfd(x,y)/2$ and define
\[
f(s):=\frac{\cos(\sqrt{\kappa(a^2+s^2)})}{\cos(\sqrt\kappa a)},\quad 0\le s\le \eps_0.
\]
Note that
\begin{eqnarray*}
1-f(\eps)\le \int_0^\eps|f'(s)|\d s \le  \int_0^\eps \frac{\kappa s\d s}{\cos(\sqrt\kappa a)}\le\frac{\kappa \eps^2}{\cos(\sqrt\kappa a)}.
\end{eqnarray*}
From this estimate, (\ref{cosinelaw}), and the fact that $a\le \lambda/2$, we have
\[
\cos(\sqrt\kappa\sfd(m,s))=f(\eps)\ge 1-\frac{\kappa \eps^2}{\cos(\sqrt\kappa \lambda/2)}.
\]
This, the elementary estimate $\displaystyle \arccos(1-t^2)\le 2t\ (0\le t<1)$ and (\ref{eps0}) then imply that 
\[
\sqrt\kappa\sfd(m,s)\le \frac{\sqrt\kappa\eps}{\sqrt{\cos(\sqrt\kappa \lambda/2)}}.
\]
This completes the proof.
\end{proof}

Being geodesic spaces, on $\Cat\kappa$ spaces it makes sense to speak about convex sets:
\begin{definition}[Convex sets and convex hull]\label{def:convex_sets}
Let $\Y$ be a $\Cat\kappa$ space. Then a set $C\subset \Y$ is said to be convex provided for any $x,y\in C$ we have that every geodesic connecting them is entirely contained in $C$.
The (closed) convex hull of a set $C\subset\Y$ is the smallest (closed) convex set containing $C$.
\end{definition}
One might define a weaker form of convexity by requiring that for every $x,y$ there exists a geodesic connecting them which is entirely contained in $C$. In $\Cat\kappa$ spaces this distinction is relevant only when $\sfd(x,y)\ge D_\kappa$, as otherwise geodesics are unique. For the purposes of the current manuscript the distinction is irrelevant.
\medskip

The following simple lemma will be useful later on:
\begin{lemma}[Separable convex hull]\label{le:sephull}
Let $\Y$ be a $\Cat\kappa$ space and $C\subset \Y$ a separable subset which is contained in a  closed ball $B$ of radius $<D_\kappa/2$.

Then the closed convex hull $\overline C_{\rm conv}$ of $C$ is separable and contained in $B$.
\end{lemma}
\begin{proof} Define the sequence $(C_n)$ of subsets of $\Y$ recursively as follows. Set $C_0:=C$, then iteratively let $C_{n+1}$ be the union of the images of geodesics whose endpoints are in $C_n$. It is clear that the convex hull of $C$ must contain $\cup_nC_n$ and thus $\overline C_{\rm conv}\supset \overline{\cup_nC_n}$. 

To conclude the proof it is therefore enough to show that $ \overline{\cup_nC_n}$ is convex and separable. The convexity of $ {\cup_nC_n}$ is a straightforward consequence of the definition using induction. Since $B$ is convex  we see that ${\cup_nC_n}\subset B$. Hence we have that $\sup_{x,y\in {\cup_nC_n}}\sfd(x,y)<D_\kappa$. By Lemma \ref{le:quantmidpoint}, the geodesic connecting two points $x,y\in \cup_nC_n$ depends continuously on $x$ and $y$. In particular, the separability of $C_{n+1}$ follows from that of $C_{n}$ (and the uniqueness of geodesics). Thus $ \overline{\cup_nC_n}$ is separable. By the continuous dependence of the (unique) geodesics and the convexity of $ {\cup_nC_n}$ the convexity of $ \overline{\cup_nC_n}$ follows.
\end{proof}
We conclude the section with the following result, taken from \cite[Part II, Lemma 3.20]{BH99}:
\begin{lemma}\label{le:contr}
Let $(\Y,\sfd)$ be a $\Cat\kappa$ space and $x\in\Y$. Then
there exists a function $C$ defined on a right neighbourhood
of $0$ such that $\lim_{r\downarrow 0}C(r)=1$ and
\begin{equation}
\frac{\sfd\big(({\sf G}_x^y)_\eps,({\sf G}_x^z)_\eps\big)}{\eps}
\leq C(r)\,\sfd(y,z)\qquad\text{ for every }\eps\in(0,1)
\text{ and }y,z\in B_r(x)
\end{equation}
for all $r<D_\kappa$ sufficiently small.
\end{lemma}
\subsection{Tangent cone}\label{sec:tangent} Here we define the tangent cone at a point on a $\Cat\kappa$ space and study its first properties. We refer the interested reader to the surveys \cite{BH99,BBI01, BS07} and the references therein for more details.

\bigskip

We start by describing a construction of tangent cone which is valid in every  geodesic space. Let $\Y$ be a geodesic space and $x\in \Y$.  We denote by $\Geo_x\Y$ the space of (constant speed) geodesics starting from $x$ and defined on some right neighbourhood of 0 and equip such space with the pseudo-distance $\sfd_x$ defined as:
\begin{equation}
\label{eq:defdx}
\sfd_x(\gamma,\eta):=\lims_{t\downarrow0}\frac{\sfd(\gamma_t,\eta_t)}{t}\qquad\forall \gamma,\eta\in\Geo_x\Y.
\end{equation}
Then $\sfd_x$ naturally induces an equivalence relation on $\Geo_x\Y$ by declaring $\gamma\sim\eta$ iff $\sfd_x(\gamma,\eta)=0$. The equivalence class of $\gamma\in\Geo_x\Y$ in $\Geo_x\Y/\sim$ will be denoted by $\gamma'_0$. Clearly $\sfd_x$ passes to the quotient and defines a distance -- still denoted by $\sfd_x$ -- on $\Geo_x\Y/\sim$ .

\begin{definition}[Tangent cone]
Let $\Y$ be a geodesic space and $x\in\Y$.  The tangent cone $(\T_x\Y,\sfd_x)$ is the completion of $(\Geo_x\Y/\sim,\sfd_x)$.  We call $0\in \T_x\Y$, or sometimes $0_x\in \T_x\Y$, the equivalence class of the constant geodesic in $\Geo_x\Y$.
\end{definition}
In a general geodesic space little can be said about the structure of tangent cones, but if $\Y$ is locally a $\Cat\kappa$ space then tangent cones have  interesting geometric properties and can be used as basic tools to build a robust first-order calculus.

In order to understand the geometry of $\T_x\Y$ it is necessary to recall the notion of \emph{angle} between geodesics. To do so, let us recall the definition of modified trigonometric functions
\[
\sn^\kappa(x):=\left\{\begin{array}{ll}
\tfrac1{\sqrt\kappa}\sin(\sqrt\kappa x)&\text{if }\kappa>0\\
x&\text{if }\kappa=0\\
\tfrac1{\sqrt{-\kappa}}\sinh(\sqrt{-\kappa} x)&\text{if }\kappa<0
\end{array}\right.
\qquad\qquad
\cn^\kappa(x):=\left\{\begin{array}{ll}
\cos(\sqrt\kappa x)&\text{if }\kappa>0\\
1&\text{if }\kappa=0\\
\cosh(\sqrt{-\kappa} x)&\text{if }\kappa<0
\end{array}\right.
\]
and that in the model space $\mathbb{M}_\kappa$ the cosine law reads, for $\kappa\neq 0$, as
\[
\cos(\alpha)=\frac{\cn^\kappa(a)-\cn^\kappa(b)\,\cn^\kappa(c)}{\kappa\,\sn^\kappa(b)\,\sn^\kappa(c)}
\]
whenever $a,b,c$ are the lengths of the sides of a geodesic triangle and $\alpha$ is the angle opposite to $a$ (in the limiting case $\kappa\to 0$ this reduces to the classical Euclidean cosine law). 

Then given three points $x,y_0,y_1$ in a metric space with $\sfd(x,y_0)+\sfd(x,y_1)+\sfd(y_0,y_1)<2D_\kappa$, we define the angle between $y_0,y_1$ seen from $x$ as
\begin{equation}
\label{eq:defan}
\overline\angle^\kappa_x(y_0,y_1):=\arccos\bigg(\frac{\cn^\kappa(\sfd(y_0,y_1))-\cn^\kappa(\sfd(x,y_0))\,\cn^\kappa(\sfd(x,y_1))}{\kappa\,\sn^\kappa(\sfd(x,y_0))\,\sn^\kappa(\sfd(x,y_1))}\bigg).
\end{equation}
Notice that this is the angle in the model space $\mathbb{M}_\kappa$ at $\bar x$ of a comparison triangle $\bar\Delta (\bar x,\bar y_0,\bar y_1)$ and from this observation it is not hard to check that
\begin{equation}
\label{eq:trangle}
\overline \angle^\kappa_x(y_0,y_2)\leq \overline \angle^\kappa_x(y_0,y_1)+\overline \angle^\kappa_x(y_1,y_2)
\end{equation}
for any four points $x,y_0,y_1,y_2$ in a metric space.

A direct consequence of the definition of $\Cat\kappa$ space and of the above cosine law is that on a  $\Cat\kappa$ space $\Y$, for $x\in\Y$ and $\gamma,\eta\in \Geo_x\Y$ 
\begin{equation}
\label{eq:monotan}\begin{split}
&\text{the angle $\overline\angle^\kappa_x(\gamma_t,\eta_s)$ is non-decreasing in both $t$ and $s$}\\
&\text{provided they vary in $\big\{(t,s)\,:\,\sfd(x,\gamma_t),\sfd(x,\eta_s)<D_\kappa\big\}$.}
\end{split}\end{equation}
Hence, if $\Y$ is a local $\Cat\kappa$ space, $x\in\Y$ and  $\gamma,\eta\in \Geo_x\Y$ the joint limit
\begin{equation}
\label{eq:anglelim}
\angle^\kappa_x(\gamma,\eta):=\lim_{t,s\downarrow0}\overline\angle^\kappa_x(\gamma_t,\eta_s)
\end{equation}
exists and it is called \emph{angle between the geodesics $\gamma,\eta$}. 
\medskip

The following technical result will be useful (for the proof see \cite[Lemma 3.3.1]{AKP} and the discussion thereafter).
\begin{lemma}[Independence of the angle on $\kappa$]\label{le:indkappa}Let  $\kappa_1,\kappa_2\in\R$, $\kappa_1\geq\kappa_2$. Then there is a constant $C=C(\kappa_1,\kappa_2)$ such that the following holds: for any metric space $\Y$ and  $x,y_1,y_2\in\Y$ with $\sfd(x,y_1),\sfd(x,y_2),\sfd(y_1,y_2)<D_{\kappa_1}$ it holds that
\[
\big|\overline\angle_x^{\kappa_1}(y_1,y_2)-\overline\angle_x^{\kappa_2}(y_1,y_2)\big|\leq C\sfd(x,y_1)\,\sfd(x,y_2).
\]
\end{lemma}
In particular, the angle $\angle_{x}^\kappa(\gamma,\eta)$ between geodesics $\gamma,\eta\in \Geo_x\Y$ does not depend on $\kappa$ and we shall drop the superscript from the notation. Picking $\kappa_1=0$ we see that, for any $\kappa\in \R$, we have 
\begin{equation}
\label{eq:angle0}
\cos(\overline\angle_x^\kappa(\gamma_t,\eta_s))=\frac{\sfd^2(\gamma_t,x)+\sfd^2(\eta_s,x)-\sfd^2(\gamma_t,\eta_s)}{2\sfd(\gamma_t,x)\sfd(\eta_s,x)}+o(ts).
\end{equation}
We drop the superscript from the notation of the comparison angle as well, with the understanding that $\kappa$ is fixed in each claim.

From \eqref{eq:trangle} it is not hard to check that $\angle_x$ is a pseudo-distance on $\Geo_x\Y$ and thus defines an equivalence relation $\sim'$ by declaring $\gamma\sim'\eta$ iff $\angle_x(\gamma,\eta)=0$. It is worth noticing that the angle between two different reparametrizations of the same geodesic is 0.

We denote by $\dir_x\Y$ the quotient $\Geo_x\Y/\sim'$ and, abusing a bit the notation, we keep denoting by $\angle_x$ and $\gamma\in\dir_x\Y$ the distance induced by $\angle_x$ and the equivalence class of $\gamma\in\Geo_x\Y$, respectively. 
\begin{definition}[Space of directions] Let $\Y$ be a local $\Cat\kappa$ space and $x\in\Y$. The space of directions $(\Sigma_x\Y,\angle_x)$ is the completion of $(\dir_x\Y,\angle_x)$.
\end{definition}
Let us now recall that given a generic metric space $(\X,\sfd_\X)$, the (Euclidean) \emph{cone} over it is the metric space $(C(\X),\sfd_{C(\X)})$ defined as follows (see also e.g.\ \cite{BBI01} for further details). As a set, $C(\X)$ is equal to $\big([0,\infty)\times\X\big)/\sim$, where $(t,x)\sim(s,y)$ iff $t=s=0$ or $(t,x)=(s,y)$. The distance is defined as
\begin{equation}
\label{eq:conedist}
\sfd_{C(\X)}^2\big((t,x),(s,y)\big):=t^2+s^2-2ts\cos\big(\sfd_\X(x,y)\wedge\pi\big).
\end{equation}
On $C(\X)$ there is a natural operation of `multiplication by a positive scalar': the product $\lambda z$ of $z=(t,x)$ by $\lambda\geq 0$ is defined as $(\lambda t,x)$.
\medskip

We then have the following:
\begin{theorem}[$\T_x\Y$ as a cone over the space of directions]\label{thm:samecone}
Let $\Y$ be a local $\Cat\kappa$ space. Fix a point $x\in\Y$. 
Then the $\lims$ in \eqref{eq:defdx} is a limit. Moreover, the map sending $\gamma\in \Geo_x\Y$ to $(\Lip(\gamma),\gamma)\in[0,\infty)\times\dir_x\Y$ passes to the quotient and uniquely extends to a bijective isometry from $\T_x\Y$ to $C(\Sigma_x\Y)$.
Finally, the map $B_{\sfr_x}(x)\ni y\mapsto (\G^y_x)'_0\in\T_x\Y$ is continuous. In particular, if $\mathcal D\subset B_{\sfr_x}(x)$ is dense in a neighbourhood of $x$, then $\{\lambda (\G^y_x)'_0\,:\,\lambda\geq 0,\ y\in\mathcal D\}$ is dense in $\T_x\Y$.
\end{theorem}
\begin{proof}
For any $\gamma,\eta\in \Geo_x\Y$, by picking $t=s$ in \eqref{eq:angle0} we  see that
\[
\frac{\sfd^2(\gamma_t,\eta_t)}{t^2}=\Lip(\gamma)^2+\Lip(\eta)^2-2\Lip(\gamma)\Lip(\eta)\cos(\overline\angle_x(\gamma_t,\eta_t))+o(t^2).
\]
Since $\lim_{t\downarrow0}\overline\angle_x(\gamma_t,\eta_t)=\angle_x(\gamma,\eta)$ it follows that the limit $\displaystyle \lim_{t\downarrow 0}\frac{\sfd^2(\gamma_t,\eta_t)}{t^2}$ exists, and equals
\[\begin{split}\
\lim_{t\downarrow 0}\frac{\sfd^2(\gamma_t,\eta_t)}{t^2}&=\Lip(\gamma)^2+\Lip(\eta)^2-2\Lip(\gamma)\Lip(\eta)\cos(\angle_x(\gamma,\eta))\\
&=\sfd^2_{C(\Sigma_x\Y)}\big((\Lip(\gamma),\gamma),(\Lip(\eta),\eta)\big).
\end{split}
\]
It follows that the map $\gamma_0'\mapsto (\Lip(\gamma),\gamma)$ defines a bijective isometry $\T_x\Y\mapsto C(\Sigma_x\Y)$.

For the continuity of $y\mapsto (\G^y_x)'_0$, notice that from the monotonicity of the angle it follows that 
\[
\angle_x\big((\G^{y}_x)'_0,(\G^z_x)'_0\big)\leq \overline \angle^\kappa_x(y,z)
\]
and thus if $z\to y$ we have $\angle_x\big((\G^{y}_x)'_0,(\G^z_x)'_0\big)\to 0$. Since trivially it also holds that ${\rm Lip}({\sf G}_{x}^{z})=\sfd(x,z)\to \sfd(x,y)={\rm Lip}({\sf G}_{x}^{y})$, continuity follows.

For the last claim, notice that by the definition of tangent cone and of multiplication by a positive scalar we have that $\{\lambda (\G^y_x)'_0\ :\ \lambda\geq 0,\ y\in B_r(x)\}$ is dense in $\T_x\Y$ for any $r\in(0,\sfr_x)$. Then the continuity just proved ensures that for any $\lambda\geq 0$ the set $\{\lambda (\G^y_x)'_0\ :\ y\in\mathcal D\}$ is dense in $\{\lambda (\G^y_x)'_0\ :\ y\in B_r(x)\}$, leading to the claim.
\end{proof}
A key property of the tangent cone is the following statement, which is central for our subsequent results. For the proof we refer to \cite[Chapter II, Theorem 3.19]{BH99}.
\begin{theorem}\label{conecat}
 Let $\Y$ be a  local $\Cat\kappa$ space and  $x\in\Y$. Then the tangent cone $(\T_x\Y,\sfd_x)$ is a  $\Cat 0$-space. 
\end{theorem}
The tangent cone $\T_x\Y$ is not only a $\Cat0$ space, but also comes with an additional structure which somehow resembles that of a Hilbert space. To make this more evident, let us introduce the following notation, valid for any $v,w\in \T_x\Y$ (see \cite{BS07,Petrunin07}).
\begin{itemize}
\item[a)] \emph{Multiplication by a positive scalar}. As for general cones, for $\lambda\geq 0$ and $v=(t,\gamma)\in \T_x\Y\approx C(\Sigma_x\Y)$ we put $\lambda v:=(\lambda t,\gamma)$.
\item[b)] \emph{Norm}. $|v|_x:=\sfd_x(v,0)$. 
\item[c)] \emph{Scalar product}. $\langle v,w\rangle_x: = \tfrac12\big[|v|_x^2+|w|_x^2-\sfd_x^2(v,w)\big]$.
\item[d)] \emph{Sum}. $v\oplus w:=2 m_{v,w}$, where $m_{v,w}$ is the midpoint of $v,w$ (well-defined because $\T_x\Y$ is a $\Cat0$ space).
\end{itemize}
The basic properties of these operations are collected in the following proposition:
\begin{proposition}[Basic calculus on the tangent cone]\label{hilbertine}
	Let $\Y$ be a  local $\Cat\kappa$  space and  $x\in \Y$.  Then the four operations defined above are continuous in their variables. The `sum' and the `scalar product' are also symmetric. Moreover:
\begin{subequations}
\begin{align}
\label{eq:norm}
|\lambda v|_x&=\lambda |v|_x,\\
\label{eq:distangle}
\sfd_x^2(v,w)&=|v|_x^2+|w|_x^2-2\langle v,w\rangle_x,\\
\label{eq:angolo}
\la\gamma'_0,\eta'_0\ra_x&=|\gamma'_0|_x|\eta'_0|_x\cos(\angle_x(\gamma,\eta)),\\
\label{eq:prhom}
\langle\lambda v,w\rangle_x&= \langle v,\lambda w\rangle_x=\lambda \langle v,w\rangle_x,\\
\label{eq:CS}
|\langle v,w\rangle_x|&\le |v|_x|w|_x,\\
\label{eq:CSeq}
\langle v,w\rangle_x&= |v|_x|w|_x\quad\text{ if and only if }\quad |w|_xv=|v|_xw,\\
\label{eq:PI}
\sfd_x^2(v,w)+|v\oplus w|_x^2&\le 2(|v|_x^2+|w|_x^2),
\end{align}
\end{subequations}	
for any $v,w\in \T_x\Y$, $\gamma,\eta\in\Geo_x\Y$ and $\lambda\geq 0$.
\end{proposition}
\begin{proof} The symmetry of the `sum' and `scalar product' are obvious and so are the continuity of the `norm' and then of the `scalar product'. The continuity of $(\lambda,v)\to \lambda v$ is a direct consequence of the inequality
\[
\sfd_x(\lambda v,\lambda'v')\leq \sfd_x(\lambda v,\lambda v')+\sfd_x(\lambda v',\lambda'v')=\lambda \sfd_x(v,v')+|\lambda'-\lambda| |v'|_x,
\]
where the equality follows trivially from the definition of cone distance and Theorem \ref{thm:samecone}. For the  continuity of the `sum'  it is now sufficient to prove that the map $(v,w)\mapsto m_{v,w}$ is continuous. This follows from the bound
\[
\sfd_x(m_{v,w},m_{v',w'})\leq\frac{1}2\big( \sfd_x(v,v')+\sfd_x(w,w')\big)
\] 
which is valid in any $\Cat0$ space (see e.g.\ \cite[Proposition 1.1.5 and Theorem 1.3.3]{Bac14}).

Notice that \eqref{eq:norm}, \eqref{eq:distangle} are direct consequences of the definitions. For \eqref{eq:angolo} we observe that by definition we have
\[
\la\gamma_0',\eta_0'\ra_x=\tfrac12\big(|\gamma'_0|_x^2+|\eta'_0|_x^2-\sfd_x^2(\gamma_0',\eta_0')\big)
\]
and thus recalling \eqref{eq:defdx} (and the fact that the $\lims$ is actually a limit -- see Theorem \ref{thm:samecone}) we obtain
\[
\la\gamma_0',\eta_0'\ra_x=|\gamma'_0|_x|\eta'_0|_x\lim_{t\downarrow0}\frac{\sfd^2(\gamma_t,x)+\sfd^2(\eta_t,x)-\sfd^2(\gamma_t,\eta_t)}{2\sfd(\gamma_t,x)\sfd(\eta_t,x)}\stackrel{\eqref{eq:angle0}}=|\gamma'_0|_x|\eta'_0|_x\cos(\angle_x(\gamma,\eta)).
\]

For \eqref{eq:prhom} we note that from the definition \eqref{eq:conedist} and Theorem \ref{thm:samecone} it is clear that $\sfd_x(\lambda v,\lambda w)=\lambda \sfd_x(v,w)$ for $\lambda\geq 0$. Hence we also have $\la\lambda v,\lambda w\ra_x=\lambda^2\la v,w\ra_x$ and thus, taking into account the symmetry of the scalar product, to conclude it is sufficient to prove that 
\begin{equation}
\label{eq:perscal}
\la\lambda v, w\ra_x\geq \lambda\la v,w\ra_x\qquad\text{  for}\ \lambda\in[0,1].
\end{equation}
To see this, notice that by the 2-convexity (\ref{eq:cat0}) of squared distance functions in $\Cat 0$-spaces we have, for $\lambda\in[0,1]$ and $v,w\in \T_x\Y$, that
\[
\sfd_x^2(\lambda v,w)\leq (1-\lambda)|w|_x^2+\lambda\sfd_x^2(v,w)-\lambda(1-\lambda)|v|_x^2.
\]
The estimate \eqref{eq:perscal} follows from this and the definition of $\la\cdot,\cdot\ra_x$.

To prove \eqref{eq:CS} let   $\gamma,\eta\in \dir_x\Y$ and $t,s\ge 0$ and observe that 
\[ 
\big|2\langle (t,\gamma),(s,\eta)\rangle_x\big| =\big| t^2+s^2-\sfd_x^2\big((t,\gamma),(s,\eta)\big)\big|=2ts|\cos\angle_x(\gamma,\eta)|
\le 2ts=2|(t,\gamma)|_x|(s,\eta)|_x.
\]
Since elements of the form $(t,\gamma)$ are dense in $\T_x\Y$, we just proved \eqref{eq:CS}. 

The `if' in \eqref{eq:CSeq} comes from \eqref{eq:prhom}, for the `only if' suppose that $\langle v,w\rangle_x=|v|_x|w|_x$ and take $\gamma_n,\eta_n\in \dir_x\Y$ so that $(|v|_x,\gamma_n)\to v$ and $(|w|_x,\eta_n) \to w$ in $\T_x\Y$. It follows that 
$$
|v|_x|w|_x=\langle v,w\rangle_x=\lim_{n\to\infty} \langle (|v|_x,\gamma_n),(|w|_x,\eta_n)\rangle_x=|v|_x|w|_x\lim_{n\to\infty}\cos\angle_x(\gamma_n,\eta_n),
$$ 
i.e.\ $\lim_{n\to\infty}\angle_x(\gamma_n,\eta_n)=0$. This implies that 
$$
\sfd_x(|w|_xv,|v|_xw)=\lim_n \sfd_x\big((|v|_x|w|_x,\gamma_n) ,(|v|_x|w|_x,\eta_n)\big)=|v|_x|w|_x\lim_n \sfd_x\big((1,\gamma_n) ,(1,\eta_n)\big)=0,
$$ 
which was the claim.

Finally, \eqref{eq:PI} is also a direct consequence of  the 2-convexity of the squared distance from a point, which gives
\[
\Big|\frac12(v\oplus w)\Big|_x^2\leq\frac12\big(|v|_x^2+|w|_x^2\big)-\frac14\sfd_x^2(v,w).
\]
Taking into account the proved homogeneities, this is the claim.
\end{proof}
It is worth underlying that, in general, $\oplus$ is not associative.
\begin{lemma}\label{le:formula_oplus}
Let $(\Y,\sfd)$ be a $\Cat\kappa$ space. Fix a point $x\in\Y$. Then for any
$y,z\in B_{D_\kappa}(x)\setminus\{x\}$ and $\alpha,\beta>0$ it holds that
\begin{equation}\label{eq:formula_oplus}
\big(\alpha({\sf G}_x^y)'_0\big)\oplus\big(\beta({\sf G}_x^z)'_0\big)
=\lim_{\eps\downarrow 0}\frac{2({\sf G}_x^{m_\eps})'_0}{\eps}\in\T_x\Y,
\end{equation}
where $m_\eps$ denotes the midpoint between $({\sf G}_x^y)_{\eps\alpha}$
and $({\sf G}_x^z)_{\eps\beta}$.
\end{lemma}
\begin{proof}
Let us call $p_\eps:=({\sf G}_x^y)_{\eps\alpha}$ and
$q_\eps:=({\sf G}_x^z)_{\eps\beta}$.
For $\eps,\delta>0$ small it holds that
$({\sf G}_x^y)_{\delta\eps\alpha}=({\sf G}_x^{p_\eps})_\delta$, whence by
using \eqref{eq:defdx} and Lemma \ref{le:contr} we obtain that
\[\begin{split}
\sfd_x\Big(\alpha({\sf G}_x^y)'_0,\frac{({\sf G}_x^{m_\eps})'_0}{\eps}\Big)
&=\frac{1}{\eps}\,\sfd_x\big(\eps\alpha({\sf G}_x^y)'_0,
({\sf G}_x^{m_\eps})'_0\big)=\frac{1}{\eps}\lim_{\delta\downarrow 0}
\frac{\sfd\big(({\sf G}_x^y)_{\delta\eps\alpha},
({\sf G}_x^{m_\eps})_\delta\big)}{\delta}\\
&=\frac{1}{\eps}\lim_{\delta\downarrow 0}
\frac{\sfd\big(({\sf G}_x^{p_\eps})_\delta,({\sf G}_x^{m_\eps})_\delta\big)}{\delta}
\leq\frac{1}{\eps}\,\sfd(p_\eps,m_\eps)=\frac{1}{2\eps}\,\sfd(p_\eps,q_\eps).
\end{split}\]
Similarly, we have that
$\sfd_x\big(\beta({\sf G}_x^z)'_0,\eps^{-1}({\sf G}_x^{m_\eps})'_0\big)
\leq\sfd(p_\eps,q_\eps)/(2\eps)$. Choosing $\eps_n\downarrow 0$ so that
\[\Big|\frac{\sfd(p_{\eps_n},q_{\eps_n})}{\eps_n}-
\sfd_x\big(\alpha({\sf G}_x^y)'_0,\beta({\sf G}_x^z)'_0\big)\Big|<\frac{2}{n}
\qquad\text{ for every }n,\]
we deduce that $\eps_n^{-1}({\sf G}_x^{m_{\eps_n}})'_0$ is a
$\frac{1}{n}$-approximate midpoint between $\alpha({\sf G}_x^y)'_0$ and
$\beta({\sf G}_x^z)'_0$. This yields \eqref{eq:formula_oplus}
by Lemma \ref{le:quantmidpoint}, as required.
\end{proof}
We close this section with the following important formula:
\begin{proposition}[First variation formula]\label{firstvar}
	Let $\Y$ be a $\Cat\kappa$ space, $x\in\Y$ and $\gamma,\eta\in  \Geo_x\Y$ with $\eta$ defined on $[0,1]$ and such that $\sfd(x,\eta_1)< D_\kappa$.  Then 
\begin{equation}
\label{eq:firstvar}
{\la\gamma'_0,\eta'_0\ra_x}=-\Lip(\eta)\lim_{t\downarrow 0}\frac{\sfd(\gamma_t,\eta_1)-\sfd(\gamma_0,\eta_1)}{t}. 
\end{equation}
\end{proposition}
\begin{proof}
We know from \eqref{eq:angolo} and \eqref{eq:angle0}  that
\[
\la\gamma'_0,\eta'_0\ra_x=\lim_{t,s\downarrow0}\frac{\sfd^2(\gamma_t,x)+\sfd^2(\eta_s,x)-\sfd^2(\gamma_t,\eta_s)}{2ts}
\]  
and by direct computation we see that
\[\begin{split}
\lim_{t\downarrow0}\frac{\sfd^2(\gamma_t,x)+\sfd^2(\eta_s,x)-\sfd^2(\gamma_t,\eta_s)}{2ts}&=\lim_{t\downarrow0}\frac{\sfd^2(\eta_s,x)-\sfd^2(\gamma_t,\eta_s)}{2ts}\\&=-\Lip(\eta)\lim_{t\downarrow0}\frac{\sfd(\gamma_t,\eta_s)-\sfd(x,\eta_s)}{t}.
\end{split}\]
Since the triangle inequality gives $\sfd(\gamma_t,\eta_1)-\sfd(x,\eta_1)\leq \sfd(\gamma_t,\eta_s)-\sfd(x,\eta_s)$, from the above we deduce
\[
\la\gamma'_0,\eta'_0\ra_x\leq-\Lip(\eta)\lim_{t\downarrow 0}\frac{\sfd(\gamma_t,\eta_1)-\sfd(\gamma_0,\eta_1)}{t}.
\]
Now notice that from   \eqref{eq:angolo}, Lemma \ref{le:indkappa},  the assumption  $\sfd(x,\eta_1)<D_\kappa$ and the monotonicity property \eqref{eq:monotan} we get
\[
\la\gamma'_0,\eta'_0\ra_x\geq\Lip(\gamma)\Lip(\eta)\lims_{t\downarrow0}\frac{\cn^\kappa(\sfd(\gamma_t,\eta_1))-\cn^\kappa(\sfd(\gamma_t,x))\,\cn^\kappa(\sfd(\eta_1,x))}{\kappa\,\sn^\kappa(\sfd(\gamma_t,x))\,\sn^\kappa(\sfd(\eta_1,x))}.
\]
Thus using the expansions
\[
\begin{split}
\cn^\kappa(\sfd(\gamma_t,x))&=1+O(t^2),\\
\cn^\kappa(\sfd(\gamma_t,\eta_1))&=\cn^\kappa(\sfd(x,\eta_1))-\kappa\,\sn^\kappa(\sfd(x,\eta_1))\big(\sfd(\gamma_t,\eta_1)-\sfd(x,\eta_1)\big)+O(t^2),\\
\sn^\kappa(\sfd(\gamma_t,x))&=t\Lip(\gamma)+O(t^2),
\end{split}
\]
we get the inequality $\geq$ in \eqref{eq:firstvar} and the conclusion.
\end{proof}

\subsection{Differential of locally semiconvex Lipschitz functions}\label{se:diff}
In this section we see that for Lipschitz and locally semiconvex functions there is a well-behaved notion of differential defined on the tangent cone of every point in the domain of the function itself. See \cite{Petrunin07,Plaut02} for the lower curvature bound case, and \cite{Lyt04} for more general classes of metric spaces.

We start by recalling the following notion:
\begin{definition}[Locally semiconvex function]\label{semiconvex}
Let $\Y$ be a geodesic metric space and $f:\Y\to\R$. We say that $f$ is semiconvex if there exists $K\in \R$ so that the inequality 
\[
f(\gamma_t)\leq(1-t)f(\gamma_0)+tf(\gamma_1)-\frac K2t(1-t)\sfd^2(\gamma_0,\gamma_1)
\]
holds for any geodesic $\gamma:[0,1]\to\Y$.

A function $f:\Omega\to\R$, with $\Omega\subset\Y$ open connected set, is called locally semiconvex if every point $x\in \Omega$ has a neighbourhood $U$ such that the inequality above holds for all geodesics $\gamma:[0,1]\to\Omega$ with endpoints in $U$.
\end{definition}
For locally semiconvex functions it is possible to define directional derivatives, which we do in the setting of $\Cat\kappa$ spaces:
\begin{definition}[Directional derivative]
Let $\Y$ be a local $\Cat\kappa$ space, $x\in \Y$,  $U\subset B_{\sfr_x}(x)$ a neighbourhood of $x$ and  $f:U\to\R$ locally semiconvex. The directional derivative of $f$ at $x$ is the map $\sigma_x f:\Geo_x\Y\to \R\cup\{-\infty\}$ defined as
\[
\sigma_x f(\gamma):=\lim_{h\downarrow0}\frac{f(\gamma_h)-f(\gamma_0)}{h}.
\]
\end{definition}
Notice that the monotonicity of incremental ratios of convex functions ensures that the limit above exists. Still, in general it is not clear if $\sigma_x f$ passes to the quotient $\Geo_x\Y/\sim$ nor if it is real-valued. In the next proposition we see that this is the case if we further assume that $f$ is  Lipschitz in a neighbourhood of $x$.

Recall that given $f:\Y\to\R$ the asymptotic Lipschitz constant $\lip_af: \Y\to[0,+\infty]$ is defined as
\[
\lip_af(x):=\lims_{y,z\to x}\frac{\big|f(y)-f(z)\big|}{\sfd(y,z)}=\lim_{r\downarrow0}\Lip(f\restr{B_r(x)})=\inf_{r>0}\Lip(f\restr{B_r(x)}).
\]

\begin{proposition}[Differentials of locally Lipschitz and  semiconvex functions]\label{lem:exdiff} Let $\Y$ be a  local $\Cat\kappa$ space, $\Omega\subset\Y$ open and $f:\Omega\to\R$ be locally semiconvex and Lipschitz.

Then for each $x\in \Omega$ there exists a unique continuous map $\d_xf:\T_x\Y\to\R$, called the differential of $f$ at $x$, such that 
\begin{equation}\label{eqn:diff}
	\d_xf(\gamma'_0)=\sigma_xf(\gamma)\qquad\forall \gamma\in \Geo_x\Y.
\end{equation}
Moreover, $\d_xf$ is convex, $\lip_a f(x)$-Lipschitz and positively 1-homogeneous, i.e\ $\d_xf(\lambda v)=\lambda \d_xf(v)$ for any $v\in\T_x\Y$ and $\lambda\geq 0$.
\end{proposition}
\begin{proof}
	Fix $x\in \Omega$ and let $r>0$ be such that $B_r(x)\subset \Omega$. Then for every $\gamma,\eta\in \Geo_x\Y$ we have $\gamma_t,\eta_t\in B_r(x)$ for $t\ll1$ and thus
	\[
\big|\sigma_x f(\gamma)-\sigma_x f(\eta)\big|\leq \lim_{h\downarrow0}\Big|\frac{f(\gamma_h)-f(\eta_h)}{h}\Big|\leq \Lip(f\restr{B_r(x)})\limi_{h\downarrow0}\frac{\sfd(\gamma_h,\eta_h)}{h}\leq \Lip(f\restr{B_r(x)})\sfd_x(\gamma,\eta).
	\]
	 This shows  that $\sigma_x f$ passes to the quotient and defines a $\Lip(f\restr{B_r(x)})$-Lipschitz map on $\Geo_x\Y/\sim$. Existence and uniqueness of the continuous extension $\d_x f$ to the whole $\T_x\Y$ are then obvious and, letting $r\downarrow0$, it is also clear that $\d_xf$ is $\lip_af(x)$-Lipschitz. 
	 
	 For the homogeneity observe that, for $\gamma\in\Geo_x\Y$ and $\lambda\geq0$, the isometry given in Theorem \ref{thm:samecone} and the definition of multiplication by positive scalar ensure that $\lambda\gamma'_0=\tilde\gamma'_0$ in $C(\Sigma_x\Y)\approx\T_x\Y$, where $\tilde\gamma_t:=\gamma_{\lambda t}$. Then \eqref{eqn:diff} and the definition of directional derivative grant that $\d_xf(\lambda\gamma'_0)=\lambda\d_xf(\gamma'_0)$ for any $\lambda\geq 0$ and $\gamma\in\Geo_x\Y$. Since tangent vectors of the form $\gamma'_0 $ are dense in $\T_x\Y$, the claim follows by the continuity of $\d_xf$ that we already proved.

It remains to prove that $\d_xf$ is convex and, thanks to the continuity just proven, it is sufficient to show that for any $\gamma,\eta\in \Geo_xB_{\sfr_x}(x)\simeq \Geo_x\Y$, letting $m$ be the midpoint of $\gamma'_0,\eta'_0\in \T_x\Y$ it holds that
\begin{equation}
\label{eq:forconv}
\d_x f(m)\leq\frac12\big(\d_x f(\gamma'_0)+\d_x f(\eta'_0)\big).
\end{equation}
To this aim, let $\eps>0$ and use the density of $ \Geo_xB_{\sfr_x}(x)/\sim$ in $\T_x\Y$ to find $\rho\in \Geo_xB_{\sfr_x}(x)$ such that $\rho'_0$ is an approximated midpoint of $\gamma'_0,\eta'_0$ in the sense that
\[
\sfd_x^2(\gamma'_0,\rho'_0),\sfd_x^2(\eta'_0,\rho'_0)\leq \tfrac14 \sfd_x^2(\gamma'_0,\eta'_0)+\eps^2.
\]
By the very definition \eqref{eq:defdx} of $\sfd_x$ we see that there exists $T>0$ such that
\[
\sfd^2(\gamma_t,\rho_t),\sfd^2(\eta_t,\rho_t)\leq\tfrac14 \sfd^2(\gamma_t,\eta_t)+2\eps^2t^2\qquad\forall t\in[0,T].
\]
Up to taking $T$ smaller, we can assume that $\sfd(\gamma_t,\eta_t)\leq  \frac12D_\kappa$, thus we are in a position to apply Lemma \ref{le:quantmidpoint} (in $\T_x\Y$ and $B_{\sfr_x}(x)$) to deduce that
\begin{equation}
\label{eq:withlip}
\begin{split}
\sfd_x(\rho'_0,m)&\leq C\eps,\\
\sfd(\rho_t,m_t)&\leq C\eps t\qquad\forall t\in[0,T],
\end{split}
\end{equation}
for some $C>0$ independent on $t$, where $m_t$ is the midpoint of $\gamma_t,\eta_t$.   Now let $V$ be a neighbourhood of $x$ where $f$ is $K$-semiconvex and $L$-Lipschitz and notice that what previously proved grants that $\d_x f$ is $L$-Lipschitz as well. Then $\gamma_t,\eta_t\in V$ for $t\ll1$ and the $K$-semiconvexity gives
\[
f(m_t)\leq \frac12\big(f(\gamma_t)+f(\eta_t)\big)-\frac K8\sfd^2(\gamma_t,\eta_t)\leq \frac12\big(f(\gamma_t)+f(\eta_t)\big)+\frac{| K|}8t^2\big(\Lip^2(\gamma)+\Lip^2(\eta)\big)
\]
and thus
\begin{equation}
\label{eq:withkconv}
\lims_{t\downarrow0}\frac{f(m_t)-f(x)}{t}\leq\frac12\lims_{t\downarrow0}\frac{f(\gamma_t)-f(x)}{t}+\frac12\lims_{t\downarrow0}\frac{f(\eta_t)-f(x)}{t}=\frac12\big(\d_xf(\gamma'_0)+\d_xf(\eta'_0)\big).
\end{equation}
Hence taking into account \eqref{eq:withlip} and the $L$-Lipschitz property of $f$ and $\d_xf$ we get
\[
\begin{split}
\d_xf(m)&\stackrel{\eqref{eq:withlip}}\leq CL\eps+\d_xf(\rho'_0)=CL\eps+\lim_{t\downarrow0}\frac{f(\rho_t)-f(x)}{t}\\
&\stackrel{\eqref{eq:withlip}}\leq 2CL\eps+\lims_{t\downarrow0}\frac{f(m_t)-f(x)}{t}\stackrel{\eqref{eq:withkconv}}\leq 2CL\eps+\frac12\big(\d_xf(\gamma'_0)+\d_xf(\eta'_0)\big).
\end{split}
\]
The conclusion follows letting $\eps\downarrow0$.
\end{proof}

In the model space $\M_\kappa$, if $\gamma$ is a geodesic on $[0,1]$ and $x\in\M_\kappa$ is such that $\sfd(\gamma_0,\gamma_1)+\sfd(\gamma_0,x)+\sfd(\gamma_1,x)<2D_\kappa$ we have that $t\mapsto\sfd(\gamma_t,x)$ is semiconvex. Hence if $\Y$ is a local $\Cat\kappa$ space and $x\in\Y$, for any $y\in B_{\sfr_x}(x)$ the function $\dist_y$ is semiconvex on $B_{\sfr_x}(x)$. 
\medskip

We collect below the main properties of the differential:
\begin{proposition}[Differentials of distance functions]\label{prop:diffdist}
Let $\Y$ be a $\Cat\kappa$ space and $x\in\Y$. Then:
\begin{itemize}
\item[i)] For $y\in \Y$ with $\sfd(x,y)<D_\kappa$ we have
\begin{equation}
\label{eq:ddist}
\d_x\dist_y(v)=\left\{\begin{array}{ll}
-\frac1{\Lip(\eta)}{\la v,\eta'_0\ra_x}&\quad\text{if }\ y\neq x,\\
\phantom{-}|v|_x&\quad\text{if }\ y=x,
\end{array}\right.\qquad\qquad\forall v\in\T_x\Y,
\end{equation}
where $\eta\in\Geo_x\Y$ is any geodesic passing through $y$.
\item[ii)] For $\mathcal D\subset B_{D_\kappa}(x)$ dense in a neighbourhood of $x$ we have
\[
|v|_x=\sup_{y\in \mathcal D}\big[-\d_x\dist_y(v)\big]\qquad\forall v\in\T_x\Y.
\]
\item[iii)] Let  $v,w\in\T_x\Y$ and $\mathcal D\subset B_{D_\kappa}(x)$ be dense. Assume that
\[
\d_x\dist_y(v)\le \d_x\dist_y(w) \qquad \forall y \in \mathcal{D}.
\]
Then $|w|^2_x \leq \la v,w\ra_x $ and in particular $|w|_x \leq |v|_x$.

If moreover either $|v|_x \leq |w|_x$ or $x \in \mathcal{D}$, then we also have $v=w$.
\end{itemize}
\end{proposition}
\begin{proof}\ \\
\noindent{\bf (i)} By the continuity of $v\mapsto\d_x\dist_y(v)$ and of the stated expression it is sufficient to check \eqref{eq:ddist} for $v$ of the form $v=\gamma'_0$ for arbitrary $\gamma\in\Geo_x\Y$. Then keeping in mind the identity \eqref{eqn:diff} and the definition of directional derivative we see that the case $y=x$ is obvious. For the case $y\neq x$ we notice that \eqref{eq:prhom} ensures  that $-\frac1{\Lip(\eta)}{\la v,\eta'_0\ra_x}$ does not depend on the particular choice of $\eta$. We choose $\eta$ to be defined on $[0,1]$ and such that $\eta_1=y$ and conclude noticing that the formula is a restatement of the first variation formula in Proposition \ref{firstvar}.

\noindent{\bf (ii)} Inequality $\geq$ follows from point (i) and the `Cauchy-Schwarz inequality' \eqref{eq:CS}.  The opposite inequality is trivial if $v=0$. If not, we use the density result in Theorem \ref{thm:samecone} to find $(y_n)\subset \mathcal D$ such that, letting $\gamma_n:[0,1]\to\Y$ be the geodesic from $x$ to $y_n$, we have $\frac1{\sfd(x,y_n)}\gamma_{n,0}'\to\frac{1}{|v|_x}v$. By point (i) (and recalling the calculus rules in Proposition \ref{hilbertine}) we have that
\[
-\d_x\dist_{y_n}(v)=\frac1{\sfd(x,y_n)}\la v,\gamma'_{n,0}\ra_x\quad \to\quad \frac{1}{|v|_x} \la v,v\ra_x=|v|_x.
\]
\noindent{\bf (iii)} If $w=0$ the first claim is obvious. Otherwise use the density result in Theorem \ref{thm:samecone} to find $(y_n)\subset\mathcal D$ such that, letting $\gamma_n:[0,1]\to\Y$ be the geodesic from $x$ to $y_n$, we have $\frac{1}{\sfd(x,y_n)}\gamma_{n,0}'\to \frac1{|w|_x}w$. Then point (i) and our assumption give $-\frac1{\sfd(x,y_n)}\la v,\gamma'_{n,0}\ra\leq-\frac1{\sfd(x,y_n)}\la w,\gamma'_{n,0}\ra$ and passing to the limit (using the calculus rules in Proposition \ref{hilbertine}) we get the first claim.

For the second claim, notice that if $x\in\mathcal D$, picking $y:=x$ in our assumption and using again point (i) we deduce $|v|_x\leq|w|_x$ (and thus $|v|_x=|w|_x$). Hence from what previously proved we obtain $\la v,w\ra_x\geq |w|_x^2\geq |v|_x|w|_x$, so that from \eqref{eq:CSeq} we conclude $|w|_xv=|v|_xw$ and from the equality of norms we conclude $v=w$, as desired.
\end{proof}

\subsection{Velocity of absolutely continuous curves} 

Recall that a curve $\gamma:[0,1]\to\Y$ with values in a metric space is said to be  \emph{absolutely continuous} provided there is $f\in L^1(0,1)$ such that
\begin{equation}
\label{eq:defac}
\sfd(\gamma_t,\gamma_s)\leq\int_t^s f(r)\,\d r\qquad\forall t,s\in[0,1],\ t<s.
\end{equation}
It is well-known that to any absolutely continuous curve we can associate a function $|\dot\gamma|\in L^1(0,1)$, called \emph{metric speed}, which plays the role of the  modulus of the derivative. The following proposition recalls the main properties of $|\dot\gamma|$; for the proof we refer to \cite[Theorem 1.1.2]{AmbrosioGigliSavare08} and its proof.
\begin{proposition}\label{prop:metsp}
Let $(\Y,\sfd)$ be a separable metric space and $\gamma:[0,1]\to\Y$ be absolutely continuous. Then for a.e.\ $t\in[0,1]$ there exists the limit $|\dot\gamma_t|$ as $h\to 0$ of $\frac{\sfd(\gamma_{t+h},\gamma_t)}{|h|}$. The function $t\mapsto|\dot\gamma_t|$ belongs to $L^1(0,1)$ and is the least, in the a.e.\ sense, function $f$ for which \eqref{eq:defac} holds.

Moreover, for any $(x_n)\subset\Y$ dense, letting $f_{n,t}:=\sfd(\gamma_t,x_n)$, the following holds: for a.e.\ $t\in[0,1]$ the function $f_n$ is differentiable at $t$ for every $n\in\N$ and
\begin{equation}
\label{eq:metsp}
|\dot\gamma_t|=\sup_{n\in\N}[-f_{n,t}'].
\end{equation}
\end{proposition}
On a local $\Cat\kappa$ space more can be said: for a.e.\ time we have not only a `numerical' value for the derivative, but also right and left derivatives as elements of the tangent cone. The key lemma needed for achieving such a result is the following (see also \cite[Theorem 1.6]{Lyt04}):
\begin{lemma}\label{le:velo}
Let $\Y$ be a local $\Cat\kappa$ space and $\gamma: [0,1] \to \Y$ be an absolutely continuous curve. 
Then for almost every $t\in[0,1]$ we have that either $|\dot\gamma_t|=0$ or
\begin{equation}
\label{eq:compangle}
\lims_{\delta_2\downarrow0}\lims_{\delta_1\downarrow0}\overline\angle_{\gamma_t}^\kappa(\gamma_{t+\delta_1},\gamma_{t+\delta_2})=0.
\end{equation}
\end{lemma} 
\begin{proof} The statement is local in nature, thus up to use the compactness of $\gamma([0,1])$, its cover made of $B_{\sfr_{\gamma_t}}(\gamma_t)$ and Lemma \ref{le:sephull} we can assume that $\Y$ is a separable $\Cat\kappa$ space of diameter $<D_\kappa$.

Let $t\in[0,1]$ be such that the metric derivative $|\dot\gamma_t|$ exists, is strictly positive and \eqref{eq:metsp} holds for some fixed countable dense $(x_n)\subset\Y$.  Then for every $x\in\Y$ we have
\[
\begin{split}
\limi_{\delta\downarrow0}\cos\big(\overline\angle_{\gamma_t}^\kappa(\gamma_{t+\delta},x)\big)=\limi_{\delta\downarrow0}\frac{\cn^\kappa(\sfd(x,\gamma_{t+\delta}))-\cn^\kappa(\sfd(\gamma_t,x))\cn^\kappa(\sfd(\gamma_t,\gamma_{t+\delta}))}{\kappa\,\sn^\kappa(\sfd(\gamma_t,\gamma_{t+\delta}))\sn^\kappa(\sfd(\gamma_t,x))}
\end{split}
\]
with obvious modifications for $\kappa=0$. Using the expansions
\[
\begin{split}
\cn^\kappa(\sfd(\gamma_t,\gamma_{t+\delta}))&=1+o(\delta),\\
\cn^\kappa(\sfd(\gamma_{t+\delta},x))&=\cn^\kappa(\sfd(\gamma_{ t},x))-\kappa\, \sn^\kappa(\sfd(\gamma_{ t},x))\big(\sfd(\gamma_{t+\delta},x)-\sfd(\gamma_{t},x)\big)+o(\delta),\\
\sn^\kappa(\sfd(\gamma_{ t},\gamma_{t+\delta}))&=\delta|\dot\gamma_t|+o(\delta),
\end{split}
\]
we obtain
\[
\begin{split}
\limi_{\delta\downarrow0}\cos\big(\overline\angle_{\gamma_t}^\kappa(\gamma_{t+\delta},x)\big)=-\lims_{\delta\downarrow0}\frac{\sfd(\gamma_{t+\delta},x)-\sfd(\gamma_{t},x)}{\delta|\dot\gamma_t|}.
\end{split}
\]
Picking $x=x_n$ and recalling that by assumption $s\mapsto f_{n,s}:=\dist_{x_n}(\gamma_s)$ is differentiable at $t$, we see that $\limi_{\delta\downarrow0}\cos\big(\overline\angle_{\gamma_t}^\kappa(\gamma_{t+\delta},x_n)\big)=-\frac{f_{n,t}'}{|\dot\gamma_t|}$. Hence by triangle inequality for angles \eqref{eq:trangle} we obtain
\[
\lims_{\delta_2\downarrow0}\lims_{\delta_1\downarrow0}\overline\angle_{\gamma_t}^\kappa(\gamma_{t+\delta_1},\gamma_{t+\delta_2})\leq \lims_{\delta_1\downarrow0}\overline\angle_{\gamma_t}^\kappa(\gamma_{t+\delta_1},x_n)+
\lims_{\delta_2\downarrow0}\overline\angle_{\gamma_t}^\kappa(\gamma_{t+\delta_2},x_n)=2\arccos\Big(-\frac{f_{n,t}'}{|\dot\gamma_t|}\Big)\quad\forall n\in\N.
\]
Taking the infimum in $n$ and using \eqref{eq:metsp} we conclude the proof.
\end{proof}
We then have the following result:
\begin{proposition}[Right derivative of AC curves]\label{prop:velo}
Let $\Y$ be a local $\Cat\kappa$ space and $\gamma: [0,1] \to \Y$ be an absolutely continuous curve. For any $t,s\in[0,1]$ write, for brevity,  $\G_t^s$ in place of $\G_{\gamma_t}^{\gamma_s}$ whenever this is well-defined.

Then for  every $t\in[0,1]$ for which $|\dot\gamma_t|$ exists and the conclusion of Lemma \ref{le:velo} holds (and thus in particular for a.e.\ $t$) we have that:
\begin{itemize}
\item[i)] the limit, denoted by $\dot \gamma^+_t$ in $\T_{\gamma_t}\Y$, of $\frac{1}{s-t}(\G^s_t)'_0$ as $s\downarrow t$ exists, and
\item[ii)] for every locally Lipschitz and locally semiconvex function $f$ defined on some neighbourhood of $\gamma_t$ it holds that
\begin{equation}
\label{eq:derdiff}
\lim_{h\downarrow 0}\frac{f(\gamma_{t+h})-f(\gamma_t)}{h}=\d_{\gamma_t}f(\dot\gamma^+_t).
\end{equation}
\end{itemize}
\end{proposition} 
\begin{proof} \ \\
\noindent{\bf (i)} We shall prove that  $s\mapsto \frac{1}{s-t}(\G^s_t)'_0\in\T_{\gamma_t}\Y$ has a limit as $s\downarrow t$ for any $t$ for which $|\dot\gamma_t|$ exists and the conclusions of Lemma \ref{le:velo} hold. Notice that $| \frac{1}{s-t}(\G^s_t)'_0|_{\gamma_t}=\frac{\sfd(\gamma_t,\gamma_s)}{|s-t|}\to|\dot\gamma_t|$, thus if $|\dot\gamma_t|=0$ the conclusion follows. If $|\dot\gamma_t|>0$, by the convergence of norms  that we just proved and recalling \eqref{eq:conedist} and Theorem \ref{thm:samecone}, to conclude it is sufficient to prove that $\lim_{s_2\downarrow t}\lim_{s_1\downarrow t}\angle_{\gamma_t}((\G_t^{s_2})'_0,(\G_t^{s_1})'_0)=0$. This is a direct consequence of \eqref{eq:compangle} and the monotonicity property \eqref{eq:monotan}, which ensures that
\[
\angle_{\gamma_t}((\G_t^{s_2})'_0,(\G_t^{s_1})'_0)\leq \overline\angle_{\gamma_t}(\G_t^{s_2}(1),\G_t^{s_1}(1))=\overline\angle^\kappa_{\gamma_t}(\gamma_{t+s_2},\gamma_{t+s_1}).
\]

\noindent{\bf (ii)} If $|\dot\gamma_t|=0$ both sides of \eqref{eq:derdiff} are easily seen to be zero, so that the conclusion follows. Otherwise, for $\delta_2\geq\delta_1>0$ put for brevity $\eta_{\delta_1,\delta_2}:=(\G_{\gamma_t}^{\gamma_{t+\delta_2}})_{\delta_1/\delta_2}$ and notice that the monotonicity \eqref{eq:monotan} of angles gives
\[
\cos\big(\overline\angle^\kappa_{\gamma_t}(\gamma_{t+\delta_1},\eta_{\delta_1,\delta_2})\big)\geq\cos\big(\overline\angle^\kappa_{\gamma_t}(\gamma_{t+\delta_1},\gamma_{t+\delta_2})\big).
\]
From \eqref{eq:angle0} we have that 
\[
\cos\big(\overline\angle^\kappa_{\gamma_t}(\gamma_{t+\delta_1},\eta_{\delta_1,\delta_2})\big)=\frac{\sfd^2(\gamma_t,\gamma_{t+\delta_1})+\sfd^2(\gamma_t,\eta_{\delta_1,\delta_2})-\sfd^2(\gamma_{t+\delta_1},\eta_{\delta_1,\delta_2})}{2\sfd(\gamma_t,\gamma_{t+\delta_1})\sfd(\gamma_t,\eta_{\delta_1,\delta_2})}+O(\delta_1)
\]
and thus using the identity $\sfd(\gamma_t,\eta_{\delta_1,\delta_2})=\frac{\delta_1}{\delta_2}\sfd(\gamma_t,\gamma_{t+\delta_2})$ and passing to the limit recalling \eqref{eq:compangle} we deduce that
\begin{equation}
\label{eq:finally}
\lims_{\delta_2\downarrow0}\lims_{\delta_1\downarrow0}\frac{\sfd^2(\gamma_{t+\delta_1},\eta_{\delta_1,\delta_2})}{\delta_1^2}=0.
\end{equation}
Now let  $L$ be the Lipschitz constant of $f$ in some neighbourhood of $\gamma_t$ and notice that for $\delta_2>0$ sufficiently small we have
\[
\begin{split}
\lims_{\delta_1\downarrow0}\Big|\d_{\gamma_t}f\big(\tfrac1{\delta_2}(\G_t^{t+\delta_2})'_0\big)-\frac{f(\gamma_{t+\delta_1})-f(\gamma_t)}{\delta_1}\Big|&=\lims_{\delta_1\downarrow0}\Big|\frac{f(\eta_{\delta_1,\delta_2})-f(\gamma_t)}{\delta_1}-\frac{f(\gamma_{t+\delta_1})-f(\gamma_t)}{\delta_1}\Big|\\
&\leq L\lims_{\delta_1\downarrow0}\frac{\sfd(\eta_{\delta_1,\delta_2},\gamma_{t+\delta_1})}{\delta_1}.
\end{split}
\]
The conclusion follows by letting $\delta_2\downarrow0$ and using \eqref{eq:finally}.
\end{proof}
\begin{remark}\label{rmk:velo}{\rm
The conclusions of Lemma \ref{le:velo} and Proposition \ref{prop:velo} hold for left derivatives as well; for a.e.\ $t\in \{ |\dot\gamma_t|>0 \}$ we have
$\lim_{\delta\downarrow 0}\lim_{\eps\downarrow 0}\overline \angle_{\gamma_t}(\gamma_{t-\delta},\gamma_{t-\eps})=0$, and for every $t$ satisfying it
the limit of $\frac{1}{t-s}(\G^s_t)'_0$ as $s\uparrow t$ exists. We denote this limit by $\dot\gamma_t^-$.
\fr}\end{remark}
\begin{lemma}\label{lem:velo}
	Let $\gamma:[0,1]\to \Y$ be an absolutely continuous curve. Then, for almost every $t\in [0,1]$, the limits $\dot\gamma^+_t$ and $\dot\gamma^-_t$ are antipodal, i.e.
	\[
	\dot\gamma_t^+\oplus\dot\gamma_t^-=0.
	\]
\end{lemma}
\begin{proof}
We prove that 
\begin{equation}\label{anglepi}
\angle_{\gamma_t}(\dot\gamma_t^+,\dot\gamma_t^-)=\lim_{\delta\downarrow 0}\overline \angle_{\gamma_t}(\gamma_{t+\delta},\gamma_{t-\delta})=\pi
\end{equation}
for almost every $t\in [0,1]$. The claim follows from this.

\medskip\noindent If $|\dot\gamma_t|=0$ then $\dot\gamma_t^+=\dot\gamma_t^-=0$ and the claim is clear.

By Lemma \ref{le:velo}, Proposition \ref{prop:velo} and Remark \ref{rmk:velo}, almost every $t\in \{ |\dot\gamma_t|>0 \}$ satisfies conditions (i) and (ii) below.

(i) $\lim_{\delta\downarrow 0}\lim_{\eps\downarrow 0}\overline\angle_{\gamma_t}(\gamma_{t+\delta},\gamma_{t+\varepsilon})=0$ and $\lim_{\delta\downarrow 0}\lim_{\eps\downarrow 0}\overline\angle_{\gamma_t}(\gamma_{t-\delta},\gamma_{t-\varepsilon})=0$;

(ii) $\lim_{\delta\downarrow 0}\frac{\sfd(\gamma_{t+\delta},\gamma_{t})}{\delta}=\lim_{\delta\downarrow 0}\frac{\sfd(\gamma_{t-\delta},\gamma_{t})}{\delta}=\lim_{\delta\downarrow 0}\frac{\sfd(\gamma_{t-\delta},\gamma_{t+\delta})}{2\delta}=|\dot\gamma_t|$ (including the existence of these limits).

We fix $t\in [0,1]$ satisfying (i) and (ii).  Note that, by the monotonicity of angles (\ref{eq:monotan}), we have the estimate
\[\begin{split}
\angle_{\gamma_t}(\dot\gamma_t^+,\dot\gamma_t^-)=&\lim_{\delta\downarrow 0}\angle_{\gamma_t}((\G_t^{t+\delta})_0',(\G_t^{t-\delta})_0')=\lim_{\delta\downarrow 0}\lim_{\varepsilon\downarrow 0}\overline\angle_{\gamma_t}((\G_t^{t+\delta})_\eps,(\G_t^{t-\delta})_\eps)\\
\le &\lim_{\delta\downarrow 0}\overline \angle_{\gamma_t}(\gamma_{t+\delta},\gamma_{t-\delta}).
\end{split}\]
To prove the opposite inequality, we use the triangle inequality (\ref{eq:trangle}) to obtain
\begin{align*}
\overline\angle_{\gamma_t}((\G_t^{t+\delta})_\eps,(\G_t^{t-\delta})_\eps)&\ge \overline\angle_{\gamma_t}((\G_t^{t+\delta})_\eps,\gamma_{t-\eps\delta})-\overline\angle_{\gamma_t}(\gamma_{t-\eps\delta},(\G_t^{t-\delta})_\eps)\\
&\ge \overline\angle_{\gamma_t}(\gamma_{t+\eps\delta},\gamma_{t-\eps\delta})-\overline\angle_{\gamma_t}(\gamma_{t+\eps\delta},(\G_t^{t+\delta})_\eps)-\overline\angle_{\gamma_t}(\gamma_{t-\eps\delta},(\G_t^{t-\delta})_\eps)\\
&\ge \overline\angle_{\gamma_t}(\gamma_{t+\eps\delta},\gamma_{t-\eps\delta})-\overline\angle_{\gamma_t}(\gamma_{t+\delta},\gamma_{t+\eps\delta})-\overline\angle_{\gamma_t}(\gamma_{t-\delta},\gamma_{t-\eps\delta}).
\end{align*}
Here the last estimate follows simply by the monotonicity of angles (\ref{eq:monotan}). By (i), it follows that
\[
\angle_{\gamma_t}(\dot\gamma_t^+,\dot\gamma_t^-)=\lim_{\delta\downarrow 0}\lim_{\eps\downarrow 0}\overline\angle_{\gamma_t}((\G_t^{t+\delta})_\eps,(\G_t^{t-\delta})_\eps)\ge \lim_{\delta\downarrow 0}\lim_{\eps\downarrow 0}\overline\angle_{\gamma_t}(\gamma_{t+\eps\delta},\gamma_{t-\eps\delta})=\lim_{\delta\downarrow 0}\overline \angle_{\gamma_t}(\gamma_{t+\delta},\gamma_{t-\delta}).
\]

\medskip\noindent It remains to show that $\lim_{\delta\downarrow 0}\overline \angle_{\gamma_t}(\gamma_{t+\delta},\gamma_{t-\delta})=\pi$. By (ii) and (\ref{eq:angle0}), we have 
\begin{align*}
\lim_{\delta\downarrow 0}\cos\overline \angle_{\gamma_t}(\gamma_{t+\delta},\gamma_{t-\delta})&=\lim_{\delta\downarrow 0}\frac{\sfd^2(\gamma_{t+\delta},\gamma_t)+\sfd^2(\gamma_{t-\delta},\gamma_t)-\sfd^2(\gamma_{t+\delta},\gamma_{t-\delta})}{2\sfd(\gamma_{t+\delta},\gamma_t)\sfd(\gamma_{t-\delta},\gamma_t)}\\
&=\frac{|\dot\gamma_t|^2+|\dot\gamma_t|^2-(2|\dot\gamma_t|)^2}{2|\dot\gamma_t|^2}=-1,
\end{align*}
implying the claim and completing the proof.	
\end{proof}

\subsection{Barycenters and rigidity}\label{se:bari}
In this section we review the concept of `barycenter' of a probability measure on a $\Cat0$ space. With the exception of the rigidity statement given by Proposition \ref{lem:rigidity}, the content comes from \cite{stu03}.

\bigskip

Fix a $\Cat0$ space $\Y$ and denote by  $\mathscr P(\Y)$ the set of all Borel probability measures on $\Y$ having separable support, and  by $\mathscr P_1(\Y)\subset \mathscr P(\Y)$ the set of those with finite first moment, i.e.\ those $\mu\in\mathscr P(\Y)$ such that for some, and thus all, $y\in \Y$ it holds that $\int \sfd(\cdot,y)\,\d\mu<\infty$.

For a proof of the following result we refer to \cite[Proposition 4.3]{stu03}.
\begin{proposition}[Definition of barycenter]\label{prop:barycenter}
	Let $\Y$ be a $\Cat{0}$ space, $\mu\in\mathscr P_1(\Y)$ and $y\in \Y$. Then
	\[
	\Y\ni x\longmapsto\int\big[\sfd^2(\cdot,x)-\sfd^2(\cdot,y)\big]\,\d\mu\in\R
	\]
	admits a unique minimizer. The minimizer does not depend on $y$, is called the  \emph{barycenter} of $\nu$ and is denoted by $\Ba(\nu)\in \Y$.
\end{proposition}
The basic properties of barycenters that we shall need are collected in the following statement:
\begin{theorem}\label{thm:bar} Let $\Y$ be a $\Cat{0}$ space. Then the following holds:
\begin{itemize}
\item[i)] \textsc{Variance inequality.} For any $\mu\in\mathscr P_1(\Y)$ and $p\in \Y$ it holds
\begin{equation}\label{eq:variance_ineq}
\int\big[\sfd^2(\cdot,p)-\sfd^2(\cdot,\Ba(\mu))\big]\,\d\mu\geq \sfd^2(p,\Ba(\mu)).
\end{equation}
\item[ii)] \textsc{Jensen's inequality.} Let $\varphi:\Y\to[0,+\infty)$ be convex and lower semicontinuous. Then for every $\mu\in\mathscr P_1(\Y)$ we have
\begin{equation}\label{eq:Jensen}
\varphi(\Ba(\mu))\leq\int\varphi\,\d\mu.
\end{equation}
\end{itemize}
\end{theorem}
\begin{proof}
The variance inequality \eqref{eq:variance_ineq} is proved in  \cite[Proposition 4.4]{stu03}, while  Jensen's inequality comes from \cite[Theorem 6.2]{stu03}.
\end{proof}
Applying Jensen's inequality \eqref{eq:Jensen} to the convex and Lipschitz function $\varphi:=\sfd(\cdot,p)$ we see that the inequality
\[
\sfd(\Ba(\mu),p)\leq \int \sfd(x,p)\,\d\mu(x)
\]
holds for any $\mu\in\mathscr P_1(\Y)$ and $p\in \Y$. Our aim is now to study the equality case and in order to do so we first recall the notion of  nonbranching geodesics.
\begin{definition}[Non-branching from $p$]\label{def:nonbranch}
We say that a geodesic space $(\X,\sfd)$ is non-branching from $p\in \X$ provided the following holds: if, for given points $q,x_1,x_2\in \X$ with $q\neq p$, we have that there are geodesics $\gamma_1,\gamma_2$ starting from $p$ and passing through $q,x_1$ and $q,x_2$ respectively, then there is a geodesic $\gamma$ starting from $p$ and passing through $q,x_1,x_2$. 
\end{definition}
Here `passing through' $q,x_i$ implies nothing about the order in which these points are met. It is not hard to see that the above definition is equivalent to the more classical one  requiring for any $t\in(0,1]$  the injectivity of the map $\gamma\mapsto\gamma\restr{[0,t]}$ on the space of constant speed geodesics $[0,1]\to\X$ starting from $p$.

It is easy to verify that if $q\neq p$ and $(x_i)\subset\X$ are given points such that there are geodesics starting from $p$ and passing through $q,x_i$ for every $i$, then  there is a curve $\gamma$ starting from $p$ and passing through $q$ and all the $x_i$'s and such curve is either a geodesic or a half-line, i.e.\ a map from $[0,+\infty)$ to $\X$ such that its restriction to any compact interval is a geodesic.

The main example of space that is non-branching from one of its points is the cone over a metric space. Here the relevant point is the vertex 0.
\begin{lemma}[Tangent cones are non-branching from the origin]
Let $\X$ be any  metric space, and $C(\X)$ the Euclidean cone over $\X$. Then $C(\X)$ is non-branching from its origin 0. In particular, for a local $\Cat\kappa$ space $\Y$ we have that
\begin{equation}
\T_p \Y\text{ is non-branching from }0\qquad\text{ for every }p\in\Y.
\end{equation}
\begin{proof} By direct computation based on the definition of the cone distance we see that if $\gamma$ is a constant speed geodesic starting from the origin $0$ it must hold $\gamma_t=t\gamma_1$, where the `product' of $t$ and $\gamma_1$ is defined as before Proposition \ref{hilbertine}.  Thus for two given such curves $\gamma,\eta$ we have -- again using the definition of distance on the cone -- that $\sfd_{C(\X)}(\gamma_t,\eta_t)=\sfd_{C(\X)}(t\gamma_1,t\eta_1)=t \sfd_{C(\X)}(\gamma_1,\eta_1)$. Hence if $\gamma_1\neq\eta_1$ we also have $\gamma_t\neq \eta_t$ for every $t\in(0,1]$. This is sufficient to conclude.
\end{proof}
\end{lemma}
We now come to the rigidity statement:
\begin{proposition}[Rigidity]\label{lem:rigidity}
Let $\Y$ be a $\Cat{0}$ space and $\mu\in\mathscr P_1(\Y)$.
Assume that for some point $p$ it holds that
\begin{equation}\label{eqn:equality} 
\sfd(\Ba(\mu),p) \geq \int \sfd(x,p) \, \d \mu.
\end{equation}
Then 
\begin{equation}
\label{eq:rigid}
\sfd(x,\Ba(\mu))=\big|\sfd(x,p)-\sfd(\Ba(\mu),p)\big|\quad\text{ for }\mu\text{-a.e.\ }x\in \Y.
\end{equation}
In particular, if $\Y$ is non-branching from $p$, then the measure $\nu$ is concentrated on the image of a curve $\gamma$ starting from $p$ which is either a geodesic or a half-line.
\end{proposition}
\begin{proof}
By the discussion following Definition \ref{def:nonbranch}, we see that it is sufficient to prove \eqref{eq:rigid}. To this aim, notice that the triangle inequality gives
\begin{equation}
\label{eq:pointgeq}
\big|\sfd(x,p)-\sfd(x,\Ba(\mu))\big|^2-\sfd^2(\Ba(\mu),p)\leq 0\qquad\forall x\in \Y.
\end{equation}
On the other hand we have
\[
\begin{split}
&\int\big|\sfd(x,p)-\sfd(\Ba(\mu),p)\big|^2-\sfd^2(x,\Ba(\mu))\,\d\mu(x)\\
=&\int \sfd^2(x,p)+\sfd^2(\Ba(\mu),p)-2\sfd(x,p)\sfd(\Ba(\mu),p)-\sfd^2(x,\Ba(\mu))\,\d\mu(x)\\
\textrm{by \eqref{eq:variance_ineq}}\qquad\geq&\int2  \sfd^2(p,\Ba(\mu))-2\sfd(x,p)\sfd(x,\Ba(\mu))\,\d\mu(x)\\
=&2\sfd(p,\Ba(\mu))\Big(\sfd(p,\Ba(\mu))-\int \sfd(x,p)\,\d\mu(x)\Big)\\
\textrm{by \eqref{eqn:equality}}\qquad\geq&0.
\end{split}
\]
This inequality and \eqref{eq:pointgeq} give \eqref{eq:rigid} and the conclusion.
\end{proof}

\begin{remark}{\rm
It is easily seen that in the preceding proposition the non-branching assumption is needed. Indeed, consider the `tripod', i.e.\ the $\Cat 0$-space $\Y$ obtained as the Euclidean cone over the space $\{a,b,c\}$ equipped with the discrete metric. Then $\Y$ is \emph{not} non-branching from $a$ and, indeed, the conclusion of Proposition \ref{lem:rigidity} fails for the measure $ \mu=\frac 13(\delta_a+\delta_b+\delta_c)$, even though the identity (\ref{eqn:equality}) holds for $\mu$. Note that in this case $\Ba(\mu)=0$.
\fr}\end{remark}

\section{Geometric tangent bundle $L^2(\T_G\Y;\mu)$}\label{se:gtb}
In this section we fix a separable local $\Cat\kappa$ space $\Y$. Our first aim here is to give a measurable structure to the `geometric tangent bundle' $\T_G\Y$, i.e.\ the collection of all tangent cones on $\Y$. Once this is done, we will endow $\Y$ with a non-negative and non-zero Radon measure  $\mu$ and study the space of `$L^2$-sections' of $\T_G\Y$, which we shall denote by $L^2(\T_G\Y;\mu)$.

\bigskip\noindent As a set, the \emph{geometric tangent bundle} $\T_G\Y$ is defined as  
$$
\T_G\Y:=\big\{(x,v)\;\big|\;x\in\Y,\,v\in\T_x\Y\big\}.
$$ 
We denote by $\pi^\Y:\T_G\Y\to \Y$ the canonical projection defined by $\pi^\Y(x,v):=x$ and call \emph{section} of $\T_G\Y$ a map $v:\Y\to \T_G\Y$ such that $\pi^\Y(v(x))=x$ for every $x\in \Y$. 

We now endow $\T_G\Y$ with a $\sigma$-algebra $\mathcal B(\T_G\Y)$, defined as the smallest $\sigma$-algebra such that:
\begin{itemize}
\item[i)] The projection map $\pi^\Y:\T_G\Y\to \Y$ is measurable, $\Y$ being equipped with Borel sets.
\item[ii)] For every $x\in \Y$ and $y\in B_{\sfr_x}(x)$ the map
$\d\,\dist_y:\,\T_G\Y\to\R$, defined as
\[(\d\,\dist_y)(z,v):=\left\{\begin{array}{ll}
\d_z\dist_y(v)\\
0
\end{array}\qquad\begin{array}{ll}
\text{ if }(z,v)\in(\pi^\Y)^{-1}\big(B_{\sfr_x}(x)\big),\\
\text{ otherwise,}
\end{array}\right.\]
is measurable.
\end{itemize}
It is clear that these define a $\sigma$-algebra $\mathcal B(\T_G\Y)$, to which we shall refer as the class of Borel subsets of $\T_G\Y$, hereafter speaking about Borel (rather than measurable) maps. This is a slight abuse of terminology, since we are not defining any topology on $\T_G\Y$. The abuse of terminology is justified by the fact that if $\Y$ is a smooth Riemannian manifold, then $\mathcal B(\T_G\Y)$ coincides with the $\sigma$-algebra of Borel subsets of the tangent bundle of $\Y$.
\medskip

The following result gives a basic description of $\mathcal B(\T_G\Y)$:
\begin{proposition}\label{prop:btgy}
Let $\Y$ be a local $\Cat\kappa$ space which is also separable and $(x_n)\subset\Y$ a countable set of points such that $\bigcup_nB_{\sfr_{x_n}}(x_n)=\Y$ (these exist by the Lindel\"of property of $\Y$). For each $n$, let $(x_{n,m})\subset B_{\sfr_{x_n}}(x_n)$ be countable and dense.

Then $\mathcal B(\T_G\Y)$ coincides with the smallest $\sigma$-algebra $\mathcal B'(\T_G\Y)$  satisfying $i)$ above and 
\begin{itemize}
\item[ii')] For every $n,m\in\N$ the function
$\d\,\dist_{x_{n,m}}$ is measurable.
\end{itemize}
Moreover, for any $x\in\Y$ the measurable structure induced on $\T_x\Y$ by $\mathcal B(\T_G\Y)$ coincides with the Borel structure of $(\T_x\Y,\sfd_x)$.
\end{proposition}
\begin{proof} 
It is clear that $\mathcal B'(\T_G\Y)\subset \mathcal B(\T_G\Y)$. To prove the other inclusion start observing that the continuity of $x\mapsto\sfr_x$ grants that $B_{\sfr_x}(x)\subset\bigcup_n B_{\sfr_{x_n}}(x_n)$ if $x_n\to x$, thus to conclude it is sufficient to show that for given $n\in\N$ and $y\in B_{\sfr_{x_n}}(x_n)$ the map $(\pi^\Y)^{-1}(B_{\sfr_{x_{n}}}(x_n))\ni (x,v)\mapsto \d_x\dist_y(v)$ is $\mathcal B'(\T_G\Y)$-measurable. Keeping in mind point $(i)$ of Proposition \ref{prop:diffdist}, this will be achieved if we prove that:
\begin{itemize}
\item[a)] the map $\T_y\Y\ni v\mapsto |v|_y$ is measurable w.r.t.\ the $\sigma$-algebra induced by $\mathcal B'(\T_G\Y)$,
\item[b)] $(\pi^\Y)^{-1}(B_{\sfr_{x_{n}}}(x_n))\ni (x,v)\mapsto \la v,(\G_x^y )'_0\ra_x$ is $\mathcal B'(\T_G\Y)$-measurable.
\end{itemize}
Point (a) is a direct consequence of point $(ii)$ of Proposition \ref{prop:diffdist}. For (b), we notice that by assumption the claim is true if $y=x_{n,m}$ for some $n,m$. Then the general case follows from the continuity of the scalar product established in Proposition \ref{hilbertine} and the continuity of the map $B_{\sfr_{x_n}}(x_n)\ni y\mapsto (\G_x^y)'_0$ proved in Theorem \ref{thm:samecone}.

For the second claim, denote by $\mathcal B(\T_x\Y)$ the collection of Borel sets in $(\T_x\Y,\sfd_x)$ and by $\mathcal A_x$ the $\sigma$-algebra induced by $\mathcal B(\T_G\Y)$ on $\T_x\Y$. Then  the continuity of the `norm' and `scalar product' on $\T_x\Y$ proved in Proposition \ref{hilbertine} and the already recalled point $(i)$ of Proposition \ref{prop:diffdist} give the inclusion $\mathcal A_x\subset \mathcal B(\T_x\Y)$. For the opposite inclusion it is sufficient to prove that for any $v\in\T_x\Y$ the map $\T_x\Y\ni w\mapsto \sfd_x^2(v,w)$ is $\mathcal A_x$-measurable. Since in (a) above we have already proved that $\T_x\Y\ni w\mapsto |w|_x$ is $\mathcal A_x$-measurable, by \eqref{eq:distangle} it is sufficient to prove that $\T_x\Y\ni w\mapsto \la w,v\ra_x$ is $\mathcal A_x$-measurable as well. For $v$ of the form $(\G_x^y)'_0$ for some $y\in B_{\sfr_x}(x)$ this can be proved as in (b) above. Then the general case follows by the positive homogeneity of the scalar product given in \eqref{eq:prhom} and the density result in Theorem \ref{thm:samecone}.
\end{proof}

\begin{corollary} Let $\Y$ be a local $\Cat\kappa$ space which is also  separable. Then
$\mathcal B(\T_G\Y)$ is countably generated.
\end{corollary}
\begin{proof}
By definition, the $\sigma$-algebra $\mathcal B'(\T_G\Y)$ defined in Proposition \ref{prop:btgy} is countably generated. Thus the same holds for $\mathcal B(\T_G\Y)$.
\end{proof}

\begin{corollary}\label{cor:normborel} Let $\Y$ be a local $\Cat\kappa$ space which is also  separable. Let us denote by ${\sf Norm}:\T_G\Y\to\R^+$ the map sending $(x,v)$ to $|v|_x$. Then ${\sf Norm}$ is a Borel function.
\end{corollary}
\begin{proof}
Given that $x\mapsto \sfr_x$ is continuous, for any $z\in\Y$ we can find
$\lambda_z\in(0,\sfr_z)$ such that $B_{\lambda_z}(z)\subset B_{\sfr_x}(x)$
whenever $x\in B_{\lambda_z}(z)$. By Lindel\"{o}f property, to get
the statement it is sufficient to prove that
$(\pi^\Y)^{-1}\big(B_{\lambda_z}(z)\big)\ni(x,v)\to|v|_x$ is Borel for
every $z\in\Y$. Fix $z\in\Y$ and choose a dense
sequence $(y_n)\subset B_{\lambda_z}(z)$. We know from
item ii) of Proposition \ref{prop:diffdist} that
$|v|_x=-\inf_n\d_x\dist_{y_n}(v)$ holds for every
$(x,v)\in(\pi^\Y)^{-1}\big(B_{\lambda_z}(z)\big)$. Thus the required measurability follows from the definition of $\mathcal B(\T_G\Y)$, which grants that
$(x,v)\mapsto\d_x\dist_{y_n}(v)$ is measurable on $(\pi^\Y)^{-1}\big(B_{\lambda_z}(z)\big)$
for every choice of $n\in\N$.
\end{proof}

We shall say that a section $v:\Y\to\T_G\Y$ is \emph{simple} provided there
are $(y_n)\subset \Y$, $(\alpha_n)\subset\R^+$ and a Borel partition  $(E_n)$
of $\Y$ such that for every $n\in\N$ and $x\in E_n$ we have $y_n\in B_{\sfr_x}(x)$
and $v(x)=\alpha_n(\G_x^{y_n})'_0$. We will use the notation  $v=\sum_n\nchi_{E_n}\alpha_n(\G_\cdot^{y_n})'_0$ for simple sections. Notice that, arguing as in the proof of Proposition \ref{prop:btgy}, we see that simple sections are automatically Borel.
\medskip

The following lemma will be useful:
\begin{lemma}\label{le:simplesect}
Let $\Y$ be a separable local $\Cat\kappa$ space.
Then for every Borel section $v$ of $\T_G\Y$ and $\eps>0$ there is
a simple section $\tilde v$ such that $\sfd_x\big(v(x),\tilde v(x)\big)<\eps$
for every $x\in\Y$.
\end{lemma}
\begin{proof}
Using the Lindel\"of property of $\Y$ and the covering made by $B_{\sfr_x/2}(x)$ it is easy to see that we can reduce to the case in which $\Y$ is $\Cat\kappa$ and, for any $x,y\in\Y$, it holds that $y\in B_{\sfr_x}(x)$. Assume this is the case and let $(y_n)\subset\Y$ be countable and dense. By Theorem \ref{thm:samecone}, for every $x\in\Y$ the set $\big\{r(\G_x^{y_n})'_0\,:\,r\in\Q^+,\ n\in\N\big\}$ is dense in $\T_x\Y$. Let $i\mapsto (r_i,y_{n_i})$ be an enumeration of the couples $(r,y_n)$ with $r\in\Q^+$ and $n\in\N$.
Given a Borel section $v$, define
\[
E_i:=\Big\{x\in \Y\;\Big|\;\text{$i$ is the least index $j$ such that $\sfd_x\big(v(x),r_j(\G_x^{y_j})'_0\big)<\eps$}\Big\}
\]
and notice that, since the map $x\mapsto \sfd_x\big(v(x),r_i(\G_x^{y_i})'_0\big)$ is Borel for every $i$, the sets $E_i$ are Borel. The density result previously recalled ensures that $\bigcup_iE_i=\Y$. It follows that $\tilde v:=\sum_i\nchi_{E_i}r_i(\G_\cdot^{y_i})'_0$ fulfills the requirements.
\end{proof}
\begin{corollary}
Let $\Y$ be a local $\Cat\kappa$ space which is also separable.
Let $f:\,\Y\to\R$ be a locally semiconvex, locally Lipschitz function
and $v$ a Borel section of $\T_G\Y$. Then $\Y\ni x\mapsto\d_x f\big(v(x)\big)$
is a Borel function.
\end{corollary}
\begin{proof}
In light of Lemma \ref{le:simplesect}, it is sufficient to prove
the statement for simple sections. Let 
$v=\sum_n\nchi_{E_n}\alpha_n({\sf G}_\cdot^{y_n})'_0$ be simple and observe that, for every $x\in\Y$, one has that
\[\d_x f\big(v(x)\big)=\sum_n\nchi_{E_n}(x)\,\alpha_n\,\d_x f\big(({\sf G}_x^{y_n})'_0\big)
=\sum_n\nchi_{E_n}(x)\,\alpha_n\lim_{h\downarrow 0}
\frac{f\big(({\sf G}_x^{y_n})_h\big)-f(x)}{h}.\]
Since the function $E_n\ni x\mapsto f\big(({\sf G}_x^{y_n})_h\big)-f(x)$ is continuous
for all $h\in(0,1)$ by Lemma \ref{le:quantmidpoint}, we conclude
that $\Y\ni x\mapsto\d_x f\big(v(x)\big)$ is Borel,
thus completing the proof of the statement.
\end{proof}

The approximation result Lemma \ref{le:simplesect} also links the notion of Borel sections to Borel functions on $\Y$.
\begin{proposition}\label{prop:boroper} Let $\Y$ be a separable local $\Cat\kappa$ space. Let $v,w$ be Borel sections of $\T_G\Y$ and $\lambda\geq 0$.
Then it holds that
\[\begin{split}
\Y\ni x&\longmapsto |v|_x,\\
\Y\ni x&\longmapsto\sfd_x\big(v(x),w(x)\big),\\
\Y\ni x&\longmapsto\big\langle v(x),w(x)\big\rangle_x
\end{split}\]
are Borel functions. Moreover, $\lambda v$ and $v\oplus w$ are Borel sections of $\T_G\Y$.
\end{proposition}
\begin{proof} For the first part of the statement it is sufficient to prove that
	$x\mapsto\sfd_x\big(v(x),w(x)\big)$ is Borel, by the definition of `norm' and of `scalar product'. As for the proof of
Lemma \ref{le:simplesect} above, we use the Lindel\"of property of $\Y$
and the covering made of the balls $B_{\sfr_x/2}(x)$, $x\in \Y$, to reduce to the case of a
$\Cat\kappa$ space $\Y$ such that $y\in B_{\sfr_x}(x)$ for every $x,y\in\Y$.
By Lemma \ref{le:simplesect} it is sufficient to prove the
claim for simple sections $v,w$. Let $v=\sum_i\nchi_{E_i}\alpha_i(\G_\cdot^{y_i})'_0$
and $w=\sum_j\nchi_{F_j}\beta_j(\G_\cdot^{z_j})'_0$ be simple, and notice that
\[\sfd_x\big(v(x),w(x)\big)=\sum_{i,j}\nchi_{E_i\cap F_j}(x)\,
\sfd_x\big(\alpha_i(\G_x^{y_i})'_0,\beta_j(\G_x^{z_j})'_0\big).\]
The Borel regularity of $x\mapsto\sfd_x\big(v(x),w(x)\big)$ will
follow if we show that $x\mapsto\sfd_x\big(\alpha(\G_x^y)'_0,\beta(\G_x^z)'_0\big)$
is Borel for every $y,z\in\Y$ and $\alpha,\beta>0$. To this aim notice that,
since geodesics in $\Y$ are unique, they depend continuously (w.r.t.\ uniform
convergence) on their endpoints (see also Lemma \ref{le:quantmidpoint}).
Therefore, for every $t\in(0,1)$, we have that
$x\mapsto\sfd\big((\G_x^y)_{\alpha t},(\G_x^z)_{\beta t}\big)$ is continuous
and the conclusion follows recalling that, by \eqref{eq:conedist} 
and Theorem \ref{thm:samecone}, we have
\[\sfd_x\big(\alpha(\G_x^y)'_0,\beta(\G_x^z)'_0\big)=
\lim_{n\to\infty}\frac{\sfd\big((\G_x^y)_{\alpha t_n},(\G_x^z)_{\beta t_n}\big)}{t_n},\]
where $(t_n)$ is any sequence decreasing to $0$.

It is straightforward to see that $\lambda v$ is a Borel section of $\T_G\Y$:
the function $\Y\ni x\mapsto\d_x\dist_y\big(\lambda v(x)\big)=
\lambda\,\d_x\dist_y\big(v(x)\big)$ is Borel for every $y\in\Y$,
whence $\lambda v$ is a Borel section.

We now aim to prove that $v\oplus w$ is a Borel section of $\T_G\Y$.
By Lemma \ref{le:simplesect} it is enough to show that
$\Y\setminus\{y,z\}\ni x\mapsto\d_x\dist_p\big(\alpha({\sf G}_x^y)'_0,
\beta({\sf G}_x^z)'_0\big)$ is Borel for every $p,y,z\in\Y$ and $\alpha,\beta>0$.
By Lemma \ref{le:formula_oplus} and the properties of
$\d_x\dist_p$ we have
\[\d_x\dist_p\big(\alpha({\sf G}_x^y)'_0,\beta({\sf G}_x^z)'_0\big)
=\lim_{\eps\downarrow 0}\d_x\dist_p\big(2\eps^{-1}({\sf G}_x^{m_\eps(x)})'_0\big)\\
=\lim_{\eps\downarrow 0}\lim_{h\downarrow 0}
\frac{\sfd\big(p,({\sf G}_x^{m_\eps(x)})_{2h/\eps}\big)-\sfd(p,x)}{h},\]
where $m_\eps(x)$ stands for the midpoint between $({\sf G}_x^y)_{\eps\alpha}$
and $({\sf G}_x^z)_{\eps\beta}$. Given that the map
sending $x\in\Y$ to $({\sf G}_x^{m_\eps(x)})_{2h/\eps}\in\Y$ is continuous
(as one can see by repeatedly applying Lemma \ref{le:quantmidpoint}), we
conclude that $\Y\setminus\{y,z\}\ni x\mapsto\d_x\dist_p\big(\alpha({\sf G}_x^y)'_0,
\beta({\sf G}_x^z)'_0\big)$ is Borel, as required. This completes the proof
of the statement.
\end{proof}

We now consider the `right derivative' map ${\sf RightDer}:\,C([0,1];\Y)\times[0,1]\to\T_G\Y$, given by
\begin{equation}\label{eq:rightder}
{\sf RightDer}(\gamma,t):=
\left\{\begin{array}{ll}
\displaystyle\big(\gamma_t,\lim_{h\downarrow 0}h^{-1}(\G_{\gamma_t}^{\gamma_{t+h}})'_0\big),&\qquad\text{ if the limit in $\T_{\gamma_t}\Y$ exists,}\\
(\gamma_t,0),&\qquad\text{ otherwise}.
\end{array}\right.
\end{equation}

\begin{proposition}\label{prop:rightder}
Let $\Y$ be a separable local $\Cat\kappa$ space. Then ${\sf RightDer}:\,C([0,1];\Y)\times[0,1]\to \T_G\Y$ is a Borel map.
\end{proposition}
\begin{proof}
Let us denote by ${\rm e}:\,C([0,1];\Y)\times[0,1]\to\Y$ the evaluation
map $(\gamma,t)\mapsto\gamma_t$, which is clearly continuous. In order
to show that $\sf RightDer$ is Borel it suffices to prove that:
\begin{itemize}
\item[$\rm i)$] $\pi^\Y\circ{\sf RightDer}$ is Borel,
\item[$\rm ii)$] $\d\,\dist_y\circ{\sf RightDer}$ is Borel for every
$x\in\Y$ and $y\in B_{\sfr_x}(x)$.
\end{itemize}
Item i) trivially follows from the observation that ${\rm e}=\pi^\Y\circ{\sf RightDer}$.
To prove ii), fix $x\in\Y$ and $y\in B_{\sfr_x}(x)$.
Let us define the sets $D'$, $D$ and $S_h$ for $h\in(0,1)$ as follows:
\[\begin{split}
D'&:=\big\{(\gamma,t)\in C([0,1];\Y)\times[0,1)\;\big|\;
\gamma_t\in B_{\sfr_x}(x)\big\}={\rm e}^{-1}\big(B_{\sfr_x}(x)\big),\\
D&:=\big\{(\gamma,t)\in D'\;\big|\;\lim_{h\downarrow 0}
h^{-1}({\sf G}_{\gamma_t}^{\gamma_{t+h}})'_0\text{ exists}\big\},\\
S_h&:=\big\{(\gamma,t)\in D'\;\big|\;t+h\in[0,1),\,
\gamma_{t+h}\in B_{\sfr_{\gamma_t}}(\gamma_t)\big\},
\end{split}\]
respectively. Since  $\rm e$ and $x\mapsto \sfr_x$ are continuous, we have
that $D'$ and $S_h$ are open. Notice that
\[\begin{split}
(\d\,\dist_y\circ{\sf RightDer})(\gamma,t)&\overset{\phantom{\eqref{eqn:diff}}}=
\nchi_D(\gamma,t)\,\lim_{h\downarrow 0}\d_{\gamma_t}\dist_y
\big(h^{-1}({\sf G}_{\gamma_t}^{\gamma_{t+h}})'_0\big)\\
&\overset{\eqref{eqn:diff}}=\nchi_D(\gamma,t)\,
\lim_{h\downarrow 0}\lim_{\eps\downarrow 0}\nchi_{S_h}(\gamma,t)
\frac{\sfd\big(y,({\sf G}_{\gamma_t}^{\gamma_{t+h}})_\eps\big)-\sfd(y,\gamma_t)}{\eps h},
\end{split}\]
where the first equality stems from the continuity of $\T_{\gamma_t}\Y\ni v\mapsto \d_{\gamma_t}\dist_y(v)$.
Given any $h,\eps\in(0,1)$, the map
$(\gamma,t)\mapsto\nchi_{S_h}(\gamma,t)\big[\sfd\big(y,
({\sf G}_{\gamma_t}^{\gamma_{t+h}})_\eps\big)-\sfd(y,\gamma_t)\big]/(\eps h)$
is continuous on $S_h$ by Lemma \ref{le:quantmidpoint}. Thus, to obtain the measurability of the function $\d\,\dist_y\circ{\sf RightDer}$, it remains to show that the
set $D$ is Borel. To this aim, let us set
\[A_{h_1,h_2,\eps}:=\Big\{(\gamma,t)\in S_{h_1}\cap S_{h_2}\;\Big|\;
\sfd_{\gamma_t}\big(h_1^{-1}({\sf G}_{\gamma_t}^{\gamma_{t+h_1}})'_0,
h_2^{-1}({\sf G}_{\gamma_t}^{\gamma_{t+h_2}})'_0\big)<\eps\Big\}\]
for every $h_1,h_2\in(0,1)$ and $\eps>0$.
Given that for all $(\gamma,t)\in S_{h_1}\cap S_{h_2}$ we can write
\[\sfd_{\gamma_t}\big(h_1^{-1}({\sf G}_{\gamma_t}^{\gamma_{t+h_1}})'_0,
h_2^{-1}({\sf G}_{\gamma_t}^{\gamma_{t+h_2}})'_0\big)
\overset{\eqref{eq:defdx}}=\lim_{\delta\downarrow 0}
\frac{\sfd\big(({\sf G}_{\gamma_t}^{\gamma_{t+h_1}})_{\delta/h_1},
({\sf G}_{\gamma_t}^{\gamma_{t+h_2}})_{\delta/h_2}\big)}{\delta},\]
we can deduce (by Lemma \ref{le:quantmidpoint}) that each set
$A_{h_1,h_2,\eps}$ is Borel. Finally, observe that
\[D=\bigcap_{\eps\in\Q^+}\bigcup_{\substack{h\in\Q^+ \\ h<1}}
\bigcap_{\substack{h_1,h_2\in\Q^+ \\ h_1,h_2<h}}A_{h_1,h_2,\eps},\]
whence the set $D$ is Borel. The statement follows.
\end{proof}

We now fix a non-negative and non-zero Radon measure $\mu$ on $\Y$. We are interested in Borel sections of $\T_G\Y$ which are also in $L^2(\mu)$.
\begin{definition}[The space $L^2(\T_G\Y;\mu)$]\label{def:L2section}
	Let $\Y$ be a separable local $\Cat\kappa$-space and $\mu$ a non-negative non-zero Radon measure on $\Y$. The space $L^2(\T_G\Y;\mu)$ is defined as
\[
L^2(\T_G\Y;\mu):=\Big\{v\text{ Borel section of }\T_G\Y\ :\ \int|v|_x^2\,\d\mu(x)<\infty\Big\}/\sim,
\]
where $v\sim w$ if $\{x:v(x)=w(x)\}$ is $\mu$-negligible. We endow  $L^2(\T_G\Y;\mu)$ with the distance
\[
\sfd_{\mu}(v,w):=\sqrt{\int\sfd^2_x\big(v(x),w(x)\big)\,\d\mu(x)}.
\]
\end{definition}
Notice that, by Proposition \ref{prop:boroper}, the integrals in Definition \ref{def:L2section} are well-defined. With a (common) abuse of notation we do not distinguish between a Borel section $v$ and its equivalence class up to $\mu$-a.e.\ equality.
\medskip

We conclude the section collecting some basic properties of $\big(L^2(\T_G\Y;\mu),\sfd_\mu\big)$:
\begin{proposition}\label{thm:tangcat} Let $\Y$ be a separable local $\Cat\kappa$ space and $\mu$ a non-negative, non-zero Radon measure on it. Then $\big(L^2(\T_G\Y;\mu),\sfd_\mu\big)$ is a complete and separable $\Cat0$ space.
\end{proposition}
\begin{proof} The fact that $\sfd_\mu$ is a distance on $L^2(\T_G\Y;\mu)$ is trivial, so we turn to the other properties.

\textsc{Completeness.} The argument is standard: as it is well-known, it is sufficient to prove that any $(v_n)\subset L^2(\T_G\Y;\mu)$ such that $\sum_n\sfd_\mu(v_n,v_{n+1})<\infty$ is convergent. Then from the inequality 
\[
\Big\|\sum_n\sfd_\cdot(v_n,v_{n+1})\Big\|_{L^2(\mu)}\leq\sum_n\big\|\sfd_\cdot(v_n,v_{n+1})\big\|_{L^2(\mu)}=\sum_n\sfd_\mu(v_n,v_{n+1})<\infty
\]
we see that $\sum_n\sfd_\cdot(v_n,v_{n+1})\in L^2(\mu)$ and in particular that $\sum_n\sfd_x(v_n(x),v_{n+1}(x))<\infty$ for $\mu$-a.e.\ $x$. For any such $x$ the sequence $(v_n(x))$ is Cauchy in $\T_x\Y$ and thus has a limit $v(x)$. It is then clear that $v$ is (the equivalence class up to $\mu$-a.e.\ equality of) a Borel section. Moreover, by Fatou's lemma and the definition of $\sfd_\mu$ we see that
\[
\lims_{n}\sfd_\mu(v,v_n)\leq \lims_n\limi_m\sfd_\mu(v_n,v_m)=0,
\]
having used again the assumption that $(v_n)$ is $\sfd_\mu$-Cauchy. This proves that $v$ is the $\sfd_\mu$-limit of $(v_n)$ and, since this fact and the triangle inequality for $\sfd_\mu$ also tell that $v\in L^2(\T_G\Y;\mu)$, the claim is proved.

\textsc{Separability.} Using the Lindel\"of property of $\Y$ and the very definition of distance $\sfd_\mu$ we can reduce to the case in which $\Y$ is a separable $\Cat\kappa$ space with diameter $<D_\kappa$ and $\mu$ is a finite measure. Then taking into account Lemma \ref{le:simplesect} above it is easy to see that to conclude it is sufficient to find a countable set $\mathcal D\subset L^2(\T_G\Y;\mu)$ whose closure contains all simple sections of the form $v=\nchi_E(\G_\cdot^y)'_0$ for generic $E\subset\Y$ Borel and $y\in\Y$. Let $\mathcal D_1\subset \Y$ be countable and dense and $\mathcal D_2\subset\mathcal B(\Y)$ be countable and such that for any $E\subset\Y$ Borel and $\eps>0$ there is $U\in \mathcal D_2$ such that $\mu(U\Delta E)<\eps$ (for instance, the family of all finite unions of open balls having center in $\mathcal D_1$ and rational radius does the job -- by regularity of the measure $\mu$).

We then define 
\[
\mathcal D:=\Big\{\nchi_E(\G_\cdot^y)'_0\ :\ y\in\mathcal D_1,\ E\in\mathcal D_2 \Big\}
\]
and claim that this does the job. To see this, notice that the inequality
\[
\sfd_\mu^2\big(\nchi_E(\G_\cdot^y)'_0,\nchi_{\tilde E}(\G_\cdot^y)'_0\big)
=\int_{E\Delta \tilde E}|(\G_x^y)'_0|_x^2\,\d\mu(x)\leq \mu(E\Delta \tilde E)\,D_\kappa^2
\]
grants that the closure of $\mathcal D$ contains all the sections of the form $\nchi_E(\G_\cdot^y)'_0$ for $E\subset\Y$ Borel and $y\in\mathcal D_1$. To conclude recall the continuity of $y\mapsto (\G_x^y)'_0\in \T_x\Y$ proved in Theorem \ref{thm:samecone} and notice that an application of the dominated convergence theorem gives that $\nchi_E(\G_\cdot^{y_n})'_0\to \nchi_E(\G_\cdot^y)'_0$ in $L^2(\T_G\Y;\mu)$ if $y_n\to y$.

\textsc{$\Cat{0}$ condition.} By \cite[Remark 2.2]{stu03}, it is sufficient
to show that for any $v,v'\in L^2(\T_G\Y;\mu)$ there exists $v''\in L^2(\T_G\Y;\mu)$
such that
\[
\sfd_\mu^2(w,v'')\leq\frac{1}{2}\,\sfd_\mu^2(w,v)+
\frac{1}{2}\,\sfd_\mu^2(w,v')-\frac{1}{4}\,\sfd_\mu^2(v,v')
\qquad\text{ for every }w\in L^2(\T_G\Y;\mu).
\]
Let us define $v''\in L^2(\T_G\Y;\mu)$ as
\[
v''_x:=\frac{1}{2}\,v_x\oplus v'_x\in\T_x\Y\qquad\text{ for }\mu\text{-a.e.\ }x\in\Y.
\]
(Note that $v''_x$ is the midpoint between $v_x$ and $v'_x$.) The fact
that $v''$ is (the equivalence class of) a Borel section of $\T_G\Y$ follows
by Proposition \ref{prop:boroper}, while the integrability condition
$(x\mapsto |v''_x|_x)\in L^2(\mu)$ is implied by inequality \eqref{eq:PI}.
By \cite[Corollary 2.5]{stu03}, we have
\[
\sfd_x^2(w_x,v''_x)\leq\frac{1}{2}\,\sfd_x^2(w_x,v_x)
+\frac{1}{2}\,\sfd_x^2(w_x,v'_x)-\frac{1}{4}\,\sfd_x^2(v_x,v'_x)
\qquad\text{ for }\mu\text{-a.e.\ }x\in\Y.
\]
By integrating with respect to $\mu$ we obtain the desired inequality.
\end{proof}

\section{Normal 1-currents and the superposition principle}
In this section we recall the notion of metric 1-current as introduced by Ambrosio-Kirchheim in \cite{AmbrosioKirchheim00}  and Paolini-Stepanov's metric version of Smirnov's superposition principle. Throughout this section $(\Y,\sfd)$ is a complete and separable metric space. See also \cite{lan11} and \cite{wen14} for more on the topic.

\bigskip

We denote by $\LIP(\Y)$ the space of real-valued Lipschitz functions on $\Y$, and by $\LIP_{b}(\Y)$ the subspace of bounded Lipschitz functions.
\begin{definition}[Normal 1-currents]
A (metric) $1$-current of finite mass on $\Y$ is a bilinear functional 
$$
T:\,\LIP_{b}(\Y)\times \LIP(\Y)\to \R
$$ 
satisfying  the following conditions:
\begin{itemize}
	\item[(a)] $T(g,f)=0$ if the function $f \in\LIP(\Y)$ is constant on the support of $g\in \LIP_{b}(\Y)$,
	\item[(b)] $T(g,f_n)\to T(g,f)$ whenever $f_n\to f$ pointwise  and $\sup_n\Lip(f_n)<\infty$,
	\item[(c)] there exists a finite Borel measure $\nu$ on $\Y$ satisfying 
	\begin{equation}
\label{eq:masscurr}
	\big|T(g,f)\big|\le \Lip(f) \int |g|\,\d\nu\qquad\forall g\in\LIP_{b}(\Y),\ f\in \LIP(\Y).
\end{equation}
\end{itemize}
A normal 1 current is a 1-current of finite mass such that there is a finite Borel measure $\mu$ (called boundary of $T$ and denoted by $\partial T$) such that
\[
T(1,f)=\int f\,\d\mu\qquad\forall f\in \LIP_b(\Y).
\]
\end{definition}
It is not hard to check that if $T$ has finite mass, there is a minimal (in the sense of partial ordering of measures) Borel measure for which \eqref{eq:masscurr} holds: it will be denoted by $\|T\|$ and called mass measure of $T$. We set $\M(T):=\|T\|(\Y)$.

A prototypical example is the normal 1-current $[\![\gamma]\!]$ induced by an absolutely continuous curve $\gamma:[0,1]\to\Y$ via the formula
\[
[\![\gamma]\!](g,f):=\int_0^1 g(\gamma_t) (f\circ\gamma)'_t\,\d t,\quad (g,f)\in \LIP_b(\Y)\times\LIP(\Y).
\]
Its mass measure is given by $\gamma_*\big(|\dot\gamma|\mathcal L^1\restr{[0,1]}\big)$ and its boundary is given by  $\partial[\![\gamma]\!]=\delta_{\gamma_1}-\delta_{\gamma_0}$.

Notice that the current $[\![\gamma]\!]$ remains unchanged if we change the parametrization of $\gamma$. This makes it natural to consider the space of `curves up to reparametrization' as follows (here we only consider non-decreasing reparametrizations).

\subsection{Reparametrizations of curves}
A reparametrization $\alpha:[0,1]\to[0,1]$ is a non-decreasing continuous surjection. If $\gamma,\eta\in C([0,1];\Y)$, we say that $\eta$ is a reparametrization of $\gamma$ if there is a reparametrization $\alpha$ satisfying  $\gamma\circ \alpha=\eta$.

\begin{remark}\label{regularparam}{\rm
Given $\gamma\in C([0,1];\Y)$, there exists a curve $\eta\in C([0,1];\Y)$ which is not constant on any open interval, and is a reparametrization of $\gamma$, cf.\ \cite[Proposition 3.6]{fah07}.
\fr}\end{remark}

Define an equivalence relation on $C([0,1];\Y)$ by declaring $\gamma\sim \eta$ if there is a curve $\varphi\in C([0,1];\Y)$ which is a reparametrization of both $\gamma$ and $\eta$. It is easy to see that this indeed defines an equivalence relation. Let $\overline\Gamma(\Y):=C([0,1],\Y)/\sim$ be the quotient space. We define a distance function on $\overline\Gamma(\Y)$ by
\[
\sfd_{\overline\Gamma}([\gamma],[\eta])=\inf\big\{\sfd_\infty(\gamma\circ\alpha,\eta\circ\beta):\ \alpha,\beta\textrm{ reparametrizations}\big\},\quad [\gamma],[\eta]\in \overline\Gamma(\Y),
\]
where
\[
\sfd_\infty(\gamma,\eta):=\sup_{0\le t\le 1}\sfd(\gamma_t,\eta_t).
\]
This is clearly symmetric, and satisfies the triangle inequality. Consequently it defines a pseudometric on $\overline\Gamma(\Y)$. It follows from Lemma \ref{lem:miniparam} below that $\sfd_{\overline\Gamma}$ defines a metric on $\overline \Gamma(\Y)$.

Note that, since a non-decreasing surjection $[0,1]\to [0,1]$ may be approximated uniformly by increasing homeomorphisms of $[0,1]$, it easily follows that $\sfd_{\overline\Gamma}([\gamma],[\eta])$ has the representation
\[
\sfd_{\overline\Gamma}([\gamma],[\eta])=\inf\big\{ \sfd_\infty(\gamma\circ\phi,\eta):\ \phi\textrm{ increasing homeomorphism of }[0,1]\big\}.
\]
\begin{lemma}\label{lem:miniparam}
	Let $\gamma,\eta\in C([0,1];\Y)$ be such that $\sfd_{\overline\Gamma}([\gamma],[\eta])=0$. Then $[\gamma]=[\eta]$.
\end{lemma}
\begin{proof}
	Let $\gamma,\eta\in C([0,1];\Y)$. By Remark \ref{regularparam} we may assume that $\gamma$ is not constant on any non-trivial interval. We will prove that there is a reparametrization $\phi$ such that $\gamma\circ\phi=\eta$. 
	
	Let $\phi_n:[0,1]\to[0,1]$ be a sequence of increasing homeomorphisms minimizing $\sfd_{\overline\Gamma}([\gamma],[\eta])$:
	\[
	\lim_{n\to\infty}\sfd_\infty(\gamma\circ\phi_n,\eta)=0.
	\]
	Denote $\psi_n=\phi_n\inv$. For each $n\in\N$, $\psi_n$ is also an increasing homeomorphism. Thus, $\phi_n$ and $\psi_n$ are of bounded variation and their distributional derivatives $\phi'_n$ and $\psi_n'$ (which are positive measures on $[0,1]$) satisfy
	\begin{align*}
	\int_0^1|\phi_n'|\,\d t=1=\int_0^1|\psi_n'|\,\d t
	\end{align*}
	for all $n\in\N$. By Helly's selection principle (see \cite{nat55}), there are subsequences (labeled here by the same indices) and functions $\phi,\psi:[0,1]\to [0,1]$ of bounded variation so that $\phi_n\to\phi$ and $\psi_n\to\psi$ pointwise. Clearly, $\phi$ and $\psi$ are non-decreasing and satisfy $\phi(0)=\psi(0)=0$, $\phi(1)=\psi(1)=1$. Since $\gamma$ is continuous, and $\phi_n\to \phi$ pointwise, we have the estimate
	\begin{align*}
	\sfd(\gamma_{\phi(t)},\eta_t)=\lim_{n\to\infty}\sfd(\gamma_{\phi_n(t)},\eta_t)\le\lim_{n\to\infty}\sfd_\infty(\gamma\circ\phi_n,\eta)=0
	\end{align*}
	for all $t\in[0,1]$. Thus
	\[
	\gamma\circ\phi=\eta.
	\]
	Similarly we obtain
	\[
	\gamma=\eta\circ\psi.
	\]
	
For any $0\le a<b\le 1$, the pointwise convergence $\phi_n\to \phi$ implies that 
\[
\bigcup_{k= 1}^\infty\bigcap_{n\ge k}\phi_n\inv[a,b]\subset \phi\inv[a,b]. 
\]
Moreover, we have
\[
(\psi(a),\psi(b))\subset \bigcup_{k=1}^\infty\bigcap_{n\ge k}\phi_n\inv[a,b].
\]
To see this, let $x\in (\psi(a),\psi(b))$. For all large enough $n\in\N$, we have $\psi_n(a)<x<\psi_n(b)$ or, equivalently, $a<\phi_n(x)<b$. Thus $\phi(x)=\lim_{n\to\infty}\phi_n(x)\in [a,b]$. The inclusions above imply
\begin{equation}\label{inclusion}
(\psi(a),\psi(b))\subset \phi\inv[a,b].
\end{equation}
It follows that $\phi$ is continuous. Indeed, if a non-decreasing function $\phi:[0,1]\to[0,1]$ has a point of discontinuity, it must omit some non-trivial interval $[a,b]\subset [0,1]$, i.e.\ $\phi\inv[a,b]=\varnothing$, implying $\psi(a)=\psi(b)$. By (\ref{inclusion}), we have
\[
\sfd(\gamma_a,\gamma_b)=\sfd(\eta_{\psi(a)},\eta_{\psi(b)})=0,
\]
which contradicts the fact that $\gamma$ is not constant on any non-trivial interval. 
\end{proof}

\begin{remark}{\rm
	Since $C([0,1];\Y)$ is complete and separable, we have that $\overline\Gamma(\Y)$ is complete and separable.
\fr}\end{remark}

\bigskip We will denote by $\Gamma(\Y)\subset \overline\Gamma(\Y)$ the image of $AC([0,1];\Y)\subset C([0,1];\Y)$ under the quotient map $q:\,C([0,1];\Y)\to\overline\Gamma(\Y)$
given by $\gamma\mapsto[\gamma]$, i.e.\ $\Gamma(\Y)=q\big(AC([0,1];\Y)\big)$.
Notice that since $AC([0,1];\Y)$ is a Borel subset of $C([0,1];\Y)$ (see for instance \cite[Section 2.2]{AmbrosioGigliSavare11}), we have that $\Gamma(\Y)$ is a Suslin subset of  $\overline\Gamma(\Y)$ and thus universally measurable.

Recall that a curve $\gamma\in C([0,1];\Y)$ is called rectifiable, if it has \emph{finite length}:
\[
\ell(\gamma):=\sup\left\{ \sum_{i=1}^m\sfd(\gamma_{t_i},\gamma_{t_{i-1}}) \right\}<\infty,
\]
where the supremum is taken over all partitions $0=t_0<\ldots<t_m=1$ of $[0,1]$. Note that the length $\ell(\gamma)$ is independent of reparametrization and, for absolutely continuous curves, is given by
\[
\ell(\gamma)=\int_0^1|\dot\gamma_t|\,\d t;
\]
see \cite{Haj03} for these statements, as well as the proposition below.

\begin{proposition}[Reparametrization with constant speed]\cite[Theorem 3.2 and Corollary 3.8]{Haj03}\label{prop:repar}
	Let $\Y$ be a complete and separable space and $\gamma\in C([0,1];\Y)$ a non-constant rectifiable curve. Define $\phi:[0,1]\to[0,1]$ by 
	\[
	\phi(t):=\inf\Big\{s\in[0,1]\ :\ \ell(\gamma\restr{[0,s]})=t\ell(\gamma)\Big\}, \quad t\in[0,1].
	\]
	Then $\bar\gamma:=\gamma\circ\phi$ is a $\ell(\gamma)$-Lipschitz curve and has constant metric speed
	$\displaystyle |\dot{\bar\gamma}_t|=\ell(\gamma)$ for almost every $t$. Moreover
	\[
	\bar\gamma\big(\ell(\gamma\restr{[0,t]})/\ell(\gamma)\big)=\gamma(t),\quad t\in[0,1],
	\] 
	i.e.\ $\bar\gamma$ is a reparametrization of $\gamma$.
	
	If $\gamma\sim\eta$ are two absolutely continuous curves, then their reparametrizations with constant speed coincide.
\end{proposition}

We shall denote by ${\sf ConstSpRep}:\Gamma(\Y)\to C([0,1],\Y)$ the map sending the equivalence class of $\gamma\in AC([0,1];\Y)$ to the constant speed reparametrization of any element in the class. Proposition \ref{prop:repar} implies that this map is well-defined. Also, we have:
\begin{proposition}\label{prop:repbor}
	Let $\Y$ be a complete separable space. Then ${\sf ConstSpRep}:\Gamma(\Y)\to C([0,1];\Y)$ is a Borel map.
\end{proposition}
\begin{proof}
	Throughout the proof we use the shorthand $\bar\theta:={\sf ConstSpRep}(\theta)$. Introduce a new metric $\sfd_0$ on $\Gamma(\Y)$ by setting
	\[
	\sfd_0(\theta_1,\theta_2):=\max\big\{\sfd_\Gamma(\theta_1,\theta_2),\big|\ell(\theta_1)-\ell(\theta_2)\big|\big\},\quad\theta_1,\theta_2\in\Gamma(\Y).
	\]
	Since the length functional $\ell:\Gamma(\Y)\to [0,\infty]$ is lower semicontinuous, and since
	\[
	B_{\sfd_0}(\theta,r)=B_\Gamma(\theta,r)\cap \ell^{-1}(\ell(\theta)-r,\ell(\theta)+r),
	\]
	it follows that $\sfd_0$-balls $B_{\sfd_0}(\theta,r)$ are Borel in $(\Gamma(\Y),\sfd_\Gamma)$. Consequently the identity
	\[
	I:(\Gamma(\Y),\sfd_\Gamma)\to(\Gamma(\Y),\sfd_0)
	\]
	is a Borel map.	Define
	\[
	h_0:(\Gamma(\Y),\sfd_0)\to C([0,1];Y),\quad \theta\mapsto \bar\theta
	\]
	and note that ${\sf ConstSpRep}=h_0\circ I$. Thus, it suffices to prove that $h_0$ is continuous. We thank Stefan Wenger for providing the elegant argument presented below.
	
	Suppose $\sfd_0(\theta_n,\theta)\to 0$ as $n\to\infty$. Then there are nondecreasing bijections $\varphi_n:[0,1]\to [0,1]$ such that
	\[
	\sfd_\infty(\theta_n\circ\varphi_n,\bar\theta)\to 0\quad\text{ as }n\to\infty.
	\]
	Denote $\gamma_n:=\theta_n\circ\varphi_n$. We have $\displaystyle \lim_{n\to\infty}\ell(\gamma_n)= \ell(\theta)$. Moreover, for any $t\in [0,1]$, we have
	\begin{align}\label{stefan1}
	\ell(\bar\theta)=\ell(\bar\theta\restr{[0,t]})+\ell(\bar\theta\restr{[t,1]})\le & \liminf_{n\to\infty}\ell(\gamma_n\restr{[0,t]})+\liminf_{n\to\infty}\ell(\gamma_n\restr{[t,1]})\\
	\le & \liminf_{n\to\infty}\big(\ell(\gamma_n\restr{[0,t]})+\ell(\gamma_n\restr{[t,1]})\big)=\ell(\bar\theta).\nonumber
	\end{align}
	Since
	\begin{align}\label{stefan2}
	\ell(\bar\theta\restr{[0,t]})\le \liminf_{n\to\infty}\ell(\gamma_n\restr{[0,t]})\quad\text{ and }\quad\ell(\bar\theta\restr{[t,1]})\le \liminf_{n\to\infty}\ell(\gamma_n\restr{[t,1]}),
	\end{align}
	it follows from \eqref{stefan1} that the inequalities in \eqref{stefan2} are in fact equalities, and we may pass to a subsequence (not relabeled) so that, for a countable dense set $D\subset [0,1]$, we have
	\begin{equation}\label{stefan3}
	\ell(\bar\theta\restr{[0,t]})=\lim_{n\to\infty}\ell(\gamma_n\restr{[0,t]}), \quad t\in D.
	\end{equation}
	Set
	\[
	\ell_n(t):=\frac{\ell(\gamma_n\restr{[0,t]})}{\ell(\gamma_n)},\quad t\in [0,1],
	\]
	whence, by \eqref{stefan3} and the fact that $\bar\theta$ is constant speed parametrized, we have
	\begin{equation}\label{stefan4}
	\lim_{n\to\infty}\ell_n(t)=\frac{\ell(\bar\theta\restr{[0,t]})}{\ell(\bar\theta)}=t,\quad t\in D.
	\end{equation}
	The sequence $(\bar\gamma_n)_n$ of constant speed parametrizations of $\gamma_n$ is uniformly Lipschitz and thus, after passing to a subsequence, it has a uniform limit $\beta:[0,1]\to\Y$ which is a Lipschitz curve. Note that $\bar\theta_n=\bar\gamma_n$.
	
	By the constant speed parametrization, we have
	\begin{equation}\label{stefan5}
	\bar\gamma_n\circ\ell_n=\gamma_n.
	\end{equation}
	For each $t\in D$ we have, by (\ref{stefan4}) and (\ref{stefan5}),
	\[
	\beta(t)=\lim_{n\to\infty}\bar\gamma_n(\ell_n(t))=\lim_{n\to\infty}\gamma_n(t)=\bar\theta(t).
	\]
	Since the equality holds on a dense set of points, we conclude that $\beta=\bar\theta$.
	
	By repeating this argument for any subsequence of $\theta_n$ we have that, if $\theta_n\to \theta$ in $\sfd_0$, then $\bar\theta_n\to\bar\theta$. Thus $h_0$ is continuous, and this completes the proof of the claim.
\end{proof}

\subsection{The superposition principle} We shall consider finite Borel measures $\pi$ on $\overline\Gamma(\Y)$  concentrated on $\Gamma(\Y)$ and typically denote by $[\gamma]$ their `integration variable'. In doing this, we always implicitly assume that $\gamma$ is absolutely continuous for $\pi$-a.e.\ $[\gamma]$ (i.e.\ we select an element in $[\gamma]$ which is absolutely continuous -- see also Proposition \ref{prop:repbor}  above).

\begin{lemma}\label{borel}
	For any $(g,f)\in\LIP_b(\Y)\times\LIP(\Y)$, the map
	$A:\,C([0,1];\Y)\to\R\cup\{+\infty\}$ given by
	\begin{equation}
	\label{eq:equivcurr}
	A(\gamma):=\left\{\begin{array}{ll}
	[\![\gamma]\!](g,f)\\
	+\infty
	\end{array}\quad\begin{array}{ll}
	\text{ if }\gamma\in AC([0,1];\Y),\\
	\text{ otherwise}
	\end{array}\right.
	\end{equation}
	is a Borel map.
\end{lemma}
\begin{proof}
	Since $AC([0,1];\Y)\subset C([0,1];\Y)$ is Borel, it suffices to show that the map
	\[
	A\restr{AC}:\,AC([0,1];\Y)\to\R,\quad \gamma\mapsto [\![\gamma]\!](g,f)
	\]
	is Borel. Let $\bar\gamma$ denote the constant speed parametrization given by Proposition \ref{prop:repar}. Let $q:\,C([0,1];\Y)\to\overline\Gamma(\Y)$ be the quotient map $\gamma\mapsto [\gamma]$. Since $\bar\gamma={\sf ConstSpRep}(q(\gamma))$, Proposition \ref{prop:repbor} implies that $AC([0,1];\Y)\ni\gamma\mapsto \bar\gamma$ is Borel. Consequently the map
	\[
	I_n:AC([0,1];\Y)\to\R,\quad \gamma\mapsto \int_0^{1-1/n}g(\bar\gamma_t)f(\bar\gamma_{t+1/n})\,\d t
	\]
	is Borel for each $n\in \N\cup\{\infty\}$. To show that $A\restr{AC}$ is Borel it suffices to see that
	\begin{equation}\label{eq:limitborel}
	A\restr{AC}=\lim_{n\to\infty}n(I_n-I_\infty).
	\end{equation}
	For each $\gamma\in AC([0,1];\Y)$ we have that $f\circ\bar\gamma$ is Lipschitz and thus, by the dominated convergence theorem,
	\begin{align*}
	A(\gamma)=[\![\gamma]\!](g,f)&=[\![\bar\gamma]\!](g,f)=\int_0^1g(\bar\gamma_t)\lim_{n\to\infty}n\big[f(\bar\gamma_{t+1/n})-f(\bar\gamma_t)\big]\,\d t\\
	&=\lim_{n\to\infty}n\left(\int_0^{1-1/n}g(\bar\gamma_t)\big[f(\bar\gamma_{t+1/n})-f(\bar\gamma_t)\big]\,\d t\right)=\lim_{n\to\infty}n\big(I_n(\gamma)-I_\infty(\gamma)\big),
	\end{align*}
	establishing \eqref{eq:limitborel}.
\end{proof}

By \eqref{eq:equivcurr}, and the fact that $[\![\gamma]\!]=[\![\eta]\!]$ if $\gamma\sim\eta$, we see that for any finite non-negative Borel measure $\pi$ on $\overline\Gamma(\Y)$ concentrated on $\Gamma(\Y)$ the functional $[\![\pi]\!]:\LIP_{b}(\Y)\times \LIP(\Y)\to \R$ given by
\[
[\![\pi]\!](g,f):=\int [\![\gamma]\!](g,f)\,\d\pi([\gamma])\qquad\forall (g,f)\in\LIP_{b}(\Y)\times \LIP(\Y)
\]
is well-defined and a normal 1-current: for its mass we have the bound
\begin{equation}
\label{eq:boundmass}
\int g\,\d\|[\![\pi]\!]\|\leq \iint g\,\d\gamma_*(|\dot\gamma|\mathcal L^1\restr{[0,1]})\,\d\pi([\gamma])=\iint_0^1g(\gamma_t)|\dot\gamma_t|\,\d t\,\d\pi([\gamma])
\end{equation}
for every  non-negative $ g\in \LIP_b(\Y)$; notice that $\gamma_*(|\dot\gamma|\mathcal L^1\restr{[0,1]})$ is independent on the parametrization of $\gamma$ -- see also Proposition \ref{prop:repar} below. For its boundary we have
\[
\int f\,\d\partial[\![\pi]\!]=\int f\,\d\big(({\rm e}_1)_*\pi-({\rm e}_0)_*\pi\big)=\int f(\gamma_1)-f(\gamma_0)\,\d\pi([\gamma])\qquad\forall f\in \LIP_b(\Y)
\]
(notice that $\gamma_0,\gamma_1$ are independent on the parametrization of $\gamma$). Observe that picking $g\equiv 1$ in \eqref{eq:boundmass} we obtain
\begin{equation}
\label{eq:boundmass2}
\M([\![\pi]\!])\leq\iint_0^1|\dot\gamma_t|\,\d t\,\d \pi([\gamma]).
\end{equation}
The superposition principle states that \emph{every} normal 1-current is of the form $[\![\pi]\!]$ for some $\pi$ as above, and moreover $\pi$ can be chosen so that equality holds in \eqref{eq:boundmass2}. For the proof of the following result we refer to \cite[Corollary 3.3]{PaolStep13}:
\begin{theorem}[Superposition principle]\label{thm:sup1}
Let $\Y$ be a complete and separable space and $T$ a normal 1-current. Then there is a finite non-negative Borel measure $\pi$ on $\overline\Gamma(\Y)$ concentrated on $\Gamma(\Y)$ such that
\[
\begin{split}
T&=[\![\pi]\!],\\
\M(T)&=\iint_0^1|\dot\gamma_t|\,\d t\,\d \pi([\gamma]).
\end{split}
\] 
\end{theorem}
For our applications it will be more convenient to deal with measures on $C([0,1];\Y)$ rather than on $\overline\Gamma(\Y)$. Using Lemma \ref{borel} and Proposition \ref{prop:repbor}, we can reformulate Theorem \ref{thm:sup1} as follows:
\begin{theorem}[Superposition principle - equivalent formulation]\label{thm:sup2}
Let $\Y$ be a complete and separable space and $T$ a normal 1-current. Then there is a finite non-negative Borel measure $\pi$ on $C([0,1];\Y)$ concentrated on the set of non-constant absolutely continuous curves with constant speed such that
\begin{equation}
\label{eq:sup}
\begin{split}
T(g,f)&=\iint_0^1  g(\gamma_t)(f\circ\gamma)'_t\,\d t\,\d\pi(\gamma),\\
\int g \,\d\|T\|&=\iint_0^1g(\gamma_t)|\dot\gamma_t|\,\d t\,\d\pi(\gamma)
\end{split}
\end{equation}
for any $g\in\LIP_b(\Y)$ and $f\in\LIP(\Y)$.
\end{theorem}
\begin{proof}
By Theorem \ref{thm:sup1}, there is a finite measure $\overline\eta \in \mathscr M(\Gamma(\Y))$ for which
\begin{align*}
T(g,f)&=\int\!\!\!\int_0^1g(\theta_t)(f\circ\theta)'_t\,\d t\, \d\overline\eta(\theta)\quad\textrm{ for all }(g,f)\in \LIP_b(\Y)\times\LIP(\Y),\\
\mathbb M(T)&=\int \ell(\theta)\,\d\overline \eta(\theta).
\end{align*}
We define $\pi\in\mathscr M(C([0,1];\Y))$ as \[ \pi:={\sf ConstSpRep}_\ast\overline\eta. \] Since both $$\ell(\theta)=\int_0^1|\dot\theta_t|\,\d t\quad\text{ and }\quad\int_0^1g(\theta_t)(f\circ\theta)'_t\,\d t$$ are independent of parametrization, we have the identities
\begin{align}
T(g,f)&=\iint_0^1g(\gamma_t)(f\circ\gamma)'_t\,\d t\,\d\pi(\gamma),\label{ident1}\\
\mathbb M(T)&=\iint_0^1|\dot\gamma_t|\,\d t\,\d\pi(\gamma)
\end{align}
for all $(g,f)\in\LIP_b(\Y)\times\LIP(\Y)$.

It remains to prove the second identity in the claim. It suffices to prove it for $g=\nchi_E$ for Borel sets $E\subset \Y$. 
It follows from \eqref{ident1} that \[\big|T(g,f)\big|\le\iint_0^1|g|(\gamma_t)\lip_af(\gamma_t)|\dot\gamma_t|\,\d t\,\d\pi(\gamma), \] whence $\|T\|\le \nu$, where $\nu$ is defined by \[ \nu(E):=\iint_0^1\nchi_E(\gamma_t)|\dot\gamma_t|\,\d t\,\d\pi(\gamma),\quad E\subset \Y\textrm{ Borel}. \]

By the characterisation of mass (see \cite[Proposition 2.7]{AmbrosioKirchheim00}) it follows that, for every $\eps>0$, there are functions $(g_\eps,f_\eps)\in \LIP_b(\Y)\times\LIP(\Y)$ such that $|g_\eps|\le 1$ and $\Lip f_\eps\le 1$, and for which
\[
\M(T)-\eps<T(g_\eps,f_\eps).
\]
Using \eqref{ident1} and the identity $1=\nchi_E(\gamma_t)+\nchi_{\Y\setminus E}(\gamma_t)$, we have
\begin{align*}
&\iint_0^1\nchi_E(\gamma_t)|\dot\gamma|_t\,\d t\,\d\pi(\gamma)+\iint_0^1\nchi_{\Y\setminus E}(\gamma_t)|\dot\gamma|_t\,\d t\,\d\pi(\gamma)-\eps=\M(T)-\eps<T(g_\eps,f_\eps)\\
=&\iint_0^1\nchi_E(\gamma_t)g_\eps(\gamma_t)(f_\eps\circ\gamma)'_t\,\d t\,\d\pi(\gamma)+\iint_0^1\nchi_{\Y\setminus E}(\gamma_t)g_\eps(\gamma_t)(f_\eps\circ\gamma)'_t\,\d t\,\d\pi(\gamma)\\
\le & \iint_0^1\nchi_E(\gamma_t)g_\eps(\gamma_t)(f_\eps\circ\gamma)'_t\,\d t\,\d\pi(\gamma)+\iint_0^1\nchi_{\Y\setminus E}(\gamma_t)|\dot\gamma_t|\,\d t\,\d\pi(\gamma),
\end{align*}
which implies
\[
\iint_0^1\nchi_E(\gamma_t)|\dot\gamma|_t\,\d t\,\d\pi(\gamma)-\eps\le \iint_0^1\nchi_E(\gamma_t)g_\eps(\gamma_t)(f_\eps\circ\gamma)'_t\,\d t\,\d\pi(\gamma)\le \|T\|(E)
\]
for every $\eps>0$. It follows that $\|T\|=\nu$, and this completes the proof of the last identity in (\ref{eq:sup}).
\end{proof}

\section{Metric measure spaces}\label{se:mms}
For our purposes, a metric measure space is a triple $(\Y,\sfd,\mu)$ where $(\Y,\sfd)$
is a complete separable metric space and $\mu$ a Borel measure on $\Y$ that is finite
on bounded sets.

\subsection{Derivations and the space $\Dertwotwo$}
We introduce \emph{derivations} and their basic properties,
based on the presentation in \cite{DiM14a,DiMarino14}.
This notion of derivation has been inspired by a similar concept
introduced by N.\ Weaver in \cite{Weaver01}.
\medskip

Let us denote by $L^0(\mu)$ the set of equivalence classes of $\mu$-measurable maps on $\Y$ (without any integrability assumptions).

\begin{definition}\label{der}
	A derivation $b$ on $\Y$ is a linear map $b:\,\LIP_b(\Y)\to L^0(\mu)$ satisfying the following two conditions:
	\begin{itemize}
		\item[(1)](Leibniz rule) $b(fh)=fb(h)+hb(f)$ $\mu$-a.e.\ for
		all $f,h\in \LIP_b(\Y)$. 
		\item[(2)](Weak locality) There is $g\in L^0(\mu)$ such that $\big|b(f)\big|\le g\,\lip_af$ $\mu$-a.e.\ for all $f\in \LIP_b(\Y)$.
	\end{itemize}
\end{definition}
We denote the set of derivations on $\Y$ by ${\rm Der}(\Y)$.
The space ${\rm Der}(\Y)$ is a $\LIP_b(\Y)$-module: given a Lipschitz function
$\varphi\in\LIP_b(\Y)$ and a derivation $b\in{\rm Der}(\Y)$, the linear map
\[
\varphi b:\,\LIP_b(\Y)\to L^0(\mu),\qquad f\mapsto\varphi b(f)
\]
is again a derivation; see \cite{DiM14a}.
\begin{remark}\label{derlip}{\rm
	By weak locality, we may extend a derivation $b\in {\rm Der}(\Y)$ to act on $\LIP(\Y)$. Indeed, given $f\in \LIP(\Y)$ and an open ball $B\subset \Y$, we have
	\[
	\nchi_Bb(f)=\nchi_B b(\tilde f)
	\]
	for any $\tilde f\in \LIP_b(\Y)$ for which $f\restr B=\tilde f\restr B$. Thus, for any $f\in \LIP(\Y)$ (and some fixed $x_0\in \Y$), the function 
	\[
	b(f)=\lim_{n\to\infty}\nchi_{B_n(x_0)}b\big((1-\dist(\cdot,B_n(x_0)))_+f\big)
	\]
	is well-defined, and $\LIP(\Y)\ni f\mapsto b(f)$ satisfies (1) and (2) above.
\fr}\end{remark}
\medskip

Given a derivation $b\in{\rm Der}(\Y)$, we define
\[ |b|=:{\rm ess\,sup}\big\{b(f)\;\big|\;f\in\LIP_b(\Y),\,\Lip(f)\le 1\big\}. \]
\begin{lemma}\label{ptnorm}
	Let $b\in{\rm Der}(\Y)$. Then $|b|$ satisfies (2) in Definition \ref{der}.
	Moreover, $|b|$ is the least function satisfying (2) in Definition \ref{der}.
\end{lemma}
\begin{proof}
	Let $f\in \LIP_b(\Y)$. For $x\in\Y$ and $r>0$, set $L_r=L_r(x):=\Lip\big(f\restr{B_r(x)}\big)$. Consider the McShane extension $g_r$  of $f\restr{B_r(x)}$; in particular, $g_r/L_r$ is $1$-Lipschitz and so we have
	$$\big|b(g_r/L_r)\big| \leq |b| \qquad \mu\text{-almost everywhere.}$$
Since  $b(g_r)=b(f)$ in $B_r(x)$, we deduce that $\big|b(f/L_r)\big|\le |b| $ holds $\mu$-almost everywhere on $B_r(x)$. Thus for each $x\in\Y$ and $r>0$ we have \[\big|b(f)\big| (y)\le |b|(y) \cdot L_r(x)\le |b|(y) \cdot L_{2r}(y)\qquad\mu\text{-a.e.\ }y\in B_r(x). \] Using this reasoning for a countable dense set $(x_n)\subset\Y$, we deduce that for every $r>0$ \[\big|b(f)\big|(y)\le |b|(y) \cdot  L_r(y) \qquad \mu\text{-a.e. }y\in\Y.\]
 The conclusion now follows by taking a sequence $r_n\downarrow 0$ and taking the limit as $n\to \infty$.
\end{proof}

A derivation $b\in{\rm Der}(\Y)$ is said to \emph{have divergence} if there exists a function $h\in L^1_b(\mu)$ (that is, $h$ is integrable on bounded sets) so that
\begin{equation}\label{div}
\int b(f)\,\d\mu=-\int f h\,\d\mu\qquad\text{ for all }f\in \LIP_b(\Y)
\end{equation}
(whenever this makes sense).
If such a function $h$ exists, it is unique and we denote it by $\div(b)$ or $\div\,b$.
The set of $b\in{\rm Der}(\Y)$ that have divergence is denoted by $D(\div)$.  
\medskip

For $1\le p\le \infty$ we set \[{\rm Der}^p_b(\Y;\mu):=\big\{b\in{\rm Der(\Y)}\,:\,|b|\in L^p_b(\mu)\big\}; \qquad{\rm Der}^p(\Y;\mu):=\big\{b\in{\rm Der}(\Y)\,:\,|b|\in L^p(\mu)\big\} \]
and, for $1\le p,q<\infty$,
\begin{align*}
{\rm Der}^{p,q}(\Y;\mu)&:=\big\{b\in{\rm Der}^p(\Y)\cap D(\div)\,:\,\div\,b\in L^q(\mu)\big\}.
\end{align*}

\begin{lemma}\label{weakco}
	Let $b\in D(\div)$. Assume $(f_n)$ is a sequence in $\LIP_b(\Y)$ converging to $f\in\LIP_b(\Y)$ pointwise and with $\sup_n\Lip(f_n)<\infty$.
	\begin{itemize}
		\item[(1)] Then
		\begin{equation}\label{eq:weakco}
		\int\varphi b(f_n)\,\d\mu\to\int\varphi b(f)\,\d\mu
		\end{equation}
		for each $\varphi\in\LIP_b(\Y)$.
		
		\item[(2)] If, in addition, $b\in{\rm Der}^p_b(\Y)$ for some $1<p<\infty$, then the convergence (\ref{eq:weakco}) holds for all $\varphi\in L^q(\mu)$ with bounded support. Here $q$ is the conjugate exponent of $p$, i.e.\ $1/p+1/q=1$.
	\end{itemize}
\end{lemma}
\begin{proof}
	By linearity it suffices to prove the claims when $f=0$. The Leibniz rule implies $$\varphi b(f_n)=b(\varphi f_n)-f_nb(\varphi).$$ Thus
	\begin{align*}
	\int\varphi b(f_n)\,\d\mu&=\int b(\varphi f_n)\,\d\mu-\int f_n b(\varphi)\,\d\mu=-\int f_n\big[\varphi\div\,b+b(\varphi)\big]\,\d\mu.
	\end{align*}
	Since $f_n\to 0$ pointwise and $\sup_n\Lip(f_n)<\infty$ it follows -- using the dominated convergence theorem -- that $\int\varphi b(f_n)\,\d\mu\to 0$ for all $\varphi\in \LIP_b(\Y)$. This proves (1).
	
	Let $\varphi\in L^q(\mu)$ have bounded support $B'\subset\Y$, and consider the set
	$B=\big\{x\,:\,\dist(B',x)\le 1\big\}$. Take a sequence $(\varphi_m)\subset \LIP_b(\Y)$ with supports in $B$ such that $$\lim_{m\to \infty}\int|\varphi_m-\varphi|^q\,\d\mu=0.$$ Denote $$L=\sup_n\Lip(f_n).$$ Then, for each $m,n\in \N$ we may estimate
	\begin{align*}
	\left|\int\varphi b(f_n)\,\d\mu\right|\le &\left|\int\varphi_m b(f_n)\,\d\mu\right|+\left|\int(\varphi-\varphi_m) b(f_n)\,\d\mu\right|\\
	\le & \left|\int\varphi_m b(f_n)\,\d\mu\right|+\left(\int|\varphi_m-\varphi|^q\,\d\mu\right)^{1/q}\left(\int_B\big|b(f_n)\big|^p\,\d\mu\right)^{1/p}\\
	\le & \left|\int\varphi_m b(f_n)\,\d\mu\right|+L\left(\int|\varphi_m-\varphi|^q\,\d\mu\right)^{1/q}\left(\int_B|b|^p\,\d\mu\right)^{1/p}
	\end{align*}
	Taking first $\limsup_{n\to\infty}$ and then $\limsup_{m\to\infty}$ we obtain $$\lim_{n\to\infty}\int_B\varphi b(f_n)\,\d\mu=0,$$
thus proving (2).
\end{proof}
In order to prove the next proposition, we recall the notion of \emph{strong locality,} cf.\ \cite[Lemma 7.13]{DiM14a}: if $b\in D(\div)$, then for every $f,g\in \LIP_b(\Y)$ we have \[
b(f)=b(g)\qquad\mu\text{-almost everywhere on }\{ f=g \}
\]
and, moreover,
\[
\big|b(f)\big|\le |b|\,\lip_a(f\restr C)\qquad\mu\text{-almost everywhere on }C,
\]
for every closed set $C\subset\Y$.

\begin{proposition}\label{prop:bnorm} Let $(\Y,\sfd,\mu)$ be a metric measure space, $\Omega\subset\Y$ an open set, and $b \in{\rm Der}^1_b(\Y)\cap D(\div)$. Let $\mathcal D=(x_n)\subset\Y$ be countable and dense in $\Omega$. Define $f_n(x):=\sfd(x_n,x)$. Then we have
	$${\rm ess\,sup}_n\big\{b(f_n)\big\}=|b|,\qquad{\rm ess\,inf}_n\big\{b(f_n)\big\}=-|b|$$
$\mu$-almost everywhere in $\Omega$.	
	\end{proposition}
\begin{proof}
	Denote $h_+={\rm ess\,sup}_nb(f_n)$ and $h_-={\rm ess\,inf}_nb(f_n)$. It suffices to prove that $|b|\le h_+$ and $-|b|\ge h_-$ $\mu$-almost everywhere on $\Omega$.\\
	
	{\bf Claim}\emph{
	Consider the countable set \[ \mathscr A=\big\{(\dist_{x_1}+q_1)\wedge\cdots\wedge(\dist_{x_k}+q_k)\;\big|\;x_1,\ldots,x_k\in\mathcal D,\,q_1,\ldots,q_k\in\Q,\,k\in\N\big\}.\] The set of restrictions $\{ g\restr\Omega\,:\,g\in\mathscr A \}$ is dense in $\LIP_1(\Omega):=\big\{f\in\LIP(\Omega)\,:\,\Lip(f)\leq 1\big\}$ in the topology of pointwise convergence.
	} 
\begin{proof}[Proof of Claim]
	Let $f\in \LIP_1(\Y)$. Since $\mathcal D$ is dense in $\Omega$, it is easy to see that
	\[
	f(x)=\inf\big\{g(x)\,:\,g\in\mathscr A,\,g\ge f\big\}
	\]
	for every $x\in\Omega$. For each $x_k\in \mathcal D$, let $(g_k^m)_m$ be a sequence in $\mathscr A$ satisfying \[ g_k^m(x_k)-f(x_k)<1/m \] for all $m\in\N$. Set
	\[
	g_n=g_1^n\wedge\cdots\wedge g_n^n.
	\]
	Then $(g_n)$ is a sequence in $\mathscr  A$, and $$\lim_{n\to \infty}g_n(x_k)= f(x_k)$$ for every $x_k\in\mathcal D$. Since $g_n$ and $f$ are 1-Lipschitz functions, it follows that $$\lim_{n\to \infty}g_n(x)= f(x)$$ for every $x\in\Omega$.
\end{proof}
	
	For any $f\in \LIP_1(\Y)$, let $(g_j)\subset \mathscr A$ be a sequence such that $g_j\restr{\Omega}$ converges to $f\restr\Omega$ pointwise. By passing to a subsequence we may assume that $g_j$ converges pointwise to some 1-Lipschitz function $f'$ in $\Y$ (in this case $f'\restr\Omega=f\restr\Omega$). For each $j$ write $g_j$ as $$g_j=f_1^j\wedge\cdots\wedge f_{k_j}^j,$$ with $f_n^j=\dist_{x_n^j}+q_n^j$. Define the sets $B^j_n$, $n=1,2,\ldots$ as $B^j_n:=\{ g_j=f _n^j \} $ if $1\le n\le k_j$ and $B^j_n:=\varnothing$ if $n>k_j$; also set $$C_n^j:=B_n^j\setminus \bigcup_{m<n}B_m^j.$$ Note that $\{C_n^j\}_n$ is a partition of $\Y$. By the strong locality of $b$ we have \[ b(g_j)=b(f_n^j) \] $\mu$-almost everywhere on $B_n^j$. Thus the identity \[ b(g_j)=\sum_n\chi_{C_n^j}b(f_n^j) \] is valid $\mu$-almost everywhere. For any non-negative $\eta\in \LIP_b(\Y)$ with bounded support, we then have
	\begin{align*}
	\int\eta b(g_j)\,\d\mu=\sum_n\int\eta\nchi_{C_n^j}b(f_n^j)\,\d\mu\leq
	\sum_n\int\eta\nchi_{C_n^j}h_+\,\d\mu=\int\eta h_+\,\d\mu.
	\end{align*}
	It follows that $b(g_j)\le h_+$ $\mu$-a.e.\ and, by Lemma \ref{weakco},
	that $b(f')\le h_+$ $\mu$-almost everywhere. Since $f'\restr\Omega=f\restr\Omega$ the locality of $b$ implies that $b(f)\le h_+$ $\mu$-almost everywhere on $\Omega$. Since $f$ is arbitrary it follows that $|b|\le h_+$ $\mu$-almost everywhere on $\Omega$. 
	
	The inequality $-|b|\ge h_-$ ($\mu$-almost everywhere) on $\Omega$ is proven analogously, using the identity $$-|b|={\rm ess\,inf}_{\Lip(f)\le 1}b(f).$$
\end{proof}

\begin{definition}\label{def:dermod}
	Given $p\ge 1$, we define the norm $\|\cdot\|_{p,p}$ on ${\rm Der}^{p,p}(\Y;\mu)$ as
	
	$$\|b\|_{p,p}:=\left(\int|b|^p\,\d\mu+\int(\div\ b)^p\,\d\mu\right)^{1/p}.$$
The normed space $\big({\rm Der}^{p,p}(\Y;\mu),\|\cdot\|_{p,p}\big)$ is a Banach space; see \cite{DiM14a}. We shall also use the norm
$$\|b\|_p:=\left(\int|b|^p\,\d\mu\right)^{1/p}.$$
\end{definition}

\subsection{The Sobolev space $W^{1,2}(\Y,\sfd,\mu)$}
In order to define Sobolev spaces on metric measure spaces, we adopt the approach in \cite{DiM14a} using derivations with divergence.
\begin{definition}[Sobolev space]\label{def:Sobolev}
Let $(\Y,\sfd,\mu)$ be a metric measure space and $p\in(1,\infty)$.
Let $q$ be the conjugate exponent of $p$. A function $f\in L^p(\mu)$ belongs to the Sobolev space $W^{1,p}(\Y,\sfd,\mu)$
provided there exists a $\LIP_b(\Y)$-linear continuous map $L_f:\,{\rm Der}^{q,q}(\Y;\mu)\to L^1(\mu)$ such that
\begin{equation}
\int L_f(b)\,\d\mu=-\int f\,\div\,b\,\d\mu\qquad\text{ for every }b\in{\rm Der}^{q,q}(\Y;\mu).
\end{equation}
\end{definition}
Whenever such a map $L_f$ exists, it is unique (cf.\ \cite[Remark 7.1.5]{DiM14a}).
\begin{theorem}[$p$-weak gradient]
Let $f\in W^{1,p}(\Y,\sfd,\mu)$. Then there is a function
$g_f\in L^p(\mu)$ such that
\begin{equation}\label{eq:def_pwg}
\big|L_f(b)\big|\leq g_f\,|b|\;\;\;\mu\text{-a.e.}
\qquad\text{ for every }b\in{\rm Der}^{q,q}(\Y;\mu).
\end{equation}
The least function $g_f$ (in the $\mu$-a.e.\ sense) that realises \eqref{eq:def_pwg}
is called \emph{$p$-weak gradient} of $f$ and denoted by $|Df|$.
\end{theorem}
For a proof of the previous result we refer to \cite[Theorem 7.1.6]{DiM14a}.
We point out that the $p$-weak gradient $|Df|$ might depend on $p$
(this dependence is omitted in our notation). Thus, the $p$-weak gradient and the $p'$-weak gradient of a function in $W^{1,p}(\Y,\sfd,\mu)\cap W^{1,p'}(\Y,\sfd,\mu)$ can be different.
\medskip

The space $W^{1,2}(\Y,\sfd,\mu)$ equipped with the norm
\begin{equation}
{\|f\|}_{W^{1,2}(\Y,\sfd,\mu)}:=
\left(\int|f|^2\,\d\mu+\int|Df|^2\,\d\mu\right)^{1/2}
\end{equation}
is a Banach space. In general it is not a Hilbert space.
There are alternative (equivalent) ways to define Sobolev spaces on metric measure spaces,
namely the approaches that have been proposed in \cite{Cheeger00,Shanmugalingam00,AmbrosioGigliSavare11};
see also \cite{HKST15,Bjorn-Bjorn11} and the monographs \cite{Heinonen01, Haj03}
for related discussions.
\medskip

By combining \cite[Theorem 7.2.5]{DiM14a} with the results of \cite{AmbrosioGigliSavare11-3},
one gets the ensuing approximation theorem: 
\begin{theorem}\label{thm:relax}
	Let $f\in W^{1,2}(\Y,\sfd,\mu)$ be given. Then there exists a sequence $(f_n)\subset \LIP_b(\Y)$ such that $f_n\to f$ and $\lip_af_n\to|Df|$ in $L^2(\mu)$.
\end{theorem}

The following identity expresses a duality between $W^{1,2}(\Y,\sfd,\mu)$
and $\Dertwotwo$.
\begin{proposition}\label{prop:duality}
	Let $f\in W^{1,2}(\Y,\sfd,\mu)$ be a given Sobolev function.
Let us denote by $\mathbb B$ the normed dual of
$\big(\Dertwotwo,\|\cdot\|_2\big)$.
	We define the element
	$\mathscr L_f\in\mathbb B$ as $\mathscr L_f(b):=\int L_f(b)\,\d\mu$
	for every $b\in\Dertwotwo$. Then
	\begin{equation}\label{eq:Int_isom}
	{\|\mathscr L_f\|}_{\mathbb B}={\big\||Df|\big\|}_{L^2(\mu)}.
	\end{equation}
\end{proposition}
To prove Proposition \ref{prop:duality}, we use the following well-known lemma. Let ${\rm LIP}_{bs}(\Y)$ be the space of all Lipschitz functions on $\Y$ with bounded support.
\begin{lemma}\label{lem:approx}
	Let $f\in L^\infty(\mu)$ be given. Then there exists a sequence
	$(f_n)_n\subseteq{\rm LIP}_{bs}(\Y)$ such that
	$f_n\to f$ pointwise $\mu$-a.e.\ and ${\|f_n\|}_{L^\infty(\mu)}\leq{\|f\|}_{L^\infty(\mu)}$
	for every $n\in\N$.
\end{lemma}

\begin{proof}[Proof of Proposition \ref{prop:duality}]
	\textsc{Step 1.} First of all, we claim that
	\begin{equation}\label{eq:Int_isom_claim}
	{\|\mathscr L_f\|}_{\mathbb B}=
	\sup\Big\{\int\big|L_f(b)\big|\,{\rm d}\mu\;\Big|
	\;b\in\Dertwotwo,\,\|b\|_2\leq 1\Big\}.
	\end{equation}
	Call $C$ the right hand side of \eqref{eq:Int_isom_claim}.
	Recall that by definition of dual norm we have
	\[{\|\mathscr L_f\|}_{\mathbb B}=
	\sup\Big\{\int L_f(b)\,{\rm d}\mu\;\Big|
	\;b\in\Dertwotwo,\,\|b\|_2\leq 1\Big\},\]
	whence trivially ${\|\mathscr L_f\|}_{\mathbb B}\leq C$. To show the converse
	inequality, fix $b\in\Dertwotwo$ with $\|b\|_2\leq 1$.
	By Lemma \ref{lem:approx}, we can choose $(g_n)_n\subseteq{\rm LIP}_{bs}(\Y)$
	such that $\sup_n|g_n|\leq 1$ and $g_n\to{\rm sgn}\,L_f(b)$ hold $\mu$-a.e..
	Hence by applying the dominated convergence theorem we get
	\[\int L_f(g_n b)\,{\rm d}\mu=\int g_n\,L_f(b)\,{\rm d}\mu
	\longrightarrow\int\big({\rm sgn}\,L_f(b)\big)L_f(b)\,{\rm d}\mu=
	\int\big|L_f(b)\big|\,{\rm d}\mu.\]
	Since $g_n b\in\Dertwotwo$ and $\|g_n b\|_2\leq\|b\|_2\leq 1$
	for all $n\in\N$, we have
	$\int\big|L_f(b)\big|\,{\rm d}\mu\leq{\|\mathscr L_f\|}_{\mathbb B}$
	and accordingly $C\leq{\|\mathscr L_f\|}_{\mathbb B}$.
	This proves \eqref{eq:Int_isom_claim}.\\
	\textsc{Step 2.} It can be readily checked that
	\[|Df|=\underset{b\in\Dertwotwo}{\rm ess\,sup\,}
	\nchi_{\{|b|>0\}}\,\frac{L_f(b)}{|b|}\quad\text{ in the }\mu\text{-a.e.\ sense.}\]
	This means that there exists a sequence
	$(b_i)_i\subseteq\Dertwotwo$ such that
	\[|Df|=\sup_{i\in\N}\,\nchi_{\{|b_i|>0\}}\frac{L_f(b_i)}{|b_i|}
	\quad\text{ in the }\mu\text{-a.e.\ sense.}\]
	For any $n\in\N$, we can pick pairwise disjoint Borel subsets $A^n_1,\ldots,A^n_n$
	of $\Y$ such that $|b_i|>0$ $\mu$-a.e.\ on $A^n_i$ for all $i\leq n$ and
	\[\sup_{i\leq n}\,\nchi_{\{|b_i|>0\}}\frac{L_f(b_i)}{|b_i|}=
	\sum_{i=1}^n\nchi_{A^n_i}\frac{L_f(b_i)}{|b_i|}\quad\text{ in the }\mu\text{-a.e.\ sense}.\]
	Notice that $\lim_n\mu\big(D\setminus\bigcup_{i\leq n}A^n_i\big)=0$, where
	we set $D:=\big\{|Df|>0\big\}$. Moreover, by monotone convergence theorem we see that
	$\sum_{i=1}^n\nchi_{A^n_i}\,L_f(b_i)/|b_i|\to|Df|$ in $L^2(\mu)$ as $n\to\infty$.
	Let us choose Borel subsets $B^n_i\subseteq A^n_i$ such that
	\begin{equation}\label{eq:Int_isom_claim2}\begin{split}
	&\nchi_{B^n_i}/|b_i|\in L^\infty(\mu)\quad\text{ for every }i\leq n,\\
	&\mu\big({\textstyle\bigcup_{i\leq n}A^n_i\setminus B^n_i}\big)\leq 1/n
	\quad\text{ for every }i\leq n,\\
	&\sum_{i=1}^n\nchi_{B^n_i}\,\frac{L_f(b_i)}{|b_i|}\to|Df|\;\;\text{ in }L^1(\mu)
	\quad\text{ as }n\to\infty.
	\end{split}\end{equation}
	Now let $n\in\N$ be fixed. Lemma \ref{lem:approx} grants for all $i\leq n$ the existence of
	$(g^i_k)_k\subseteq{\rm LIP}_{bs}(\Y)$ such that
	$\sup_k{\|g^i_k\|}_{L^\infty(\mu)}<+\infty$ and $g^i_k\to\nchi_{B^n_i}/|b_i|$ $\mu$-a.e.\ in
	$\Y$. Then an application of the dominated convergence theorem yields
	\[\begin{split}
	\sum_{i=1}^n g^i_k\,L_f(b_i)\overset{k}\longrightarrow
	\sum_{i=1}^n\nchi_{B^n_i}\frac{L_f(b_i)}{|b_i|}&\quad\text{ in }L^1(\mu),\\
	\bigg|\sum_{i=1}^n g^i_k\,b_i\bigg|\overset{k}\longrightarrow\nchi_{\bigcup_{i\leq n}B^n_i}&
	\quad\text{ in }L^2(\mu).
	\end{split}\]
	Hence for $k$ sufficiently big we have that the derivation
	$\tilde b_n:=\sum_{i=1}^n g^i_k\,b_i\in\Dertwotwo$
	is such that the $L^1(\mu)$-norm of
	$L(\tilde b_n)-\sum_{i=1}^n\nchi_{B^n_i}\,L_f(b_i)/|b_i|$
	and the $L^2(\mu)$-norm of $|\tilde b_n|-\nchi_{\bigcup_{i\leq n}B^n_i}$ are
	smaller than $1/n$. By recalling \eqref{eq:Int_isom_claim2}, we thus deduce
	that $L(\tilde b_n)\to|Df|$ in $L^1(\mu)$ and $|\tilde b_n|\to\nchi_D$ in $L^2(\mu)$
	as $n\to\infty$. Possibly passing to a not relabeled subsequence, we can assume
	that there exists $G\in L^1(\mu)$ such that
	\begin{equation}\label{eq:Int_isom_claim3}\begin{split}
	\big|L_f(\tilde b_n)\big|,|\tilde b_n|^2\leq G
	&\;\;\;\mu\text{-a.e.}\quad\text{ for every }n\in\N,\\
	L_f(\tilde b_n)\to|Df|
	&\;\;\;\mu\text{-a.e.}\quad\text{ as }n\to\infty,\\
	|\tilde b_n|\to\nchi_D
	&\;\;\;\mu\text{-a.e.}\quad\text{ as }n\to\infty.
	\end{split}\end{equation}
	\textsc{Step 3.} We can finally prove \eqref{eq:Int_isom}.
	For any $b\in\Dertwotwo$ with $\|b\|_2\leq 1$ it holds that
	\[\int\big|L_f(b)\big|\,{\rm d}\mu\leq\int|Df||b|\,{\rm d}\mu
	\leq{\big\||Df|\big\|}_{L^2(\mu)}\]
	by H\"{o}lder inequality, whence ${\|\mathscr L_f\|}_{\mathbb B}
	\leq{\big\||Df|\big\|}_{L^2(\mu)}$ by \eqref{eq:Int_isom_claim}. For the
	converse inequality, fix $h\in{\rm LIP}_{bs}(\Y)$.
	By recalling \eqref{eq:Int_isom_claim3} and using the dominated convergence theorem, we get
	\begin{equation}\label{eq:Int_isom_claim4}\begin{split}
	\int|h||Df|\,{\rm d}\mu&=\lim_{n\to\infty}\int|h|\big|L_f(\tilde b_n)\big|\,{\rm d}\mu
	=\lim_{n\to\infty}\int\big|L_f(h\tilde b_n)\big|\,{\rm d}\mu
	\overset{\eqref{eq:Int_isom_claim}}\leq\lim_{n\to\infty}{\|h\tilde b_n\|}
	{\|\mathscr L_f\|}_{\mathbb B}\\
	&={\|\nchi_D h\|}_{L^2(\mu)}{\|\mathscr L_f\|}_{\mathbb B}
	\leq{\|h\|}_{L^2(\mu)}{\|\mathscr L_f\|}_{\mathbb B}.
	\end{split}\end{equation}
	Now choose any sequence $(h_i)_i\subseteq{\rm LIP}_{bs}(\Y)$
	such that $h_i\to|Df|$ pointwise $\mu$-a.e.\ and (dominated) in $L^2(\mu)$.
	By writing \eqref{eq:Int_isom_claim4} with $h_i$ in place of $h$ and
	then letting $i\to\infty$, we conclude that
	${\big\||Df|\big\|}_{L^2(\mu)}\leq{\|\mathscr L_f\|}_{\mathbb B}$, as required.
\end{proof}
\section{Proof of the main result in the separable case}
In this section we assume that $(\Y,\sfd)$ is a complete and separable local $\Cat\kappa$ space  equipped with a Borel measure $\mu$ that is finite on bounded sets. As discussed in the introduction, the crucial step in the proof of Theorem \ref{main} is the construction of an embedding of the `abstract analytical object'  $\Dertwotwo$ into the `concrete and geometric bundle' $L^2(\T_G\Y;\mu)$ that preserves distances on fibres. The construction of such embedding is the scope of this section.
\medskip

We start by recalling the following general fact (see also \cite[Theorem~3.7]{SchioppaCurrents} for the general module homomorphism between derivations and $1$-currents; notice that it is obvious that the boundary operation and the divergence operator are in correspondence under this homomorphism).

\begin{lemma}[From derivations to currents]\label{le:dercurr}
Let $(\Y,\sfd,\mu)$ be a metric measure space. Fix any
derivation $b\in{\rm Der}^{1,1}(\Y;\mu)$.
Then the functional $T_b:\,\LIP_b(\Y)\times\LIP(\Y)\to\R$ defined by
\begin{equation}\label{eq:deftb}
T_b(g,f):=\int g b(f)\,\d\mu
\end{equation}
is a normal $1$-current and the mass measure $\|T_b\|$ satisfies
	\begin{equation}\label{ident0}
	\|T_b\|=|b|\mu.
	\end{equation}
\end{lemma}
See Remark \ref{derlip} for extending derivations to act on $\LIP(\Y)$.
\begin{proof}
	By Lemma \ref{ptnorm} we get the estimate
	\begin{equation}\label{a}
	\big|T_b(g,f)\big|\le \int|g||b|\lip_a f\,\d\mu\leq \Lip(f)\int |g||b|\,\d\mu
	\end{equation}
	and thus taking into account Lemma \ref{weakco} we see that $T_b$ is a finite mass 1-current, with 
	\begin{equation}
\label{eq:mass1}
\|T_b\|\le |b|\mu.
\end{equation}
It is moreover normal, since
\[
T_b(1,f)=\int b(f)\,\d\mu=-\int f \div(b)\,\d\mu,\qquad \forall f\in \LIP_b(\Y).
\]
We are left with proving \eqref{ident0}. By \eqref{eq:mass1}, it suffices to show that $\M(T_b)=\int |b|\,\d\mu$. Let $(x_n)\subset\Y$ be countable and dense, and let $f_n$ be the function $x\mapsto \sfd(x_n,x)$, for each $n\in\N$. For $\eps>0$ and $n\in\N$, set
\[
A_n:=\big\{x\in\Y:\ b(f_n)(x)\ge |b|(x)-\eps\big\}\quad \textrm{ and } B_n:=A_n\setminus\bigcup_{m<n}A_m.
\]
Since, by Proposition \ref{prop:bnorm}, $|b|=\sup_nb(f_n)$ $\mu$-almost everywhere, we have that the sets $A_n$ cover $\Y$ up to a set of $\mu$-measure zero. Thus the collection $(B_n)$ is a countable Borel partition of $\Y$ up to a $\mu$-null set. Let $B\subset\Y$ be a ball, and estimate
\begin{align*}
\int_B|b|\,\d\mu&=\sum_n\int_{B_n\cap B}|b|\,\d\mu\le\sum_n\int_{B_n\cap B}\big(b(f_n)+\eps\big)\,\d\mu=\eps\mu(B)+\sum_n\int_{B_n\cap B}b(f_n)\,\d\mu\\
&=\eps\mu(B)+\sum_nT_b(\nchi_{B_n\cap B},f_n).
\end{align*}
By the characterization of mass (cf.\ \cite[Proposition 2.7]{AmbrosioKirchheim00}), we obtain
\[
\int_B|b|\,\d\mu\le\eps\mu(B)+\|T_b\|(B).
\]
Since $\eps>0$ and $B$ are arbitrary, the claim follows.
\end{proof}

We now come to the construction of the embedding.
\begin{theorem}[Embedding of $\Dertwotwo$ into $L^2(\T_G\Y;\mu)$]\label{thm:main}
Let $(\Y,\sfd,\mu)$ be a complete and separable local $\Cat\kappa$ space 
equipped with a Borel measure $\mu$ which is finite on bounded sets, and let $b\in\Dertwotwo$. Then there exists a unique $v\in L^2(\T_G\Y;\mu)$ such that for any $\bar x\in\Y$ and $y\in B_{\sfr_{\bar x}}(\bar x)$ it holds that
\begin{equation}
\label{eq:claim}
\d_x\dist_y(v(x))=b(\dist_y)(x)\qquad \mu\text{-a.e.\ }x\in  B_{\sfr_{\bar x}}(\bar x).
\end{equation}
Moreover, $v$ satisfies
\begin{equation}
\label{eq:normv}
|v(x)|_x=|b|(x)\qquad\mu\text{-a.e.\ }x\in\Y.
\end{equation}
\end{theorem}
\begin{proof}\ \\
\noindent{\sc Borel regularity.}
Taking into account Proposition \ref{prop:diffdist}(i), we can rewrite \eqref{eq:claim} as
\begin{equation}
\label{eq:claim1}
\la v(x),(\G_x^y)'_0\ra_x=-\sfd(x,y)\,b(\dist_y)(x)\qquad \mu\text{-a.e.\ }x\in  B_{\sfr_{\bar x}}(\bar x).
\end{equation}
Thus taking into account the continuity of $y\mapsto (\G_x^y)'_0$, established in Theorem \ref{thm:samecone}, and the weak continuity of $y\mapsto b(\dist_y)$, given by Lemma \ref{weakco}, we see that \eqref{eq:claim} holds for every $y\in B_{\sfr_{\bar x}}(\bar x)$ if and only if it holds for a countable and dense set of $y\in B_{\sfr_{\bar x}}(\bar x)$. Since the continuity of $x\mapsto\sfr_x$ grants that $B_{\sfr_x}(x)\subset\cup_n B_{\sfr_{x_n}}(x_n)$ if $x_n\to x$, using an argument based on the Lindel\"of property of $\Y$, we can reduce the claim to checking \eqref{eq:claim} for a countable and dense set of $\bar x$'s.

Now for given $\bar x$, and $y\in B_{\sfr_{\bar x}}(\bar x)$ running in these countable sets, fix a Borel representative $f_{\bar x,y}$ of $b(\dist_y)$ on $B_{\sfr_{\bar x}}(\bar x)$ and notice that if $v$ satisfies \eqref{eq:claim} for any $y,\bar x$ in such countable sets, there is a Borel $\mu$-negligible set $\mathcal N\subset\Y$ such that $\d_x\dist_y(v(x))=f_{\bar x,y}(x)$ for every $x\in B_{\sfr_{\bar x}}(\bar x)\setminus\mathcal N$. Thus redefining $v$ on $\mathcal N$ by setting it to 0 and recalling Proposition \ref{prop:btgy} we conclude that any $v$ for which \eqref{eq:claim} holds for any $\bar x\in\Y$ and $y\in B_{\sfr_{\bar x}}(\bar x)$ is, up to modification in a negligible set, a Borel section of $\T_G\Y$.

\noindent{\sc Integrability.} Propositions \ref{prop:bnorm} and \ref{prop:diffdist} ensure that any $v$ for which \eqref{eq:claim} holds also satisfies \eqref{eq:norm}. This, together with the Borel measurability proved above, implies that any $v$ satisfying \eqref{eq:claim} belongs to $L^2(\T_G\Y;\mu)$.

\noindent{\sc Uniqueness.} Let $v_1,v_2\in L^2(\T_G\Y;\mu)$ satisfy \eqref{eq:claim} so that, by what we already proved, we have that $|v_1(x)|_x=|v_2(x)|_x$ for $\mu$-a.e.\ $x$. By Proposition \ref{prop:diffdist}$(iii)$, we conclude that $v_1(x)=v_2(x)$ for $\mu$-a.e.\ $x$.

\noindent{\sc Existence.} Assume at first that $b\in{\rm Der}^{1,1}(\Y;\mu)$ and let $T_b$ be defined as in Lemma \ref{le:dercurr}, so that $T_b$ is a normal 1-current. By Theorem \ref{thm:sup2}, we find a finite non-negative Borel measure $\pi$ on $C([0,1];\Y)$ concentrated on curves with constant speed for which \eqref{eq:sup} holds with $T=T_b$. Notice that, by restricting $\pi$ to the complement of the set of constant curves (this does not affect the validity of \eqref{eq:sup}), we can assume that $\pi$ gives 0 mass to constant curves. 

Let ${\rm e}:\,C([0,1];\Y)\times[0,1]\to \Y$ be the evaluation map defined as ${\rm e}(\gamma,t):=\gamma_t$ and put $\hat\pi:=\pi\times\mathcal L^1\restr{[0,1]}$ and $\nu:={\rm e}_*\hat\pi$. Since $C([0,1];\Y)\times[0,1]$ and $\Y$ are Polish spaces, we may apply the disintegration theorem (see e.g.\ \cite[Theorem 5.3.1]{AmbrosioGigli08} or \cite[Chapter 45]{Fremlin4}) to $\hat\pi$ and ${\rm e}$ to find a weakly measurable family $\{\hat\pi_x\}_{x\in\Y}$ of Borel probability measures on  $C([0,1];\Y)\times[0,1]$ such that 
\begin{equation}
\label{eq:distconc}
{\rm e}_*\hat\pi_x=\delta_x\qquad\text{ for }\nu\text{-a.e.\ }x
\end{equation}
and
\begin{equation}
\label{eq:disint}
\int \Psi(\gamma,t)\,\d\hat\pi(\gamma,t)=\int \Big(\int\Psi\,\d\hat\pi_x(\gamma,t)\Big)\,\d\nu(x)
\end{equation}
for any Borel real-valued map $\Psi$ for which any of these two integrals makes sense.

Recall that the map ${\sf RightDer}:\,C([0,1];\Y)\times[0,1]\to \T_G\Y$ defined in \eqref{eq:rightder} is Borel (Proposition \ref{prop:rightder}) and set
$\nn_x:={\sf RightDer}_*\hat\pi_x$. Notice that although, by definition, the measures $\nn_x$ are measures on $\T_G\Y$, in fact for $\nu$-a.e.\ $x$ we have that $\nn_x$ is concentrated on $\T_x\Y$ and will therefore be considered, with a slight abuse of notation, as a measure on $\T_x\Y$. To see this, let $\pi^\Y:\,\T_G\Y\to\Y$ be the canonical projection and notice that ${\rm e}=\pi^\Y\circ {\sf RightDer}$, thus \eqref{eq:distconc} gives $\pi^\Y_*\nn_x=\delta_x$ for $\nu$-a.e.\ $x$, which implies the claim.

Now observe that, for any $g\in\LIP_b(\Y)$, we have
\[
\int g|b|\,\d\mu\stackrel{\eqref{eq:sup},\eqref{ident0}}=\iint_0^1g(\gamma_t)|\dot\gamma_t|\,\d \hat\pi(\gamma,t).
\]
By Proposition \ref{prop:velo} and the definition of ${\sf Norm}:\T_G\Y\to\R^+$ given in Corollary \ref{cor:normborel} (which also grants that this map is Borel, so that the integrals below are well-defined) we have
\[
\iint_0^1g(\gamma_t)|\dot\gamma_t|\,\d \hat\pi(\gamma,t)=\iint_0^1g( {\rm e}(\gamma,t))\,{\sf Norm}({\sf RightDer}(\gamma,t))\,\d \hat\pi(\gamma,t).
\]
Therefore we have
\begin{equation}
\label{eq:formass}
\begin{split}
\int g|b|\,\d\mu&=\iint_0^1g( {\rm e}(\gamma,t))\,{\sf Norm}({\sf RightDer}(\gamma,t))\,\d \hat\pi(\gamma,t)
\\
\text{(by \eqref{eq:disint})}\qquad&=\int g(x)\int {\sf Norm}({\sf RightDer}(\gamma,t))\,\d \hat\pi_x(\gamma,t)\,\d\nu(x)\\
\text{(by definition of $\nn_x$)}\qquad&=\int g(x)\int {\sf Norm}(y,v)\,\d \nn_x(y,v)\,\d\nu(x)\\
&=\int g(x)\int |v|_x\,\d \nn_x(v)\,\d\nu(x).
\end{split}
\end{equation}
As mentioned above, in the last step we made the slight abuse of notation in considering $\nn_x$ as a measure on $\T_x\Y$. In particular, choosing $g\equiv 1$, we get
\[
\iint |v|_x\,\d \nn_x(v)\,\d\nu(x)=\int |b|\,\d\mu<\infty,
\]
which implies that $\nn_x\in \mathscr P_1(\T_x\Y)$ for $\nu$-a.e.\ $x$. Let us define
\[
{\sf B}(x):=\Ba(\nn_x)\in\T_x\Y\qquad\text{ for }\nu\text{-a.e.\ }x.
\]
Now set, for brevity,  $\Phi(x):=\int |v|\,\d \nn_x(v)$ and notice that the regularity granted by the disintegration theorem ensures that $\Phi$ is Borel. Also, the fact that $\pi$ is concentrated on curves whose speed is constant and non-zero tells that ${\sf Norm}({\sf RightDer}(\gamma,t))>0$ for $\hat\pi$-a.e.\ $(\gamma,t)$ and hence that $\Phi>0$ $\nu$-a.e.. Now notice that \eqref{eq:formass} and the arbitrariness of $g$ yield $|b|\mu=\Phi\nu$, so that the positivity of $\Phi$ implies $\nu\ll\mu$. Hence it holds that $|b|\mu=\Phi\frac{\d\nu}{\d\mu}\mu$, i.e.\
\begin{equation}
\label{eq:nb}
|b|(x)=\frac{\d\nu}{\d\mu}(x)\int |v|_x\,\d \nn_x(v)\qquad\mu\text{-a.e.\ }x. 
\end{equation}
Let $\bar x\in\Y$ and $\bar y\in B_{\sfr_{\bar x}}(\bar x)$ and denote $f:=\dist_{\bar y}$. Thus $f$ is Lipschitz and semiconvex on $B_{\sfr_{\bar x}}(\bar x)$. Then, for $g\in\LIP_b(\Y)$ with support in $B_{\sfr_{\bar x}}(\bar x)$, we have, by the same considerations as before to justify the computations (and writing $\d f(x,v)$ for $\d_xf(v)$):
\[
\begin{split}
\int gb(f)\,\d\mu&=\iint_0^1g( {\rm e}(\gamma,t))\,\d f({\sf RightDer}(\gamma,t))\,\d \hat\pi(\gamma,t)
\\
\text{(by \eqref{eq:disint})}\qquad&=\int g(x)\int \d f({\sf RightDer}(\gamma,t))\,\d \hat\pi_x(\gamma,t)\,\d\nu(x)\\
\text{(by definition of $\nn_x$)}\qquad&=\int g(x)\int \d f(y,v)\,\d \nn_x(y,v)\,\d\nu(x)\\
&=\int g(x)\int \d_xf(v)\,\d \nn_x(v)\,\d\nu(x).
\end{split}
\]
By the arbitrariness of $g$, it follows that
\begin{equation}
\label{eq:almostdone}
b(\dist_{\bar y})(x)=\frac{\d\nu}{\d\mu}(x)\int \d_x\dist_{\bar y}(v)\,\d \nn_x(v)\qquad\mu\text{-a.e.\ }x\in B_{\sfr_{\bar x}}(\bar x).
\end{equation}
By the Jensen inequality  recalled in Subsection \ref{se:bari} and the convexity and continuity of $\d_xf$ (Proposition \ref{lem:exdiff}) this gives
\begin{equation}
\label{eq:forinf}
b(\dist_{\bar y})(x)\geq \frac{\d\nu}{\d\mu}(x)\,\d_x\dist_{\bar y}({\sf B}(x))\qquad\mu\text{-a.e.\ }x\in B_{\sfr_{\bar x}}(\bar x).
\end{equation}
Now we let $\bar x$ vary in a countable set so that the balls $B_{\sfr_{\bar x}}(\bar x)$ cover the whole $\Y$ (such set can be found by the Lindel\"of property of $\Y$) and for each such $\bar x$ we let  $\bar y$ vary in a countable dense set in $B_{\sfr_{\bar x}}(\bar x)$: taking the infimum in \eqref{eq:forinf} among these $\bar x,\bar y$ and recalling Proposition \ref{prop:diffdist}$(ii)$ and Proposition \ref{prop:bnorm}, we deduce
\[
-|b|(x)\geq - \frac{\d\nu}{\d\mu}(x)|{\sf B}(x)|_x\qquad\mu\text{-a.e.\ }x\in\Y.
\]
Hence taking into account \eqref{eq:nb} we obtain
\[
|{\sf B}(x)|_x\geq \int |v|_x\,\d \nn_x(v)\qquad\nu\text{-a.e.\ }x\in\Y.
\]
Since $|v|_x=\sfd_x(v,0)$, by the rigidity statement in Proposition \ref{lem:rigidity}  we deduce that $\nn_x$ is concentrated on a half-line starting from $0\in\T_x\Y$ for $\nu$-a.e.\ $x$. For any $x$ for which this is true, it is easy to check (see also \cite[Example 5.2]{stu03}) that any positively 1-homogeneous function $h:\T_x\Y\to\R$ satisfies
\[
\int h(v)\,\d\nn_x(v)=h({\sf B}(x)).
\]
Applying this identity to $h:=\d_x\dist_{\bar y}$, from \eqref{eq:almostdone} we get
\[
b(\dist_{\bar y})(x)=\frac{\d\nu}{\d\mu}(x) \, \d_x\dist_{\bar y}({\sf B}(x))=\d_x\dist_{\bar y}\Big(\frac{\d\nu}{\d\mu}(x) {\sf B}(x)\Big)\qquad\mu\text{-a.e.\ }x\in B_{\sfr_{\bar x}}(\bar x),
\]
which by the arbitrariness of $\bar y$ means that $v(x):=\frac{\d\nu}{\d\mu}(x) {\sf B}(x)$ satisfies \eqref{eq:claim}, and thus concludes the proof for $b\in{\rm Der}^{1,1}(\Y;\mu)$.

For the case $b\in\Dertwotwo$ we argue as follows. Fix $\bar x\in\Y$ and let $(\eta_n)$ be a sequence of Lipschitz functions with bounded support such that $\eta_n\equiv 1$ on $B_n(\bar x)$. Then, by the Leibniz rule for the divergence (cf.\ \cite[Lemma 7.1.2]{DiM14a}), we see that $\eta_nb\in{\rm Der}^{1,1}(\Y;\mu)$. Thus we have the existence of $v_n\in L^2(\T_G\Y;\mu)$ satisfying \eqref{eq:claim} for $b_n$. In particular, by \eqref{eq:normv}, we  have that
\begin{equation}
\label{eq:normvn}
|v_n(x)|_x=|b_n|(x)=|b|(x)\qquad\mu\text{-a.e.\ }x\in B_n(\bar x).
\end{equation}
From the weak locality of derivations it follows that $v_n=v_m$ on $B_n(\bar x)$ for every $m\geq n$, hence the Borel section $v$ of $\T_G\Y$ given by
\[
v(x):=v_n(x)\qquad\text{ for }\mu\text{-a.e.\ }x\in B_n(\bar x),\ \forall n\in\N
\]
is well-defined and, by \eqref{eq:normvn} and the assumption $|b|\in L^2(\mu)$, belongs to $L^2(\T_G\Y;\mu)$. Then again the weak locality of derivations ensures that $v$ satisfies \eqref{eq:claim}, thus concluding the proof.
\end{proof}
\begin{definition}
For $b\in\Dertwotwo$ we shall denote by $\mathscr F(b)$ the section $v\in L^2(\T_G\Y;\mu)$ given by Theorem \ref{thm:main}. Thus we have a map
\[
\mathscr F:\,\Dertwotwo\to L^2(\T_G\Y;\mu),\quad b\mapsto \mathscr F(b).
\]
\end{definition}
Then we have:
\begin{corollary}[`Linearity' of $\mathscr F$]\label{cor:main}
Let $(\Y,\sfd,\mu)$ be a complete and separable local $\Cat\kappa$ space equipped with a Borel measure $\mu$ finite on bounded sets and $b_1,b_2\in\Dertwotwo$. Then $\mu$-a.e.\ we have
\begin{equation}
\label{eq:linearf}
\begin{split}
\mathscr F(b_1+b_2)&=\mathscr F(b_1)\oplus\mathscr F(b_2),\\
\sfd_\cdot(\mathscr F(b_1),\mathscr F(b_2))&=|b_1-b_2|,\\
|\mathscr F(b_1+b_2)|_\cdot^2+|\mathscr F(b_1-b_2)|_\cdot^2&=2\big(|\mathscr F(b_1)|_\cdot^2+|\mathscr F(b_2)|_\cdot^2\big).
\end{split}
\end{equation}
\end{corollary}
\begin{proof} The statement is local in nature, thus up to using a countable cover of $\Y$ with balls of the form $B_{\sfr_x}(x)$, we can assume that $\Y$ is a separable $\Cat\kappa$ space with diameter $<D_\kappa$.

Now let  $(y_n)\subset\Y$ be countable and dense and put for brevity $f_n:=\dist_{y_n}$. For every $n\in\N$ we have
\[
\begin{split}
|(b_1-b_2)(f_n)|&=|b_1(f_n)-b_2(f_n)|\stackrel{\eqref{eq:claim}}=|\d_\cdot f_n(\mathscr F(b_1))-\d_\cdot f_n(\mathscr F(b_2))|\leq \sfd_\cdot\big(\mathscr F(b_1),\mathscr F(b_2)\big)
\end{split}
\]
$\mu$-a.e., having used the fact that $\d_xf_n$ is 1-Lipschitz in the last step (Proposition \ref{lem:exdiff}). Passing to the supremum in $n$ we obtain
\begin{equation}
\label{eq:minus}
|b_1-b_2|\leq \sfd_\cdot\big(\mathscr F(b_1),\mathscr F(b_2)\big)\qquad\mu\text{-a.e..}
\end{equation}
On the other hand, using the convexity and positive 1-homogeneity of $\d_xf_n$ (Proposition \ref{lem:exdiff}) we have
\begin{equation}\label{4}\begin{split}
\d_x f_n\big(\mathscr F(b_1)(x)\oplus\mathscr F(b_2)(x)\big)
&\leq \d_x f_n\big(\mathscr F(b_1)(x)\big)+\d_x f_n\big(\mathscr F(b_2)(x)\big)
\\
&=b_1(f_n)(x)+b_2(f_n)(x)\\ 
&=(b_1+b_2)(f_n)(x)\\
&=\d_x f_n\big(\mathscr F(b_1+b_2)(x)\big)
\end{split}\end{equation}
for $\mu$-a.e.\ $x$. By Proposition \ref{prop:diffdist}$(iii)$ and the arbitrariness of $n$ this implies
\begin{equation}
\label{eq:plus}
\big|\mathscr F(b_1+b_2)\big|_\cdot\leq \big|\mathscr F(b_1)\oplus\mathscr F(b_2)\big|_\cdot
\qquad\mu\text{-a.e..}
\end{equation}
Therefore, $\mu$-a.e.\ we have
\begin{align*}
|b_1-b_2|^2+|b_1+b_2|^2&\leq \sfd^2_\cdot\big(\mathscr F(b_1),\mathscr F(b_2)\big)+\big|\mathscr F(b_1+b_2)\big|_\cdot^2&&\text{by \eqref{eq:minus},\eqref{eq:claim}}\\
&\leq \sfd^2_\cdot\big(\mathscr F(b_1),\mathscr F(b_2)\big)+ \big|\mathscr F(b_1)\oplus\mathscr F(b_2)\big|_\cdot^2&&\text{by \eqref{eq:plus}}\\
&\leq 2\,\big|\mathscr F(b_1)\big|_\cdot^2+2\,\big|\mathscr F(b_2)\big|_\cdot^2&&\text{by \eqref{eq:PI}}\\
&=2\,|b_1|^2+2\,|b_2|^2&&\text{by \eqref{eq:claim}}.
\end{align*}
Writing this for $b_1+b_2,b_1-b_2$ in place of $b_1,b_2$ we see that all the inequalities that we used are in fact equalities.

In particular the last inequality is an equality, thus proving the last identity in \eqref{eq:linearf}. The equality in \eqref{eq:minus} is the second in \eqref{eq:linearf}. Finally, the equality in \eqref{eq:plus} and Proposition \ref{prop:diffdist}$(iii)$ imply the first identity in \eqref{eq:linearf}. This completes the proof.
\end{proof}

We can now easily prove our main result. We restrict ourselves to the separable setting for the moment, and postpone the technical differences to deal with in non-separable spaces to the next section.

\bigskip{\bf Proof of Theorem \ref{main} for separable spaces.}
By Proposition \ref{prop:duality} we have \[\big\||Df|\big\|_{L^2(\mu)}=\|\mathscr L_f\|_{\B}=\sup\left\{ \int L_f(b)\,\d\mu\;:\;b\in\Dertwotwo,\,\|b\|_2\le 1 \right\}.  \] Since the space 
\begin{equation}
\label{eq:defD}
\D:=\big(\Dertwotwo,\|\cdot\|_2\big)
\end{equation}
is (pre)Hilbert by Theorem \ref{thm:main} and Corollary \ref{cor:main} (in particular by the third in \eqref{eq:linearf}), it follows that its dual is a Hilbert space (note that $\mathscr L_f\in \D^*=\B$ in the notation of Proposition \ref{prop:duality}). Thus
\[\begin{split}
{\big\||D(f+g)|\big\|}^2_{L^2(\mu)}+{\big\||D(f-g)|\big\|}^2_{L^2(\mu)}
&={\|\mathscr L_{f+g}\|}^2_{\mathbb B}
+{\|\mathscr L_{f-g}\|}^2_{\mathbb B}\\
&={\|\mathscr L_f+\mathscr L_g\|}^2_{\mathbb B}
+{\|\mathscr L_f-\mathscr L_g\|}^2_{\mathbb B}\\
&=2\,{\|\mathscr L_f\|}^2_{\mathbb B}
+2\,{\|\mathscr L_g\|}^2_{\mathbb B}\\
&=2\,{\big\||Df|\big\|}^2_{L^2(\mu)}
+2\,{\big\||Dg|\big\|}^2_{L^2(\mu)}.
\end{split}\]
This completes the proof.
\qed

\bigskip In fact, as we shall see shortly, the completion of the space $\D$ defined in \eqref{eq:defD} is isomorphic to the $L^2$-tangent module. This is the content of Proposition \ref{prop:isom-module} below. We briefly introduce some additional machinery before stating the proposition.

\bigskip Recall the space of $L^2$-derivations
\[
{\rm Der}^2(\Y;\mu)=\big\{ b\in {\rm Der}(\Y;\mu)\,:\,|b|\in L^2(\mu)\big\}
\]
which, by \cite[Section 7.1.1]{DiM14a}, is complete when equipped with the norm $\|\cdot\|_2$. Since $\D\subset {\rm Der}^2(\Y;\mu)$, the completion $\overline\D$ of $\D$ under $\|\cdot\|_2$ is a Banach space and satisfies $\overline\D\subset {\rm Der}^2(\Y;\mu)$. In particular, there is a pointwise norm $|\cdot|:\,\overline\D\to L^2(\mu)$ given by the norm of a derivation (see Lemma \ref{ptnorm}). Using the fact that $\D$ is a $\LIP_{bs}(\Y)$-module (cf.\ \cite[Lemma 7.1.2]{DiM14a}), Lemma \ref{lem:approx} and the dominated convergence theorem, we see that $\overline\D$ is an $L^\infty(\mu)$-module. Thus $\big(\overline\D,\|\cdot\|_2,|\cdot|\big)$ is an $L^2(\mu)$-normed $L^\infty(\mu)$-module. We refer to \cite{Gigli14} for the theory of normed $L^\infty(\mu)$-modules.

The estimate (\ref{eq:def_pwg}) implies that, given $f\in W^{1,2}(\Y,\sfd,\mu)$, the module-homomorphism $L_f:\Dertwotwo\to L^1(\mu)$ extends to a $L^\infty(\mu)$-linear bounded map $L_f:\overline\D\to L^1(\mu)$ satisfying the bound
\begin{equation}\label{mini-ineq}
\big|L_f(\bar b)\big|\le|Df||\bar b|,\quad \bar b\in \overline\D.
\end{equation}

We briefly recall that the cotangent module $L^2(\T^*\Y;\mu)$ (see \cite{Gigli14}) is an $L^2(\mu)$-normed $L^\infty(\mu)$-module, equipped with an exterior derivative
\[
\d:\,W^{1,2}(\Y,\sfd,\mu)\to L^2(\T^*\Y;\mu)
\]
whose image generates $L^2(\T^*\Y;\mu)$ as a module.  The \emph{tangent module $L^2(\T\Y;\mu)$} is defined to be the module dual of $L^2(\T^*\Y;\mu)$. A vector field $X\in L^2(\T\Y;\mu)$ is said to have Sobolev divergence if there exists a function $g\in L^2(\mu)$ such that 
\[
\int fg\,\d\mu=-\int\d f(X)\,\d\mu,\quad f\in W^{1,2}(\Y,\sfd,\mu).
\] 
The function $g$, if it exists, is unique, and denoted by $\dive_SX$. We denote by $D(\dive_S)$ the vector space of elements of $L^2(\T\Y;\mu)$ that have Sobolev divergence. See \cite{Gigli14} for the details.

\begin{proposition}\label{prop:isom-module}
Let $(\Y,\sfd,\mu)$ be an infinitesimally Hilbertian metric measure space. Then the map
\[
A:L^2(\T\Y;\mu)\to {\rm Der}(\Y;\mu),\quad X\mapsto X\circ\d\restr{\LIP_b(\Y)}
\]
takes values in $\overline\D$ and provides an isomorphism of modules between $L^2(\T\Y;\mu)$ and $\overline\D$.
\end{proposition}
\begin{proof}
It is easy to see that, if $X\in L^2(\T\Y;\mu)$ has divergence $\dive_SX\in L^2(\mu)$, then $A(X)$ has divergence in the sense of (\ref{div}), and
\[
\dive_SX=\dive A(X)
\]
$\mu$-almost everywhere. Since $W^{1,2}(\Y,\sfd,\mu)$ is a Hilbert space, \cite[Proposition 2.3.17]{Gigli14} implies that  $L^2(\T\Y;\mu)$ is a Hilbert module. As a simple consequence of \cite[Proposition 2.3.14 and (2.3.13)]{Gigli14}, the space $D(\dive_S)$ is dense in $L^2(\T\Y;\mu)$.
Thus since we already noticed that $A(D(\dive_S))\subset \D$, we also get that $A(L^2(\T\Y;\mu))\subset \overline\D$. We will prove that $A$ is a module isomorphism $L^2(\T\Y;\mu)\to \overline\D$.

For each $f\in \LIP_b(\Y)$, $g,h\in L^\infty(\mu)$ and $V,W\in L^2(\T\Y;\mu)$, we have 
\begin{align*}
A(gV+hW)(f)=(gV)(\d f)+(hW)(\d f)=gV(\d f)+hW(\d f)=\big(gA(V)+hA(W)\big)(f),
\end{align*}
establishing that $A$ is a $L^\infty(\mu)$-linear module homomorphism $L^2(\T\Y;\mu)\to \overline\D$. Note that
\[
\big|A(V)\big|={\rm ess\,sup}\big\{V(\d f):\ f\in \LIP_b(\Y),\ \Lip(f)\le 1\big\}\le |V|_\ast,
\]
so that $A$ is bounded. By definition, we have that
\begin{align*}
|V|_\ast={\rm ess\,sup}\left\{\sum_{i=1}^m\nchi_{E_i}\big|V(\d f_i)\big|\,:\,\sum_{i=1}^m\nchi_{E_i}|\d f_i|\le 1,\ f_1,\ldots,f_m\in W^{1,2}(\Y,\sfd,\mu)  \right\}.
\end{align*}
Since $W^{1,2}(\Y,\sfd,\mu)$ is a Hilbert space, using Theorem \ref{thm:relax} and Mazur's lemma, it is easy to see that  $\LIP_{bs}(\Y)$ is dense in $W^{1,2}(\Y,\sfd,\mu)$ (see also \cite[Corollary 2.9]{Gigli12}). From this and Theorem \ref{thm:relax}, it follows that
\begin{align*}
|V|_\ast={\rm ess\,sup}\left\{ \sum_{i=1}^m\nchi_{E_i}\big|V(\d f_i)\big|\,:\,\sum_{i=1}^m\nchi_{E_i}\lip_af_i\le 1,\ f_1,\ldots,f_m\in \LIP_{bs}(\Y)  \right\}.
\end{align*}
Thus we have 
\begin{align*}
|V|_\ast &={\rm ess\,sup}\left\{\sum_{i=1}^m\nchi_{E_i}\big|V(\d f_i)\big|\,:\,\sum_{i=1}^m\nchi_{E_i}\lip_af_i\le 1 \right\}\\
&={\rm ess\,sup}\left\{\sum_{i=1}^m\nchi_{E_i}\big|A(V)(f_i)\big|\,:\,\sum_{i=1}^m\nchi_{E_i}\lip_af_i\le 1 \right\}\\
&\le{\rm ess\,sup}\left\{\sum_{i=1}^m\nchi_{E_i}\big|A(V)\big|\lip_a(f_i)\,:\,\sum_{i=1}^m\nchi_{E_i}\lip_af_i\le 1\right\}=\big|A(V)\big|.
\end{align*}

We have established that $A:\,L^2(\T\Y;\mu)\to \overline\D$ is an $L^\infty(\mu)$-module homomorphism satisfying
\[
\big|A(V)\big|=|V|_\ast
\]
pointwise $\mu$-almost everywhere, for every $V\in L^2(\T\Y;\mu)$. To show it is an isometric module isomorphism, it suffices to prove that it is onto.

Let $\bar b\in\overline\D$. Define the linear map
\[
L:W^{1,2}(\Y,\sfd,\mu)\to L^1(\mu),\quad f\mapsto L_f(\bar b).
\]
By (\ref{mini-ineq}) and \cite[Proposition 1.4.8]{Gigli14}, $L$ extends to a vector field $X\in L^2(\T\Y;\mu)$ satisfying
\[X\circ\d\restr{W^{1,2}(\Y,\sfd,\mu)}=L.\]
In particular, for $f\in \LIP_{bs}(\Y)$, we have
\[
A(X)(f)=X(\d f)=L(f)=L_f(\bar b)=\bar b(f).
\]
This implies the surjectivity of $A$, and concludes the proof.
\end{proof}

See \cite{DiMarino14} for more on preduals of the Sobolev spaces.

\bigskip {\bf Proof of Theorem \ref{main2}.} Let $(\Y,\sfd)$ be a complete and separable $\Cat\kappa$-space, and $\mu$ a Borel measure on $\Y$, which is finite on bounded sets. By the proof above of Theorem \ref{main} in the separable case, we have that $W^{1,2}(\Y,\sfd,\mu)$ is a Hilbert space. From Theorem \ref{thm:main} and Corollary \ref{cor:main} it follows that the space $\overline\D$ admits an isometric embedding
\[
\mathscr F':\,\overline\D\to L^2(\T_G\Y;\mu)
\]
satisfying (\ref{eq:linearf}). Thus the claim follows directly from Proposition \ref{prop:isom-module} by precomposing $\mathscr F'$ with the isometric module isomorphism $A:L^2(\T\Y;\mu)\to \overline\D$.\qed

\section{The non-separable case}

In defining derivations and Sobolev functions we assumed, following \cite{DiM14a}, that the underlying metric space is separable. Yet, as noted in the introduction, from a purely geometric perspective it is quite unnatural to impose a separability condition when dealing with $\Cat\kappa$ spaces. In this section we discuss how to remove the condition of separability, the relevant result being Theorem \ref{thm:consistency}. Let us remark that we shall continue to assume that the measure $\mu$ on $\Y$ has separable support, or equivalently that it is tight: the discussion here concerns the definition of Sobolev functions itself, in this setting.

\bigskip

One of the reasons for the success of the theory of Sobolev calculus on metric measure spaces is that there are many different  definitions of Sobolev spaces in such environment which turn out to be equivalent. In trying to extend such an equivalence result to the non-separable setting one could either re-run all the arguments and check that they work even in the more general framework (this is possible -- and works -- but is quite tedious) or argue as below.

\bigskip

Out of the several definitions of Sobolev functions, there are two `extremal' ones introduced in \cite{AmbrosioGigliSavare11}: the one obtained by relaxation of the asymptotic Lipschitz constant (we shall denote the corresponding space and notion of minimal relaxed upper gradient by $W^{1,2}_{\rm rel}(\Y,\sfd,\mu)$ and $|\d f|_{\rm rel}$) and the one obtained by duality with test plans  (we shall denote the corresponding space and notion of minimal weak upper gradient by $W^{1,2}_{\rm tp}(\Y,\sfd,\mu)$ and $|\d f|_{\rm tp}$). These produce in some sense the `biggest' and `smallest' weak notion of upper gradient and it is easy to check from the definitions that 
\begin{equation}
\label{eq:ordersob}
W^{1,2}_{\rm rel}(\Y,\sfd,\mu)\subset W^{1,2}_{\rm tp}(\Y,\sfd,\mu)\quad\text{with}\quad |\d f|_{\rm tp}\leq |\d f|_{\rm rel}\quad \mu\text{-a.e.}\quad\forall f\in W^{1,2}_{\rm rel}(\Y,\sfd,\mu).
\end{equation}
One of the main results in \cite{AmbrosioGigliSavare11} is the proof that the two spaces and the two notions of upper gradients coincide. This fact is used by the first author in \cite{DiM14a} to prove that the notion of Sobolev space obtained by duality with derivations coincides with $W^{1,2}_{\rm rel}(\Y,\sfd,\mu)= W^{1,2}_{\rm tp}(\Y,\sfd,\mu)$ and induces the same upper gradient.

We add the following ingredient to the discussion above:
\begin{theorem}\label{thm:consistency}
Let $(\Y,\sfd,\mu)$ be a complete and separable metric space equipped with a positive Radon measure which is finite on bounded sets. Let $\Y_1,\Y_2\subset\Y$ be closed sets on which $\mu$ is concentrated. Set $\sfd_i:=\sfd\restr{\Y_i\times\Y_i}$, $\mu_i:=\mu\restr{\Y_i}$, $i=1,2$ and notice that the identity on the support of $\mu$ induces an isomorphism $\iota:L^2(\Y_1,\mu_1)\to L^2(\Y_2,\mu_2)$. Then:
\begin{itemize}
\item[i)]  $\iota$ induces an isomorphism from $W^{1,2}_{\rm rel}(\Y_1,\sfd_1,\mu_1)$ to $W^{1,2}_{\rm rel}(\Y_2,\sfd_2,\mu_2)$ which respects $|\d\cdot|_{\rm rel}$,
\item[ii)] $\iota$ induces an isomorphism from $W^{1,2}_{\rm tp}(\Y_1,\sfd_1,\mu_1)$ to $W^{1,2}_{\rm tp}(\Y_2,\sfd_2,\mu_2)$ which respects $|\d\cdot|_{\rm tp}$.
\end{itemize}
\end{theorem}
\begin{proof}
We can assume $\Y_2=\Y$.

\noindent{\bf (i)} Given that for any $f:\Y\to\R$ we have $\lip_a(f)(x)\geq \lip_a(f\restr{\Y_1})(x)$, we see that  $W^{1,2}_{\rm rel}(\Y,\sfd,\mu)\subset  W^{1,2}_{\rm rel}(\Y_1,\sfd_1,\mu_1)$ with $|\d f|_{{\rm rel},\Y_1}\leq|\d f|_{{\rm rel},\Y}$ for any $f\in W^{1,2}_{\rm rel}(\Y,\sfd,\mu)$. To prove the other inclusion and inequality, by the definition of $W^{1,2}_{\rm rel}(\Y,\sfd,\mu)$ it is sufficient to prove that for any Lipschitz function $f:\Y\to\R$ we have 
\[
|\d f|_{{\rm rel},\Y}\leq \lip_a(f\restr{\Y_1})\qquad\mu\text{-a.e..}
\] 
Fix a Lipschitz function $f:\Y\to\R$ and $\eps>0$. For any $x\in\Y_1$, let $r>0$ be such that $\Lip(f\restr{\Y_1\cap B_r(x)})\leq \lip_a(f\restr{\Y_1})+\eps$. By the McShane extension lemma there is a Lipschitz function $g:\Y\to\R$ coinciding with $f$ on $\Y_1\cap B_r(x)$ such that $\Lip(g)=\Lip(f\restr{\Y_1\cap B_r(x)}) $. By the locality property of relaxed upper gradients we see that 
\[
|\d f|_{{\rm rel},\Y}=|\d g|_{{\rm rel},\Y}\qquad\mu\text{-a.e.\ on }\{f=g\}\supset \Y_1\cap B_r(x).
\]
Keeping in mind that $|\d g|_{{\rm rel},\Y}\leq\Lip(g) $ and the construction we deduce that
\begin{equation}
\label{eq:pluseps}
|\d f|_{{\rm rel},\Y}\leq \lip_a(f\restr{\Y_1})+\eps\qquad\mu\text{-a.e.\ on }\Y_1\cap B_r(x).
\end{equation}
Repeat this argument for every $x\in\Y_1$ and then use the Lindel\"of property of $\Y_1$ to deduce that, as $x$ varies in a countable dense set, the balls $B_r(x)$ as above cover the whole $\Y_1$. Then \eqref{eq:pluseps} gives
\[
|\d f|_{{\rm rel},\Y}\leq \lip_a(f\restr{\Y_1})+\eps\qquad\mu\text{-a.e.}
\]
and the conclusion follows by letting $\eps\downarrow0$.

\noindent{\bf (ii)} It is sufficient to check that a  test plan on $\Y$ is also a test plan on $\Y_1$ and vice versa. The `vice versa' is obvious by the inclusion $C([0,1];\Y_1)\subset C([0,1];\Y)$. For the other implication it is sufficient to show that any test plan $\pi$ on $\Y$ is concentrated on $C([0,1];\Y_1)$. To see this, let ${\rm e}_t:\,C([0,1];\Y)\to \Y$ be defined by ${\rm e}_t(\gamma):=\gamma_t$ and notice that for any dense set $(t_n)\subset[0,1]$ the inclusion
\[
C([0,1];\Y)\setminus C([0,1];\Y_1)=\bigcup_n{\rm e}_{t_n}^{-1}(\Y\setminus\Y_1)
\]
holds. Since $({\rm e}_t)_*\pi\ll\mu$ and $\mu$ are concentrated on $\Y_1$, we have that $\pi\big({\rm e}_{t_n}^{-1}(\Y\setminus\Y_1)\big)=0$ for every $n$. The claim follows.
\end{proof}
Thanks to this result we can now give the following definition:
\begin{definition}[Sobolev spaces on non-separable metric spaces]\label{def:sobns}
Let $(\Y,\sfd,\mu)$ be a complete, not necessarily separable, metric space equipped with a non-negative and non-zero Radon measure $\mu$ giving finite mass to bounded sets.

Then the Sobolev space $W^{1,2}(\Y,\sfd,\mu)$ (and the corresponding notion of upper gradient $|\d f|$) is defined as $W^{1,2}(\Y_1,\sfd_1,\mu_1)$, where $\Y_1$ is any closed and separable subspace of $\Y$ on which $\mu$ is concentrated, while $\sfd_1:=\sfd\restr{\Y_1\times\Y_1}$ and $\mu_1:=\mu\restr{\Y_1}$.
\end{definition}
The role of Theorem \ref{thm:consistency} is to prove that this definition is consistent with the case of separable spaces. By the fact that most of the notions of Sobolev spaces in mm-spaces (including those of Cheeger \cite{Cheeger00}, \cite{Shanmugalingam00} and the first author \cite{DiM14a}) are naturally `chained' between $W^{1,2}_{\rm rel}$ and $W^{1,2}_{\rm tp}$ and, since these latter spaces coincide as already remarked, we see that Theorem \ref{thm:consistency} implies that all these notions remain unchanged when passing from $\Y_1$ to $\Y_2$, as in Theorem \ref{thm:consistency}. This is why we do not specify the definition of Sobolev space we are referring to in Definition \ref{def:sobns}: they all agree.

\bigskip

With this said, the proof of our main Theorem \ref{main} in the general case is a trivial consequence of the result established in the separable setting:

{\bf Proof of Theorem \ref{main} in the general non-separable setting.}
We need to prove that for any $f,g\in W^{1,2}(\Y,\sfd,\mu)$ it holds that
\begin{equation}
\label{eq:infhi}
|\d(f+g)|^2+|\d(f-g)|^2=2\big(|\d f|^2+|\d g|^2\big)\qquad\mu\text{-a.e..}
\end{equation}
Notice that the measure $\mu$ is by assumption finite on bounded sets and Radon. Hence it is concentrated on a countable union $Z$ of compact sets, which is separable.  Fix $x\in Z$. We claim that there exists $\Omega\subset\Y$ with the following properties:
\begin{align}
&\mu(\Omega)>0,\label{7.4}\\
&\bar\Omega\text{ is a separable $\Cat\kappa$ space,}\label{7.5}\\\
&\Omega\text{ contains a neighbourhood of $x$ in $\bar Z$},\label{7.6}\\
&\text{$\Omega$ is open in the space $\bar\Omega\cup\bar Z$ and in such space has  $\mu$-negligible boundary.}\label{7.7}
\end{align}
To construct such a set $\Omega$ we start by noticing that the map $r\mapsto\mu(B_r(x))$ is non-decreasing, hence continuous except at a countable number of points. Fix a continuity point $r<\sfr_x$, for which $\mu(B_r(x))>0$. Since $r$ is a continuity point, we have $\mu(\partial B_r(x))=0$. Let $C$ be the closed convex hull of $B_r(x)\cap \bar Z$ and define $\Omega$ as the interior of $C$ in $C\cup\bar Z$. (Notice that $\Omega\subset \bar\Omega\subset C$ and that, by convexity of the ball $B_r(x)$, $C\cap \bar Z=B_r(x)\cap \bar Z$.)

Since $\Omega$ is the interior of a convex set it follows that $\Omega$, and thus its closure $\bar\Omega$, is a $\Cat\kappa$-space. The set $\bar\Omega$ is separable by construction. This establishes \eqref{7.5}.

Note that $B_r(x)\cap \bar Z$ is open in $\bar Z$. Moreover, $B_r(x)\cap\bar Z\subset \Omega$. To see this, let $y\in B_r(x)\cap \bar Z$ and let $\eps>0$ be a radius for which $B_\eps(y)\subset B_r(x)$. Then
$$
B_\eps(y)\cap(C\cup\bar Z)=(B_\eps(y)\cap C)\cup(B_\eps(y)\cap\bar Z)=B_\eps(y)\cap C\subset C
$$
is a neighbourhood of $y$ in $C\cup\bar Z$. Thus $y$ is an interior point of $C$. This proves (\ref{7.4}) and (\ref{7.6}).

To show (\ref{7.7}), note that since $\Omega$ is open in $C\cup\bar Z$, it is open in $\bar\Omega\cup\bar Z$. It suffices to show that $\mu(\partial_{C\cup\bar Z}\Omega)=0$. This follows from the estimate
\begin{align*}
\mu(\partial_{C\cup\bar Z}\Omega)=&\mu(\partial_{C\cup\bar Z}\Omega\cap \bar Z)=\mu(\partial_{\bar Z}\Omega)\le \mu(\partial_{\bar Z}C)\\
=&\mu(\partial_{\bar Z}(C\cap \bar Z))=\mu(\partial_{\bar Z}(B_r(x)\cap Z))\le \mu(\partial B_r(x))=0.
\end{align*}
Thus we have constructed a set $\Omega$ with the desired properties.

By \cite[Theorem 4.19$(i)$]{AmbrosioGigliSavare11-2} applied with $\X:=\bar\Omega\cup\bar Z$ we see that $f\restr{\bar \Omega}\in W^{1,2}(\bar\Omega)$ with 
\begin{equation}
\label{eq:samegrad}
\big|\d(f\restr{\bar\Omega})\big|_{\bar\Omega}=|\d f|\qquad\mu\text{-a.e.\ on }\Omega,
\end{equation}
and the same holds for $g$. Since we know that Theorem \ref{main} holds on separable $\Cat\kappa$ spaces we have (see, e.g., also \cite[Proposition 2.3.17]{Gigli14}) that
\[
\big|\d\big((f+g)\restr{\bar\Omega}\big)\big|_{\bar\Omega}^2+\big|\d\big((f-g)\restr{\bar\Omega}\big)\big|_{\bar\Omega}^2
=2\big(\big|\d(f\restr{\bar\Omega})\big|_{\bar\Omega}^2+
\big|\d(g\restr{\bar\Omega})\big|_{\bar\Omega}^2\big)\qquad\mu\text{-a.e.\ on }\Omega.
\]
Then the conclusion \eqref{eq:infhi} comes from this identity, \eqref{eq:samegrad} and the Lindel\"of property of $Z$.
\qed

\def\cprime{$'$} \def\cprime{$'$}


\end{document}